\documentclass[11pt,a4paper]{article}
\usepackage[T1]{fontenc}
\usepackage[sc]{mathpazo}
\usepackage{microtype}
\usepackage{titlesec,fancyhdr}
\titleformat{\section}{\large\bfseries}{\thesection}{.7em}{}
\titleformat{\subsection}{\normalsize\bfseries}{\thesubsection}{.7em}{}
\titlespacing*{\section}{0pt}{3.0ex plus 1ex minus .2ex}{1.4ex plus .2ex}
\titlespacing*{\subsection}{0pt}{2.5ex plus .8ex minus .2ex}{1ex plus .2ex}
\usepackage[a4paper,margin=1in]{geometry}
\usepackage{amsmath,amssymb,amsthm,mathtools}
\usepackage{tikz,booktabs,array,needspace,placeins}
\usepackage{adjustbox,caption}
\usetikzlibrary{arrows.meta,calc}
\definecolor{inkblue}{HTML}{253B53}
\definecolor{cutteal}{HTML}{007F82}
\definecolor{warmrust}{HTML}{BB542D}
\definecolor{lightgray}{HTML}{E9EDF0}
\tikzset{every picture/.style={line cap=round,line join=round},
axis/.style={-{Stealth[length=4pt]},gray!65,thin,text=inkblue}, boundary/.style={inkblue,line width=.8pt}, cut/.style={cutteal,line width=1.2pt}, figure label/.style={font=\small,text=inkblue}}
\usepackage[numbers,sort&compress]{natbib}
\usepackage[hidelinks]{hyperref}
\hypersetup{pdftitle={Split Cuts in Fixed Integer Dimensions},
pdfauthor={Avinash Bhardwaj and Vishnu Narayanan},pdfkeywords={split closure, cutting planes, integer programming, polyhedral complexity}}
\newtheorem{theorem}{Theorem}[section]
\newtheorem{lemma}[theorem]{Lemma}
\newtheorem{proposition}[theorem]{Proposition}
\newtheorem{corollary}[theorem]{Corollary}
\DeclareMathOperator{\conv}{conv}
\DeclareMathOperator{\SC}{SC}
\DeclareMathOperator{\diag}{diag}
\DeclareMathOperator{\cone}{cone}
\DeclareMathOperator{\xc}{xc}
\newcommand{\R}{\mathbb R}
\newcommand{\Q}{\mathbb Q}
\newcommand{\Z}{\mathbb Z}
\newcommand{\eps}{\epsilon}
\numberwithin{equation}{section}

\title{\LARGE\bfseries Split Cuts in Fixed Integer Dimensions}
\author{Avinash Bhardwaj\qquad Vishnu Narayanan\\[.7em]
\small Department of Industrial Engineering and Operations Research\\
\small Indian Institute of Technology Bombay, Mumbai, India\\[.3em]
\small\texttt{abhardwaj@iitb.ac.in}\qquad\texttt{vishnu@iitb.ac.in}}
\date{}

\begin{document}
\maketitle
\thispagestyle{fancy}
\begin{abstract}
We study how many split cuts are needed to reach the integer hull, and how large split closures can be. In every fixed dimension, a rational pure-integer polyhedron admits a polynomial-time constructible split proof of its integer hull with polynomial total encoding length; consequently, its full split rank is polynomial in the input length. This establishes the split-rank conjecture of Basu, Conforti, Di Summa, and Jiang for pure-integer polyhedra in fixed dimension. Central to our approach is a preservation result: in arbitrary pure-integer dimension, an inequality whose violating points map, under an affine map with integer coefficients, into the interior of a planar lattice-free polyhedron has a polynomial split proof retaining every point satisfying it. Whereas Chv\'atal--Gomory closures have polynomial-size inequality descriptions in fixed dimension, split closures need not: tetrahedra of input length $\Theta(\log N)$ can have split closures with at least $2N-1$ facets, yet extension complexity $\Theta(\log N)$. For the system $x=f+Rs$, $s\ge0$, $x\in\mathbb Z^2$ of input length $L$, the split closure has an exact extended formulation with $O(L)$ inequalities and auxiliary variables; this bound is tight. Separation over the split closure is polynomial-time in fixed total dimension. Finally, on mixed-integer feasible rational polytopes with fixed integer dimension $p$ and arbitrary continuous dimension, we optimize a primitive integral objective on the integer variables using $O_p(G+1)$ split cuts. Each cut removes a current optimal LP vertex, and the algorithm returns an optimal LP vertex that is mixed-integer feasible. $G$ is the integrality gap of the LP relaxation, rounded down.
\end{abstract}
\begin{center}\small
\end{center}
\clearpage
\begingroup
\small
\tableofcontents
\endgroup
\clearpage

\section{Introduction}
\suppressfloats[t]

Split cuts are among the most useful cutting planes in integer programming. The Gomory mixed-integer and mixed-integer rounding cuts that modern solvers generate routinely are split cuts, and much of the theory of cutting planes is, in one way or another, a study of what split cuts can and cannot do. The idea behind them is simple. In a mixed-integer linear program, a linear function with integer coefficients involving only the integer variables takes an integer value at every feasible point, and so it never lies strictly between two consecutive integers. For such a function $\ell$ and an integer $k$, every feasible point therefore satisfies the \emph{split disjunction} $\ell(x)\le k$ or $\ell(x)\ge k+1$. Taking the convex hull of the points of the linear programming (LP) relaxation that satisfy either inequality can remove part of the relaxation without losing a single feasible point. Every linear inequality valid for this hull is a \emph{split cut}, and the \emph{split closure} is the intersection of these hulls over all split disjunctions, each applied to the same relaxation.

Two questions about split cuts motivate this paper. The first concerns the split closure itself: when only a few variables are integer, how many linear inequalities are needed to describe it, and can far fewer suffice if we allow additional continuous variables and then project? The second concerns sequences of cuts. Here the cuts are added one at a time, each derived from a split disjunction of the current relaxation, and we ask how many are needed to reach the integer hull, the convex hull of the integer feasible points. For optimization one may ask for more, namely that every cut remove a current optimal LP vertex.

Our main result answers the second question for pure-integer problems in fixed dimension. Basu, Conforti, Di Summa, and Jiang~\cite{BCDJComplexity} conjectured that the split rank of a rational pure-integer polyhedron is polynomial in its input length, and noted that even dimension three was open. We show that, in every fixed dimension, the integer hull of a rational polyhedron is reached by a sequence of split cuts that can be constructed in polynomial time and has polynomial total encoding length (Theorem~\ref{thm:fixedproof}); polynomially many rounds of the split closure therefore suffice. Central to our approach is a preservation property: certain inequalities can be derived by split cuts without cutting off a single point, integer or fractional, that already satisfies them (Theorem~\ref{thm:planarsimulation}). This preservation lets us apply induction on lower-dimensional sections while retaining the points contributed by the other sections.

The answer to the first question is rather different. Whereas Chv\'atal--Gomory closures have polynomial-size descriptions in fixed dimension~\cite{BE}, split closures need not: in dimension three we construct tetrahedra whose split closures have exponentially many facets in the input length, although their extension complexity is only linear in it (Theorem~\ref{thm:contrast}). For two-row systems the split closure always has an extended formulation of size linear in the input length (Theorem~\ref{thm:uniformcorner}). We also give polynomial-time separation over the split closure when the total dimension is fixed, and a strong split-cut algorithm whose number of cuts depends only on the number of integer variables and the objective gap (Theorems~\ref{thm:separation} and~\ref{thm:gapalgorithm}).

\subsection{Split cuts and notation}
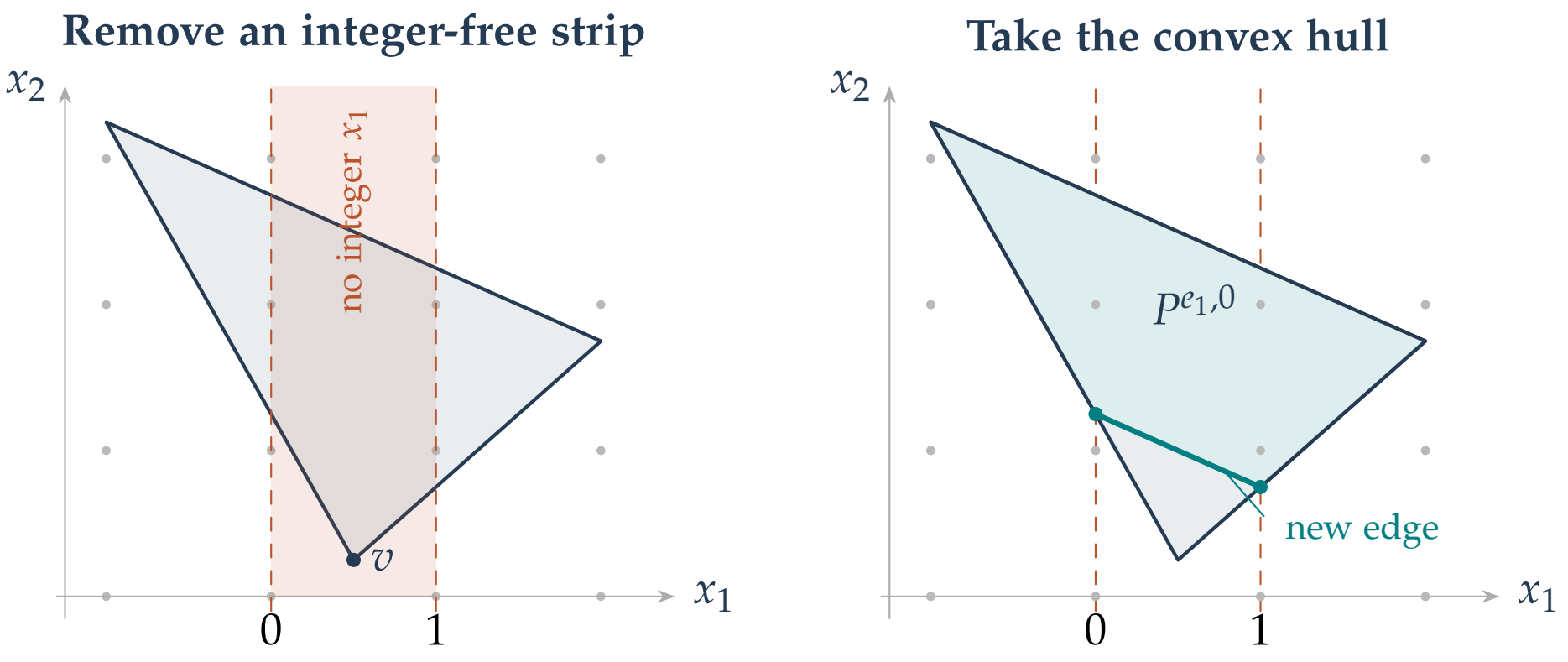
\begin{figure}[htbp]
\centering
\begin{adjustbox}{max width=\linewidth}
\begin{tikzpicture}[x=1.3cm,y=1.15cm,font=\small]
\foreach \shift/\ttl in {0/{Remove an integer-free strip},6.5/{Take the convex hull}} {
  \begin{scope}[xshift=\shift cm]
    \node[font=\small\bfseries,text=inkblue] at (.5,3.85) {\ttl};
    \fill[lightgray] (-1,3.25)--(.5,.25)--(2,1.75)--cycle;
    \foreach \x in {-1,0,1,2} {\foreach \y in {0,1,2,3} {\fill[gray!55] (\x,\y) circle(1pt);}}
    \draw[axis] (-1.3,0)--(2.45,0) node[right] {$x_1$};
    \draw[axis] (-1.25,-.15)--(-1.25,3.5) node[left] {$x_2$};
    \draw[boundary] (-1,3.25)--(.5,.25)--(2,1.75)--cycle;
    \draw[warmrust,dashed] (0,-.1)--(0,3.5);
    \draw[warmrust,dashed] (1,-.1)--(1,3.5);
    \node[below] at (0,0) {$0$}; \node[below] at (1,0) {$1$};
  \end{scope}}
\fill[warmrust,opacity=.12] (0,0)--(1,0)--(1,3.5)--(0,3.5)--cycle;
\node[text=warmrust,font=\scriptsize,rotate=90] at (.5,2.65) {no integer $x_1$};
\fill[inkblue] (.5,.25) circle(1.6pt) node[right] {$v$};
\begin{scope}[xshift=6.5cm]
  \fill[cutteal!13] (-1,3.25)--(0,1.25)--(1,.75)--(2,1.75)--cycle;
  \draw[boundary] (-1,3.25)--(0,1.25)--(1,.75)--(2,1.75)--cycle;
  \draw[cut] (0,1.25)--(1,.75);
  \fill[cutteal] (0,1.25) circle(1.6pt) (1,.75) circle(1.6pt);
  \foreach \x in {-1,0,1,2} {\foreach \y in {0,1,2,3} {\fill[gray!55] (\x,\y) circle(1pt);}}
  \node[figure label] at (.6,2) {$P^{e_1,0}$};
  \node[text=cutteal,anchor=west,font=\scriptsize] at (1.05,.45) {new edge};
  \draw[cutteal,thin] (1.02,.55)--(.78,.86);
\end{scope}
\end{tikzpicture}
\end{adjustbox}
\caption{Removing an integer-free strip and taking the convex hull creates a new edge, but leaves fractional points inside the strip. The dots mark the integer lattice.}
\label{fig:split}
\end{figure}

Let $P\subseteq\R^d$ be a rational polyhedron, and let $I\subseteq\{1,\ldots,d\}$ index the coordinates designated integer; the mixed-integer feasible points of $P$ are those $x\in P$ with $x_i\in\Z$ for every $i\in I$. Put $\Z_I^d=\{\pi\in\Z^d:\pi_j=0\text{ for }j\notin I\}$. For $\pi\in\Z_I^d\setminus\{0\}$ and $k\in\Z$, we call $\{x\in\R^d:k<\pi^Tx<k+1\}$ the \emph{open split strip} and $\{x\in\R^d:k\le\pi^Tx\le k+1\}$ the \emph{closed split strip}. No mixed-integer feasible point lies in the open strip, so each satisfies $\pi^Tx\le k$ or $\pi^Tx\ge k+1$, and the convex hull of the points of $P$ satisfying at least one of these branches is
\[
 P^{\pi,k}=\conv\bigl((P\cap\{\pi^Tx\le k\})\cup(P\cap\{\pi^Tx\ge k+1\})\bigr).
\]
A \emph{split cut} is any linear inequality valid for this hull, and the split closure is
\begin{equation}\label{eq:definition}
 \SC_I(P)=\bigcap_{\pi\in\Z_I^d\setminus\{0\},\ k\in\Z}P^{\pi,k}.
\end{equation}
We write $\SC(P)$ when all variables are integer, and set $\SC_\varnothing(P)=P$. For any polyhedron $Q\subseteq\R^d$, we write $Q_I=\conv(Q\cap\Z^d)$ for its integer hull. In the pure-integer case a single round of the split closure need not reach $P_I$, since convexification can put fractional points back inside the removed strip, as Figure~\ref{fig:split} shows.

Throughout, input length means the total binary encoding length of the supplied data, and a polyhedron is supplied by rational linear inequalities unless another representation is specified; a positive integer $N$ takes $O(\log(N+1))$ bits to write. We keep counts of inequalities and variables separate from total encoding length, which also accounts for their coefficients.

A \emph{polytope} is a bounded polyhedron, and we write $\R_+=[0,\infty)$ and $\Z_+=\Z\cap\R_+$. An integer vector is \emph{primitive} if the greatest common divisor (gcd) of its coordinates is one. An affine map with integer coefficients, including its constant terms, maps integer points to integer points. A unimodular matrix is an integer matrix of determinant $\pm1$, so that the corresponding change of coordinates preserves the integer lattice. Finally, $\{a\}=a-\lfloor a\rfloor$ denotes the fractional part of $a\in\R$, and subscripts on $O(\cdot)$ and $\operatorname{poly}(\cdot)$ indicate the parameters on which the constants or the polynomial may depend.

\subsection{Split-closure formulations}
Cook, Kannan, and Schrijver~\cite{CKS} proved that split closures are rational polyhedra, and Vielma~\cite{V07} gave a constructive finite description. For Chv\'atal--Gomory closures much more is known: Bockmayr and Eisenbrand~\cite{BE} showed that, in fixed dimension, they have inequality descriptions of size polynomial in the input length. Whether the same holds for split closures is a question that Dash~\cite{DashSurvey} attributes to Eisenbrand, and Basu, Conforti, Di Summa, and Jiang~\cite{BCDJ} answered it positively in the plane.

In dimension three the answer is negative, as our tetrahedral family shows. Perhaps surprisingly, the split closures in this family nevertheless have small extended formulations. Recall that a \emph{linear extended formulation} of a polyhedron is a linear system in the original variables and auxiliary continuous variables whose projection onto the original variables is that polyhedron. Its size is the number of inequalities, equations not being counted, and the minimum size is the extension complexity $\xc$~\cite{Y}.

\begin{theorem}\label{thm:contrast}
For every integer $N\ge2$, there is a rational tetrahedron $P_N\subseteq\R^3$ with an inequality description of total encoding length $\Theta(\log N)$ whose split closures $\SC(P_N)$ and $\SC_{\{1,2\}}(P_N)$ each have at least $2N-1$ facets. Both closures have extension complexity $\Theta(\log N)$ and admit exact rational extended formulations with $O(\log N)$ inequalities and auxiliary variables, constructible in time polynomial in $\log N$ with polynomial total encoding length.
\end{theorem}

The idea behind Theorem~\ref{thm:contrast} is geometric. The supporting hyperplanes of these facets meet an affine plane in tangents to a parabola, and all the inequalities needed for the closures fall into a constant number of families whose coefficients are quadratic functions of an integer parameter. Bogomolov, Fiorini, Maksimenko, and Pashkovich~\cite{BFMP} gave small extended formulations for polygons whose vertices lie on a parabola, and these yield a small formulation for each family. The mixed-integer example also answers, negatively, a question discussed by Kalantarzadeh~\cite[Section~1.3]{K}: whether split closures with two integer variables and one continuous variable always have polynomial-size inequality descriptions in the original variables.

A linear bound on formulation size also holds for the following system of two affine equations, with two coordinates designated integer and any number of nonnegative continuous variables.

\begin{theorem}\label{thm:uniformcorner}
Let $n\ge0$ be arbitrary, let $f\in\Q^2$ and $R\in\Q^{2\times n}$, and let $L$ be the total binary encoding length of $(f,R)$. Define
\[
Q(f,R)=\{(x,s)\in\R^2\times\R_+^n:x=f+Rs\}.
\]
With only $x_1,x_2$ designated integer, one can construct an exact rational linear extended formulation of the split closure of $Q(f,R)$ with $O(L)$ inequalities and auxiliary variables, uniformly in $n$. The construction time and total encoding length are polynomial in $L$.
\end{theorem}

This bound cannot be improved in general: already for $n=3$, there are examples requiring $\Omega(L)$ inequalities in every extended formulation of their split closures. The proof first selects split disjunctions sufficient to define the closure and groups their normal vectors into arithmetic progressions. Within each progression, either the cuts at the endpoints suffice, or one additional split cut makes redundant all cuts sufficiently far from the endpoints. The remaining intervals admit the same parabolic formulations as in Theorem~\ref{thm:contrast}, and bounding the sum of their sizes yields the linear bound.

To pass from $Q(f,R)$ to a general inequality system, we use a result of Andersen, Cornu\'ejols, and Li~\cite{ACL}, by which the split closure is the intersection of the split closures of the relaxations obtained by retaining maximal linearly independent sets of constraint rows. Applying our formulation to each of these relaxations yields the following bound.

\begin{theorem}\label{thm:twointeger}\label{thm:rowexcess}
Let $P=\{x\in\R^d:Ax\ge b\}$ be a rational polyhedron, where $A$ has $m$ rows and rank $r$. Let $I\subseteq\{1,\ldots,d\}$ index the coordinates designated integer, with $|I|\le2$, and let $L$ be the total input length of $(A,b,I)$. There is an absolute polynomial $p$ such that an exact rational linear extended formulation of $\SC_I(P)$ is constructible with total encoding length and running time at most $\binom mr p(L)$.
\end{theorem}

The construction is therefore polynomial-time whenever either $d$ or the difference $m-r$ is fixed. Since a full-dimensional simplex has a facet description with $m=d+1$ and $r=d$, simplices are covered even when $d$ varies.

The restriction to two integer variables can also be relaxed for a polyhedron defined by inequalities with linearly independent normal vectors. Let $V$ be its lineality space, the set of directions along which translation leaves the polyhedron unchanged, and consider the lattice of integer vectors that are supported on the integer coordinates and orthogonal to $V$. If this lattice has rank at most two, then the split closure has an extended formulation constructible in polynomial time with polynomial total encoding length. This allows any number of integer variables and does not fix the ambient dimension.

\subsection{Split proofs}\label{sec:introproofs}
A \emph{split proof} starts from the input relaxation and adds one split cut at a time, each justified by a specified split disjunction of the current relaxation; earlier cuts are never discarded. A split proof of an inequality ends with a relaxation on which that inequality is valid. Its length is the number of cuts, and its encoding includes the cut inequalities together with their split disjunctions. Nothing in this definition requires a cut to remove an optimal LP vertex. When we do impose this, we call a sequence of split cuts, or a split proof, \emph{strong} for a prescribed objective: every cut removes a vertex that is optimal for that objective over the current relaxation, the algorithm being free to choose among tied optimal vertices. The adjective \emph{weak}, used in contrast with strong, indicates that the requirement is not imposed. The distinction is the one between cutting-plane proofs and cutting-plane algorithms in~\cite[Definitions~1.3 and~1.6]{BCDJComplexity}, where an algorithm chooses a current optimal solution and adds a cut removing it.

For a pure-integer polyhedron $P$, the \emph{full split rank} is the least $r$ with $\SC^r(P)=P_I$, where $\SC^0(P)=P$ and $\SC^{r+1}(P)=\SC(\SC^r(P))$; if no such $r$ exists, the rank is infinite. Basu, Conforti, Di Summa, and Jiang~\cite[Section~4]{BCDJComplexity} conjecture that the split rank of a rational pure-integer polyhedron is polynomially bounded in its input length, and note that this was open even in dimension three. The following theorem establishes the conjecture in every fixed dimension.

\begin{theorem}\label{thm:fixedproof}
Fix $d$. Let $P\subseteq\R^d$ be a rational polyhedron with all coordinates designated integer, and let $P_I=\conv(P\cap\Z^d)$, possibly empty or lower-dimensional. One can construct in polynomial time a split proof ending exactly at $P_I$, with total binary encoding length polynomial in the input length. The full split rank of $P$ is also polynomially bounded in the input length. The polynomial may depend on $d$.
\end{theorem}

Bounds on rank follow from bounds on proof length by monotonicity. If $Q_j$ denotes the relaxation after $j$ proof cuts, then
\[
 \SC^j(P)\subseteq Q_j
 \quad\Longrightarrow\quad
 \SC^{j+1}(P)\subseteq\SC(Q_j)\subseteq Q_{j+1},
\]
where the last inclusion holds because the next inequality is a split cut of $Q_j$. Every iterate retains $P_I$, so a proof ending at $Q_k=P_I$ gives $\SC^k(P)=P_I$, which bounds the rank in terms of the binary input length. The resulting description of the hull also permits polynomial-time optimization of any supplied rational linear objective. This is, of course, already possible by fixed-dimensional integer programming; what the theorem adds is a short split proof of the entire hull.

Let us describe the idea of the proof. We derive the defining inequalities of the integer hull one at a time, so consider one of them, written as $h\le0$ with $h$ affine and with integer coefficients. First suppose $P$ is bounded and has points on both sides of $h=1/2$. The polytope $C=P\cap\{h\ge1/2\}$ contains no integer point. We reduce $C$ to its face $h=1/2$ while retaining that face. Each step extends to $P$ while retaining all points with $h\le1/2$, and the final split $h\le0$ or $h\ge1$ proves the desired inequality.

To reduce $C$ to this face, flatness gives an integer direction in which only a constant number of sections at integer values can be nonempty. One might hope to treat the sections separately, but taking the convex hull with the other sections can bring back points already removed. We control this effect by bounding $2h-1$ on each section and interpolating these bounds to obtain a concave piecewise-linear function. Induction lets us lower each bound while retaining the points contributed by the other sections. A complete pass decreases the upper bound by a fixed factor; after polynomially many passes, a split on $2h-1$ leaves exactly the desired face. Integer-free polytopes are refuted by first reducing to a section with one integer coordinate fixed and then using induction.

For an unbounded $P$, we first project along any recession direction on which $h$ is constant. This reduces dimension, and every integer point of the projection has an integer lift to $P$. When no such direction remains, $P\cap\{h\ge1/2\}$ is bounded, so the preceding construction applies. The proof thus uses induction on dimension for both the sections and the recession directions.

For this argument to apply recursively, the split cuts on a section must retain the points contributed by the other sections, fractional ones included. The preservation property we need comes from a result of independent interest. A planar polyhedron is \emph{lattice-free} if its interior contains no integer point. Suppose that every point violating an inequality $h\le0$ maps, under two affine functions with integer coefficients, into the interior of such a polyhedron. The following theorem then gives a short split proof of the inequality, even when the original polytope has arbitrarily many variables.

\begin{theorem}\label{thm:planarsimulation}
Let $n$ be arbitrary and let $P\subseteq\R^n$ be a rational polytope with all coordinates designated integer. Suppose the input supplies an affine map $u:\R^n\to\R^2$ with integer coefficients, a rational affine function $h$, and a full-dimensional rational lattice-free polyhedron $K\subseteq\R^2$ such that
\[
h(z)>0\ \Longrightarrow\ u(z)\in\operatorname{int}K\qquad(z\in P).
\]
One can construct a split proof ending exactly at $P\cap\{h\le0\}$, with construction time, cut count, and total binary encoding length polynomial in $(P,u,h,K)$. Every original point satisfying the inequality, including fractional points, is retained.
\end{theorem}

The hypothesis explains why $h\le0$ is valid for the integer points of $P$: their images under $u$ are integer points, whereas every violating point maps into an interior containing no integer point. The theorem turns this geometric reason for validity into a polynomially encoded split proof. When $K$ is a strip between consecutive integer hyperplanes, a single split cut suffices. A general planar lattice-free polyhedron, however, need not fit inside such a strip, and the content of the theorem is that a polynomial sequence of split cuts is still enough.

Let us describe some consequences of Theorem~\ref{thm:planarsimulation}. First, consider the two-row system $Q(f,R)$ of Theorem~\ref{thm:uniformcorner}. Suppose a pure-integer rational polytope $P\subseteq\R^n$ is supplied with an affine map $u:P\to\R^2$ with integer coefficients and a rational affine map $s$ such that
\[
 s(z)\ge0,\qquad u(z)=f+Rs(z)\qquad(z\in P).
\]
Thus $(u(z),s(z))\in Q(f,R)$. We can enforce membership in the mixed-integer hull of $Q(f,R)$, the convex hull of its points with $x\in\Z^2$, by split cuts in the original variables $z$. The construction time and total proof encoding are polynomial in the supplied data, without fixing $n$ or the number of components of $s$.

The same theorem lets us replace a \emph{cross cut} by split cuts. Here a cross cut is an inequality valid on each of the four subsets obtained by choosing one branch from each of two split disjunctions. More generally, a supplied proof using cross cuts, the mixed-integer hull operations above, or the planar lattice-free certificates of Theorem~\ref{thm:planarsimulation} can be converted into a split proof in polynomial time. The input encoding includes the maps and disjunctions defining each step.

In fixed dimension, we can also construct the planar lattice-free polyhedron needed in the theorem. Let $P$ be a pure-integer rational polytope and let $u$ consist of two supplied affine functions with integer coefficients. We can derive any inequality valid at every point of $P$ with $u(z)\in\Z^2$, retaining every point that already satisfies it. We can also obtain exactly
\[
 U=\conv\{z\in P:u(z)\in\Z^2\}.
\]
If a rational polytope $G\subseteq P$ is to be retained, the construction instead gives $\conv(G\cup U)$. These constructions take polynomial time using classical fixed-dimensional integer programming and integer-hull algorithms~\cite{L83,H,CHKM}. This is the form of preservation that allows the induction to treat one section while retaining the points contributed by the others.

\subsection{Separation and optimization}
A split closure may have many facets even when membership can be tested efficiently. Given a rational point, \emph{strong separation} means either certifying that the point belongs to the split closure or returning a linear inequality valid for the closure and violated by that point. This use of ``strong'' has nothing to do with requiring each cut in a sequence to remove a current optimal LP vertex. We obtain polynomial-time separation when the total number of variables is fixed.

\begin{theorem}\label{thm:separation}
Fix $d$ and let $P\subseteq\R^d$ be a rational polyhedron. Given a rational inequality description of $P$, a set $I\subseteq\{1,\ldots,d\}$ of coordinates designated integer, and a rational query point, one can perform exact strong separation over $\SC_I(P)$ in time polynomial in the total encoding length of these data.
\end{theorem}

The algorithm reduces separation to polynomially many problems of minimizing a quadratic function over the integer points of a bounded rational polyhedron. Each has at most $d+1$ variables, so the fixed-dimensional algorithm of Ari and Hildebrand~\cite[Theorem~1.2]{AH} applies. None of the formulation or split-proof results above depends on this quadratic-optimization algorithm.

Separation, however, does not tell us which sequence of cuts reaches an optimal integer solution. For this sequential question, Basu, Conforti, Di Summa, and Jiang~\cite{BCDJ} give a polynomial strong split-cut algorithm in the plane; see also~\cite{BCDJComplexity}. Our algorithm fixes only the number of integer variables, allows arbitrary continuous dimension, and requires the objective to depend only on the integer variables. Its running time is pseudopolynomial: polynomial in the input length and linear in one plus the gap between the rounded LP optimum and the mixed-integer optimum.

\begin{theorem}\label{thm:gapalgorithm}
Fix $p\ge1$ and let $q\ge0$ be arbitrary. Let $P\subseteq\R^p\times\R^q$ be a rational polytope with $P\cap(\Z^p\times\R^q)\ne\varnothing$, and let $c\in\Z^p$ be nonzero and primitive. Let $L$ be the total input length of $(P,c)$ and define
\[
h(x,y)=c^Tx,\qquad z^*=\max_{(x,y)\in P\cap(\Z^p\times\R^q)}h(x,y),\qquad G=\lfloor\max_P h\rfloor-z^*.
\]
Given only $(P,c)$, an algorithm adds and retains $O_p(G+1)$ split cuts in a strong sequence for maximizing $h$. It returns a vertex of its final LP relaxation that is mixed-integer feasible and attains objective value $z^*$. Its running time is $(G+1)\operatorname{poly}_p(L)$, and each generated cut inequality and split disjunction has encoding length $\operatorname{poly}_p(L)$.
\end{theorem}

The algorithm exploits the fact that the objective takes integer values at mixed-integer feasible points. It first rounds the LP bound down to an integer $K$ and considers the points with $h=K$. Either it finds a mixed-integer feasible vertex at that level, which is then optimal, or it uses cuts to eliminate the entire level and impose $h\le K-1$. Inequalities constructed on the level are extended to the full relaxation so that each retained cut removes a current optimal vertex, while every original point with $h\le K-1$ is preserved. Consequently, once a level has been eliminated, the remaining relaxation is exactly $P\cap\{h\le K-1\}$, and the next step can start afresh from the original data and the new bound. This keeps coefficient lengths polynomial in the original input length.

Underlying the algorithm is a feasibility result that does not prescribe which objective the cuts must follow. It finds a mixed-integer feasible vertex, or proves that none exists, using a number of split cuts that depends only on the number of integer variables.

\begin{proposition}\label{prop:feasiblevertex}
For each fixed $p\ge0$, there is a constant $T_p$ and a polynomial-time algorithm with the following property. Given a rational polytope $Q\subseteq\R^p\times\R^q$, with arbitrary $q\ge0$ and only the first $p$ coordinates designated integer, it adds at most $T_p$ mixed-integer split cuts and either obtains an empty relaxation or returns a mixed-integer feasible vertex of the final relaxation. All coefficient encodings are polynomial in the full input length.
\end{proposition}

In particular, if $Q\cap(\Z^p\times\R^q)=\varnothing$, then
\[
\SC_I^{T_p}(Q)=\varnothing,\qquad I=\{1,\ldots,p\}.
\]
Indeed, after $j$ cuts the relaxation contains $\SC_I^j(Q)$. The bound is independent of the continuous dimension and of the numerical data. That some dimension-only bound on infeasibility proofs exists already follows from Cook, Coullard, and Tur\'an~\cite{CCT87} by projection onto the integer coordinates; Section~5.2 explains what is gained algorithmically here.

What if a short split proof is already available? Its cuts may leave every current optimal vertex untouched, so it need not give a strong sequence. We can convert it into one, again with a dependence on the objective gap.

Let $P$ be a nonempty pure-integer rational polytope, let $c$ be a nonzero integer vector, and let $\tau\in\Q$. Suppose a supplied weak split proof establishes $c^Tx\le\tau$, or emptiness, in $k$ cuts. Set
\[
 G=\max\{0,\lfloor\max_P c^Tx\rfloor-\lfloor\tau\rfloor\}.
\]
We can construct a strong sequence for maximizing $c^Tx$ that establishes $c^Tx\le\lfloor\tau\rfloor$, or emptiness, in at most $1+G(k+1)$ cuts. Here $G$ counts the integer objective levels above the target; the target need not be optimal. The running time is $(G+1)$ times a polynomial in the supplied encoding, and every generated cut and split has polynomial encoding length.

There is one setting in dimension three where we can remove this dependence on the gap: all points of the original relaxation that violate the desired inequality lie in a single closed split strip.

\begin{corollary}\label{cor:capstrip}
Let $P\subseteq\R^3$ be a rational polytope with all coordinates designated integer, let $h$ be rational affine, and let $\tau\in\Q$ satisfy $h\le\tau$ on $P\cap\Z^3$. Given primitive $\pi\in\Z^3$ and $k\in\Z$ such that
\begin{equation}\label{eq:originalstrip}
P\cap\{h>\tau\}\subseteq\{k\le\pi^Tx\le k+1\},
\end{equation}
one can establish $h\le\tau$, or emptiness, by a strong split-cut algorithm for maximizing $h$, with running time and coefficient encoding polynomial in the total encoding length of $(P,h,\tau,\pi,k)$, independently of the numerical objective gap.
\end{corollary}

A polynomial-time strong split-cut algorithm for general three-dimensional pure-integer polytopes remains open here. The corollary applies only under the stated condition on the original relaxation. For objectives involving continuous variables, even a finite split proof of the optimal objective bound need not exist, as the classical example of Cook, Kannan, and Schrijver~\cite{CKS} shows.

\paragraph{Organization.}
Section~\ref{sec:formulations} proves Theorems~\ref{thm:contrast}--\ref{thm:twointeger} on split-closure size and extended formulations; Appendix~\ref{sec:beyond} explores extensions and limits of this method. Section~\ref{sec:cornersimulation} proves Theorem~\ref{thm:planarsimulation} on planar lattice-free cuts, which Section~\ref{sec:fixedproof} uses to prove Theorem~\ref{thm:fixedproof}. Section~\ref{sec:algorithms} proves Theorem~\ref{thm:separation}, Proposition~\ref{prop:feasiblevertex}, Theorem~\ref{thm:gapalgorithm}, and Corollary~\ref{cor:capstrip}; its separation algorithm uses the bounds on split normals from Section~\ref{sec:normals}. A reader interested mainly in Theorem~\ref{thm:fixedproof} can read Sections~\ref{sec:introproofs}, \ref{sec:cornersimulation}, and~\ref{sec:fixedproof}, needing only Lemma~\ref{lem:bezouthalving} from Section~\ref{sec:formulations}.

\section{Split closures and extended formulations}\label{sec:formulations}\label{sec:cones}

The facet count alone does not determine the size of an extended formulation. For our tetrahedra, the split inequalities can be grouped into families indexed by consecutive integers. After clearing denominators, their coefficients are quadratic in this index, so repeated reflection gives a small formulation for each family. We begin with the geometric construction that makes many of these inequalities define distinct facets.

\subsection{A tetrahedral family}

On a translated simplicial cone, the effect of a split can be read directly from the points where its rays meet the boundary of the split strip.

Let $v\in\Q^d$ and let $R\in\Q^{d\times d}$ be nonsingular, with columns $r_1,\ldots,r_d$. Consider
\[
C=v+R\R_+^d.
\]
The point $v$ is the vertex of $C$, and every $x\in C$ has a unique representation $x=v+Rs$ with $s\ge0$. Here $s_j$ is the coefficient of $r_j$, and the corresponding ray is $v+t r_j$, $t\ge0$. The variables $s$ are continuous: integrality is imposed on the original variables $x$, not on $s$.

We first designate all coordinates of $x$ as integer. Fix a split normal $\pi\in\Z^d\setminus\{0\}$ and define
\[
z=R^T\pi,\qquad f=\{\pi^Tv\}.
\]
If $f=0$, every split with normal $\pi$ leaves $C$ unchanged. If $0<f<1$, the only integer right-hand side that can change $C$ is $k=\lfloor\pi^Tv\rfloor$. The resulting disjunctive hull has the description
\begin{equation}\label{eq:cut}
s\ge0,\qquad \alpha(\pi)^Ts\ge1,
\qquad
\alpha_j(\pi)=\max\left\{\frac{z_j}{1-f},-\frac{z_j}{f}\right\}.
\end{equation}
We call this inequality the cut with normal $\pi$; when $f=0$, no cut is required. To see the geometry behind the formula, follow the ray $v+t r_j$ from $v$. If $z_j>0$, it reaches the boundary of the split strip at $t=(1-f)/z_j$; if $z_j<0$, it reaches the other boundary at $t=-f/z_j$. A ray with $z_j=0$ is parallel to both boundary hyperplanes and contributes coefficient zero. The reciprocals of these parameters give the cut formula of Balas~\cite{B71}. For a split, \eqref{eq:cut} is in fact the complete description of the disjunctive hull, by~\cite[Lemma~1]{ACL}, and the same lemma shows that splits whose open split strips miss $v$ leave the cone unchanged. The following two properties will be used repeatedly:
\begin{equation}\label{eq:abs}
\alpha_j(\pi)\ge|z_j|\quad(0<f<1),
\qquad \alpha(-\pi)=\alpha(\pi).
\end{equation}

For example, take $v=(1/2,1/4)$, $r_1=(-1,2)$, $r_2=(1,1)$, and the split $x_1\le0$ or $x_1\ge1$ from Figure~\ref{fig:split}. Here $z=(-1,1)$ and $f=1/2$. Both rays leave the open split strip at parameter $t=1/2$, so the cut is $2s_1+2s_2\ge1$. Since $x_1+2x_2=1+3s_1+3s_2$, this is exactly $x_1+2x_2\ge5/2$. Formula~\eqref{eq:cut} expresses the same parameter calculation in any dimension.

Multiplying the cut inequality by $f(1-f)$ removes its denominators, which is convenient whenever the cut coefficients vary with the normal. For $0<f<1$, write the contributions from positive and negative scalar products as $U=\sum_{z_j\ge0}z_js_j$ and $V=\sum_{z_j<0}(-z_j)s_j$. Then
\begin{equation}\label{eq:clearedsplit}
\frac{U}{1-f}+\frac{V}{f}\ge1
\quad\Longleftrightarrow\quad
fU+(1-f)V-f(1-f)\ge0.
\end{equation}
For fixed $s$, if $f,U,V$ are affine in one index, the expression on the right of~\eqref{eq:clearedsplit} is a polynomial of degree at most two in that index. Likewise, on a region where each $\pi^Tr_j$ is nonnegative throughout or nonpositive throughout and $k<\pi^Tv<k+1$, substitute $f=\pi^Tv-k$: the expression has degree at most two in $(\pi,k)$. It is negative exactly when the split inequality is violated, because the multiplier $f(1-f)$ is positive.

\paragraph{The construction.}\label{sec:construction}

How can we make many split inequalities define distinct facets? We start with a family whose members are all tight at a common point $s>0$, and move the vertex $v$ slightly. In a suitable plane, the boundaries of the resulting inequalities become distinct tangents to a parabola. Each point of tangency satisfies just one of these inequalities at equality. We must still check the other split inequalities, however. We will show that a neighborhood of each tangency point, within the corresponding hyperplane, satisfies all of them.

Fix an integer $N\ge2$, and define
\begin{equation}\label{eq:construction}
\eps=\frac1{100N^2},\qquad
v=\begin{pmatrix}1/2\\\eps\\0\end{pmatrix},\qquad
R=\begin{pmatrix}N&-N&-N\\1&1&-1\\0&0&10N\end{pmatrix},\qquad
C_N=v+R\R_+^3.
\end{equation}
The bounded example is the simplex
\begin{equation}\label{eq:box}
P_N=v+R\{s\ge0:s_1+s_2+s_3\le4\}.
\end{equation}
The first two ray directions have slopes $\pm1/N$ in the $(x_1,x_2)$ plane, and the projection of the third direction onto this plane is the negative of the first. The third coordinate $10N$ has a different purpose. As we will see, it makes every cut~\eqref{eq:cut} with $\pi_3\ne0$ strict near the points at which we establish the facets. We can therefore concentrate on normals in the first two coordinates. The constraint $\sum_j s_j\le4$ makes $P_N$ bounded. In Section~\ref{sec:bounded}, we verify its description in the original coordinates and complete the proof of the following theorem.

\begin{theorem}\label{thm:main}
For each integer $N\ge2$, let $P_N\subseteq\R^3$ be defined by
\begin{align}
100Nx_1+100N^2x_2+20Nx_3&\ge50N+1,\label{eq:poly1}\\
-100Nx_1+100N^2x_2&\ge-50N+1,\label{eq:poly2}\\
x_3&\ge0,\label{eq:poly3}\\
100N^2x_2+20Nx_3&\le400N^2+1.\label{eq:poly4}
\end{align}
Then $P_N$ is a full-dimensional tetrahedron, and the displayed system has total binary encoding length $L_N=\Theta(\log N)$. Both split closures, $\SC(P_N)$ and $\SC_{\{1,2\}}(P_N)$, are full-dimensional and have at least $2N-1$ facets.
\end{theorem}

Let us return to the cone $C_N$. Since $\det R=20N^2$, its three ray directions are linearly independent, and we may use~\eqref{eq:cut}. For an integer normal $\pi=(p,q,r)$, the scalar products and fractional part in that formula are
\begin{equation}\label{eq:products}
z=(Np+q,-Np+q,-Np-q+10Nr),\qquad
f=\{p/2+q\eps\}.
\end{equation}
For each integer $t$ with $-N\le t\le N$, the primitive normal $(1,t,0)$ gives the cut $a(t)^Ts\ge1$, where
\begin{equation}\label{eq:curve}
a(t)=\left(
\frac{N+t}{1/2-\eps t},\quad
\frac{N-t}{1/2+\eps t},\quad
\frac{N+t}{1/2+\eps t}\right).
\end{equation}
The same expression is defined for real $t\in[-N,N]$. It is instructive to see first why $v$ must be displaced, so set $\eps=0$ for the moment. The cut then reads
\[
 2N(s_1+s_2+s_3)+2t(s_1-s_2+s_3)\ge1.
\]
Every intermediate inequality is a convex combination of the two endpoint inequalities, so this family does not yet give us the desired facets. All its inequalities are tight when $s_1+s_2+s_3=1/(2N)$ and $s_2=s_1+s_3$. Among these common points, choose
\[
 s^0=\left(\frac1{8N},\frac1{4N},\frac1{8N}\right).
\]
We now restore $\eps=1/(100N^2)$ and examine how the inequalities change. Put $c_\pm=1/2\pm N\eps$. To compare the cuts as $t$ varies, multiply their slack $a(t)^Ts-1$ by the common denominator $1/4-\eps^2t^2$, which is positive on $[-N,N]$. This gives
\[
 \eps(s_1+s_2-s_3+\eps)t^2+
 [c_+(s_1-s_2)+c_-s_3]t+
 \tfrac N2(s_1+s_2+s_3)-\tfrac14.
\]
For the calculation that follows, record
\begin{equation}\label{eq:parameters}
 c_\pm=\tfrac12\pm N\eps,\qquad \kappa=\frac{\eps}{4N}+\eps^2>0.
\end{equation}
We would like the coefficient of $t^2$ to be constant, so restrict attention to the plane $s_1+s_2-s_3=1/(4N)$ through $s^0$. On this plane that coefficient is $\kappa$. Using the other two coefficients, define affine coordinates
\[
 X=-\frac{c_+(s_1-s_2)+c_-s_3}{2\kappa},\qquad
 Y=\frac{\frac N2(s_1+s_2+s_3)-\frac14}{\kappa}.
\]
Thus $(1/4-\eps^2t^2)(a(t)^Ts-1)$ equals $\kappa(t^2-2tX+Y)$, and the cuts become
\[
 Y\ge2tX-t^2.
\]
These are precisely the tangents to the parabola $Y=X^2$. For an integer $k$ with $|k|<N$, the point $(X,Y)=(k,k^2)$ gives $\kappa(t-k)^2$ for this expression, so only the $k$th cut is tight. Solving for $s$ yields the point
\begin{equation}\label{eq:witness}
 s^{(k)}_3=\frac1{8N}+\frac{\kappa k^2}{N},\quad
 s^{(k)}_1=\frac1{8N}+\frac{N\eps s^{(k)}_3-\kappa k}{c_+},\quad
 s^{(k)}_2=\frac1{4N}+s^{(k)}_3-s^{(k)}_1.
\end{equation}
Because the displacement is small, these points stay close to $s^0$, and we will locate them in the common open region
\begin{equation}\label{eq:region}
 U=\left\{s\in\R^3:\frac9{10}s_j^0<s_j<1\text{ for }j=1,2,3\right\}.
\end{equation}
The $(X,Y)$ coordinates describe only this section; the integer restrictions still apply to the original coordinates $x$.

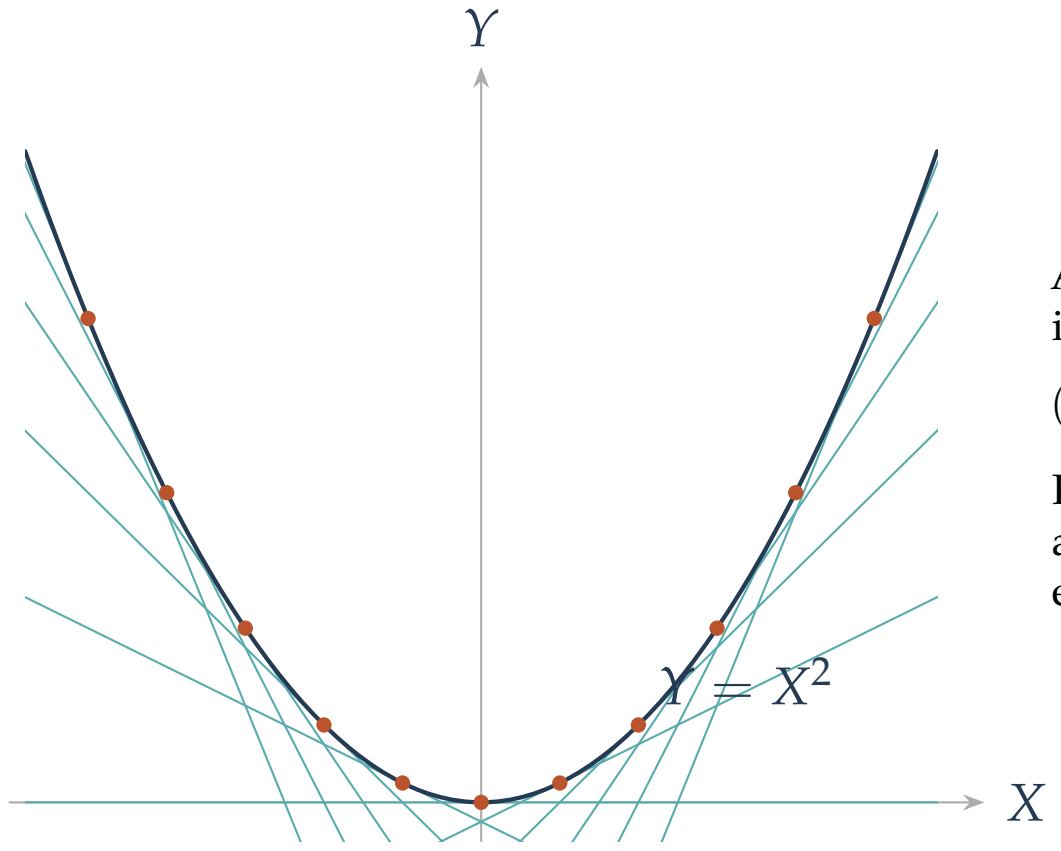
\begin{figure}[htbp]
\centering
\begin{adjustbox}{max width=\linewidth}
\begin{tikzpicture}[x=.65cm,y=.16cm,font=\small,scale=0.8]
\draw[axis] (-6,0)--(6.4,0) node[right] {$X$};
\draw[axis] (0,-2)--(0,38) node[above] {$Y$};
\begin{scope}
\clip (-5.8,-2) rectangle (5.8,36);
\foreach \k in {-5,...,5} {\draw[cutteal!65,thin] plot[domain=-6:6] (\x,{2*\k*\x-\k*\k});}
\draw[inkblue,thick] plot[domain=-5.8:5.8,samples=80] (\x,{\x*\x});
\foreach \k in {-5,...,5} {\fill[warmrust] (\k,{\k*\k}) circle(1.8pt);}
\end{scope}
\node[above right,text=inkblue] at (2,4) {$Y=X^2$};
\node[anchor=west,align=left,text width=4.4cm,font=\scriptsize] at (7,19) {At $(k,k^2)$ the slack of the $t$th line is\\[5pt] $(t-k)^2$.\\[7pt]Each marked point satisfies exactly one tangent inequality at equality.};
\end{tikzpicture}
\end{adjustbox}
\caption{The section $s_1+s_2-s_3=1/(4N)$ in its affine $(X,Y)$ coordinates, shown for $N=6$. The eleven points $(k,k^2)$, $-6<k<6$, lie on distinct tangent lines.}
\label{fig:facets}
\end{figure}

\begin{lemma}\label{lem:witness}
Let $k$ be an integer with $|k|<N$. The point $s^{(k)}$ defined in~\eqref{eq:witness} belongs to $U$, and the slack of each inequality in the family~\eqref{eq:curve} is given by
\begin{equation}\label{eq:identity}
 a(t)^Ts^{(k)}-1=\frac{\kappa(t-k)^2}{1/4-\eps^2t^2}\qquad(-N\le t\le N).
\end{equation}
\end{lemma}
\begin{proof}
The identity follows by substituting $(X,Y)=(k,k^2)$ in $\kappa(t^2-2tX+Y)$ and dividing by $1/4-\eps^2t^2$. The denominator is positive, so the slack vanishes exactly when $t=k$.

It remains to check that $s^{(k)}$ stays in $U$. We do this by comparing its coordinates with those of $s^0$. The choice of $\eps$ gives $N\eps\le1/200$ and $\kappa N\le\eps/3$. From~\eqref{eq:witness},
\[
 s^{(k)}_3\le\frac1{8N}+\frac\eps3<\frac1{6N},
 \qquad N\eps s^{(k)}_3<\frac\eps6.
\]
Using these bounds in the formula for $s^{(k)}_1$, together with $c_+\ge1/2$ and $|k|<N$, and then using the formula for $s^{(k)}_2$, we obtain
\[
\left|s^{(k)}_3-s^0_3\right|\le\eps/3,\qquad
\left|s^{(k)}_1-s^0_1\right|\le\eps,\qquad
\left|s^{(k)}_2-s^0_2\right|\le4\eps/3.
\]
Each coordinate therefore changes by less than $2\eps\le0.08\min_j s_j^0$. In particular, $s^{(k)}_j>0.9s_j^0$ for every $j$. The upper bounds are also less than one, and hence $s^{(k)}\in U$.
\end{proof}

We now turn to the remaining split cuts. There are infinitely many of them, so proving that each is strict at $s^{(k)}$ would not by itself give a neighborhood satisfying them all. The next lemma deals with this difficulty: outside a finite set of normals, every cut is strict throughout the same open set $U$.

\begin{lemma}\label{lem:finite}
Let
\begin{equation}\label{eq:list}
\mathcal F_N=\{(p,q,0)\in\Z^3:|p|\le2,\ |q|<3N\}.
\end{equation}
For every normal $\pi\notin\mathcal F_N$ with $0<f<1$, the cut~\eqref{eq:cut} is strictly satisfied at every point of $U$.
\end{lemma}

\begin{proof}
By~\eqref{eq:abs}, the absolute scalar products $|z_j|$ give lower bounds on the cut coefficients. Also, every $s\in U$ satisfies $s_j>\tfrac9{10}s_j^0$. Thus it suffices to bound $W=\sum_j|z_j|s_j^0$ from below. For $\pi=(p,q,r)$, put $t=q/N$. Substituting~\eqref{eq:products} gives
\begin{align*}
8W&=|p+t|+2|t-p|+|10r-p-t|\\
&\ge2|p|+|10r-2p|\ge10|r|.
\end{align*}
For the first inequality, pair one copy of $|t-p|$ with $|p+t|$ and the other with $|10r-p-t|$, and apply the triangle inequality twice. If $r\ne0$, it follows from~\eqref{eq:abs} that every $s\in U$ satisfies
\[
\alpha(\pi)^Ts>\frac9{10}W\ge\frac98>1.
\]
This deals with every normal having a nonzero third coordinate. When $r=0$, the same expression simplifies to
\[
W=\frac{|p+t|+|t-p|}{4}
=\frac12\max\{|p|,|t|\}.
\]
Consequently, whenever $\max\{|p|,|q|/N\}\ge3$,
\[
\alpha(\pi)^Ts>\frac9{10}W\ge\frac{27}{20}>1.
\]
These two cases cover every normal outside~\eqref{eq:list}.
\end{proof}

For the remaining normals in $\mathcal F_N$, the parity of the first coordinate determines which comparison to use. An odd first coordinate leads to comparison with the selected cuts. An even first coordinate gives a small fractional part, hence a large cut coefficient whenever the corresponding scalar product is negative.

\begin{lemma}\label{lem:all}
For every integer $k$ with $|k|<N$, the point $s^{(k)}$ satisfies every cut~\eqref{eq:cut} of $C_N$ with $0<f<1$. With all inequalities written with right-hand side one, the only inequality tight at $s^{(k)}$ is $a(k)^Ts\ge1$.
\end{lemma}

\begin{proof}
Fix an integer $k$ with $|k|<N$. Lemmas~\ref{lem:witness} and~\ref{lem:finite} already deal with every normal outside $\mathcal F_N$. We therefore consider $(p,q,0)\in\mathcal F_N$ with $0<f<1$.

If $p$ is odd, then $p\in\{-1,1\}$, and by~\eqref{eq:abs} we may negate the normal and assume $p=1$. Since $|q|<3N$,
\[
f=\frac12+q\eps\in(0,1).
\]
For $-N\le q\le N$, the assertion is exactly Lemma~\ref{lem:witness}. Outside this interval, we compare the cut with the nearer endpoint inequality. If $q>N$, its first and third coefficients are
\[
\frac{N+q}{1/2-q\eps},\qquad
\frac{N+q}{1/2+q\eps}.
\]
Both increase with $q$: their derivatives have positive numerators $c_+$ and $c_-$, respectively. The second coefficient is nonnegative, whereas the second coefficient of $a(N)$ is zero. Therefore $\alpha(1,q,0)\ge a(N)$ componentwise. If $q<-N$, the second coefficient is
\[
\frac{N-q}{1/2+q\eps},
\]
which decreases with $q$. The vector $a(-N)$ has only its second coordinate nonzero, so $\alpha(1,q,0)\ge a(-N)$ componentwise. Both endpoint inequalities are strict at $s^{(k)}$ by Lemma~\ref{lem:witness}.

Suppose next that $p$ is even, so $p\in\{-2,0,2\}$. If $q=0$, then $f=0$, contrary to the choice of a cut with $f\ne0$. Otherwise we may negate the normal to assume $q>0$. Now
\[
f=q\eps<\frac3{100N}.
\]
The small value of $f$ is useful when some $z_j<0$: since $z_j$ is an integer, the corresponding coefficient $-z_j/f$ is at least $1/f$. As $s^{(k)}\in U$, this single term gives
\[
\alpha(\pi)^Ts^{(k)}\ge\frac{s^{(k)}_j}{f}
>\frac{9/(80N)}{3/(100N)}=\frac{15}{4}>1.
\]
The only case left is that all $z_j$ are nonnegative. But $z_3=-z_1$ by~\eqref{eq:products}, so $z_1=z_3=0$ and $q=-Np$. As $q>0$, the only possibility is $(p,q)=(-2,2N)$. Its coefficient vector is
\[
\alpha(-2,2N,0)=\left(0,\frac{4N}{1-2N\eps},0\right)=a(-N),
\]
which is again strict at $s^{(k)}$.
\end{proof}

All but finitely many cuts are strict on $U$, and among the remaining inequalities only $a(k)^Ts\ge1$ is tight at $s^{(k)}$. Thus, near this point, the split closure is described by a single halfspace.

\begin{proposition}\label{prop:facets}
The split closure of $C_N$ is full-dimensional. For each integer $-N<k<N$, the inequality
\begin{equation}\label{eq:xfacet}
a(k)^TR^{-1}(x-v)\ge1
\end{equation}
defines a distinct facet of $\SC(C_N)$, with $v+Rs^{(k)}$ in its relative interior.
\end{proposition}

\begin{proof}
Fix an integer $k$ with $|k|<N$. We first find a neighborhood in which only the inequality $a(k)^Ts\ge1$ can exclude points. Since $s^{(k)}\in U$ and $U$ is open, choose an open ball about $s^{(k)}$ contained in $U$. Lemma~\ref{lem:finite} ensures that every cut with normal outside $\mathcal F_N$ is strict throughout this ball. Among the finitely many inequalities from $\mathcal F_N$, all except $a(k)^Ts\ge1$ are strict at its center by Lemma~\ref{lem:all}. By continuity, we can shrink the ball so that these inequalities remain strict throughout it.

Consequently, in the $s$ coordinates, the split closure agrees on this ball with the halfspace $a(k)^Ts\ge1$. The points in the ball with $a(k)^Ts>1$ form a nonempty open set, while its intersection with $a(k)^Ts=1$ is a relatively open two-dimensional set. The first observation proves full dimensionality; the second proves that the inequality defines a facet and that $s^{(k)}$ lies in its relative interior. The invertible change of coordinates $x=v+Rs$ gives the corresponding conclusions for~\eqref{eq:xfacet} and $v+Rs^{(k)}$.

To see that these facets are distinct, note that the first coordinate of $a(t)$ has derivative $c_+/(1/2-\eps t)^2>0$. Thus the coefficient vectors are distinct as $k$ varies. Since each inequality has right-hand side one, they define distinct hyperplanes.
\end{proof}

\subsection{Truncating the cone}\label{sec:bounded}

To pass from the cone to the tetrahedron, we use relaxations obtained by retaining inequalities with linearly independent normal vectors. Their split closures determine the original split closure. Let $P=\{x:Ax\ge b\}$ be nonempty, with $r=\operatorname{rank}A\ge1$. For every set $B$ of $r$ linearly independent rows of $A$, let $A_B,b_B$ denote the selected rows and right-hand sides. Retaining only these inequalities gives the relaxation $C_B=\{x:A_Bx\ge b_B\}$. By~\cite[Corollary~1]{ACL},
\begin{equation}\label{eq:corneridentity}
\SC_I(P)=\bigcap_B\SC_I(C_B).
\end{equation}
The intersection includes every such row set $B$, whether or not the equations $A_Bx=b_B$ have a solution in $P$. For a full-dimensional simplex, each relaxation retains the inequalities tight at one vertex; it is the \emph{supporting cone} at that vertex. In~\cite{ACL}, the family $\mathcal B_r^*$ consists of all independent $r$-row subsets, where $r=\operatorname{rank}A$; Corollary~1 is stated for arbitrary choices of the integer coordinates. Thus~\eqref{eq:corneridentity} does not require $r=d$, pointedness, or pure integrality; in particular, it applies in the presence of lineality.

We next record why a bounded truncation can preserve the split closure near the vertex, an observation that applies to either a pure-integer or a mixed-integer split closure.

\begin{lemma}\label{lem:truncation}
Let $v\in\Q^d$, let $R\in\Z^{d\times d}$ be nonsingular, and set
\[
C=v+R\R_+^d,\qquad
T=v+R\{s\ge0:\textstyle\sum_j s_j\le d+1\},\qquad
H=R\R_+^d.
\]
For every $I\subseteq\{1,\ldots,d\}$ and every $s\ge0$ with $\sum_j s_j\le d$,
\begin{equation}\label{eq:local-simplex}
v+Rs\in\SC_I(T)\quad\Longleftrightarrow\quad v+Rs\in\SC_I(C).
\end{equation}
Moreover,
\begin{equation}\label{eq:cone-recovery}
\SC_I(C)=\SC_I(T)+H.
\end{equation}
\end{lemma}

\begin{proof}
Monotonicity gives the forward implication in~\eqref{eq:local-simplex}. For the reverse, fix $x=v+Rs\in\SC_I(C)$ with $s\ge0$ and $\sum_i s_i\le d$. The relaxations in~\eqref{eq:corneridentity} are the $d+1$ cones at the vertices of $T$. Besides $C$, the cone at $v+(d+1)r_j$ has generating vectors $-r_j$ and $r_i-r_j$ for $i\ne j$. The coordinates of $x$ along these rays are, respectively,
\[
d+1-\sum_i s_i\ge1,\qquad s_i\quad(i\ne j).
\]
Fix a permitted normal $\pi\in\Z_I^d\setminus\{0\}$, and set $z=R^T\pi$ and $f=\{\pi^Tv\}$. All vertices of $T$ differ by integer vectors, so $\{\pi^Tw\}=f$ for each vertex $w$ of $T$. If $f=0$, the normal gives no cut at any vertex. If $0<f<1$ and $z_j\ne0$, the coefficient on $-r_j$ is at least $|z_j|\ge1$. Its coordinate is also at least one, so this term alone satisfies the cut. If $z_j=0$, that ray has coefficient zero and the other scalar products are $z_i-z_j=z_i$, so the cut at that vertex is exactly the cut at $v$, which $x$ satisfies. Identity~\eqref{eq:corneridentity} proves~\eqref{eq:local-simplex}.

For~\eqref{eq:cone-recovery}, take any $x=v+Rs\in\SC_I(C)$ and define $s'_j=\min\{s_j,1\}$. Every cut~\eqref{eq:cut} with $0<f<1$ remains satisfied at $s'$. Indeed, if a changed coordinate has positive cut coefficient, that coefficient is at least one because $R$ is integer, and $s'_j=1$ suffices. If all changed coordinates have zero coefficients, the left-hand side of the cut inequality is unchanged. Thus $v+Rs'\in\SC_I(C)$, and $\sum_j s'_j\le d$ gives $v+Rs'\in\SC_I(T)$ by~\eqref{eq:local-simplex}. Since $R(s-s')\in H$, this proves one inclusion. Conversely, $\SC_I(T)\subseteq\SC_I(C)$, and adding $H$ preserves every cone cut because its coefficients in $s$ are nonnegative, which gives the other inclusion. The argument also covers $I=\varnothing$, when there are no split cuts.
\end{proof}

We will use the two conclusions separately: local agreement transfers the facets just constructed to the simplex, while the global identity transfers extended formulations back to the cone by adding $d$ nonnegative variables. The matrix $R$ need not be unimodular; integrality remains in the original coordinates.

\begin{figure}[t]
\centering
\begin{adjustbox}{max width=\linewidth}
\begin{tikzpicture}[x=1.28cm,y=1.15cm,font=\small]
\fill[lightgray] (0,0)--(3,0)--(0,3)--cycle;
\fill[cutteal!16] (0,0)--(2,0)--(0,2)--cycle;
\draw[axis] (-.15,0)--(3.55,0) node[right] {$s_1$};
\draw[axis] (0,-.15)--(0,3.5) node[above] {$s_2$};
\draw[boundary] (3,0)--(0,3);
\draw[cut,dashed] (2,0)--(0,2);
\draw[cutteal,thin] (1,0)--(1,1)--(0,1);
\node[text=cutteal,align=center,font=\scriptsize] at (.62,.48) {closures agree\\for $s_1+s_2\le2$};
\node[text=inkblue,rotate=-42,font=\scriptsize] at (1.75,1.57) {$s_1+s_2=3$};
\fill[inkblue] (1,1) circle(1.8pt) node[above left] {$s'$};
\fill[inkblue] (2.6,2.1) circle(1.8pt) node[above] {$s$};
\draw[warmrust,-{Stealth[length=5pt]},line width=1pt] (2.54,2.06)--(1.07,1.05);
\node[text=warmrust,font=\scriptsize,align=left] at (2.6,2.8) {replace $s_j$ by\\$\min\{s_j,1\}$};

\end{tikzpicture}
\end{adjustbox}
\caption{The geometry of Lemma~\ref{lem:truncation}, shown in the $(s_1,s_2)$ plane. The gray triangle is $\{s\ge0:s_1+s_2\le3\}$; the teal triangle is the region of local agreement. The arrow replaces $s_j$ by $\min\{s_j,1\}$ for each $j$. The integer restrictions apply to $x=v+Rs$.}
\label{fig:simplex}
\end{figure}
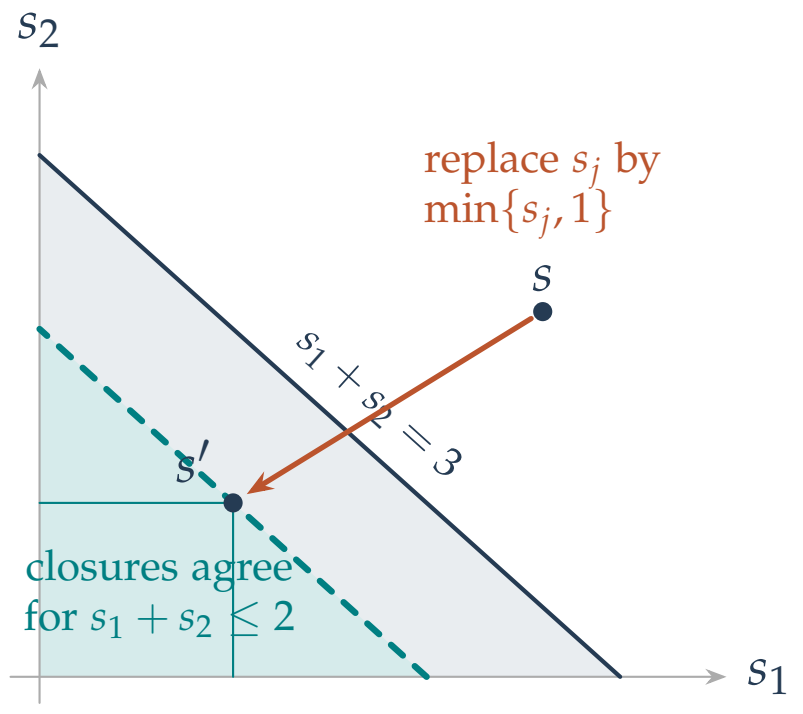

We now identify the tetrahedron in Theorem~\ref{thm:main} and transfer the facets to it. Set
\[
A=\begin{pmatrix}5&5N&1\\-1&N&0\\0&0&1\end{pmatrix}.
\]
For the $v$ and $R$ in~\eqref{eq:construction},
\begin{equation}\label{eq:matrix}
AR=\diag(10N,2N,10N),\qquad
Av=\begin{pmatrix}5/2+1/(20N)\\-1/2+1/(100N)\\0\end{pmatrix}.
\end{equation}
The first three inequalities of $P_N$ therefore say exactly that $s=R^{-1}(x-v)\ge0$. Also,
\[
s_1+s_2+s_3=x_2-\eps+\frac{x_3}{5N}.
\]
Consequently, the fourth inequality is $s_1+s_2+s_3\le4$, proving the representation~\eqref{eq:box}.

\begin{proof}[Proof of Theorem~\ref{thm:main}]
The representation~\eqref{eq:box} proves that $P_N$ is a full-dimensional tetrahedron. Since $R$ is integer, Lemma~\ref{lem:truncation} shows that $\SC(P_N)$ and $\SC(C_N)$ agree whenever $s\ge0$ and $s_1+s_2+s_3\le3$. Formula~\eqref{eq:witness} and $\kappa N\le\eps/3$ give
\[
s^{(k)}_1+s^{(k)}_2+s^{(k)}_3<\frac1N<3.
\]
Each $s^{(k)}$ has positive coordinates, so its neighborhood from Proposition~\ref{prop:facets} can be shrunk to lie in this region of agreement. The local facet descriptions therefore carry over to $\SC(P_N)$. It follows that $\SC(P_N)$ is full-dimensional and has at least $2N-1$ distinct facets.

For the mixed-integer statement, let
\[
Q_3=\SC(P_N),\qquad Q_2=\SC_{\{1,2\}}(P_N).
\]
We have $Q_3\subseteq Q_2$, since the latter closure uses fewer splits. Each inequality~\eqref{eq:xfacet} comes from a split with normal $(1,k,0)$, so it remains valid for $Q_2$. Let $H_k$ denote its halfspace. Then
\[
Q_3\subseteq Q_2\subseteq H_k.
\]
The relatively open two-dimensional subset of the facet of $Q_3$ on the boundary of $H_k$ also belongs to $Q_2$. As $Q_2$ contains the full-dimensional polyhedron $Q_3$, this proves that the same inequality defines a facet of $Q_2$. The $2N-1$ inequalities remain distinct.

All coefficients and right-hand sides in
\eqref{eq:poly1}--\eqref{eq:poly4} are integers of magnitude $O(N^2)$,
and at least one coefficient is $100N^2$. Their total binary encoding length is therefore $\Theta(\log N)$.
\end{proof}

The same containment argument, applied directly to $C_N$, proves that $\SC_{\{1,2\}}(C_N)$ also has the facets in~\eqref{eq:xfacet}. Thus the exponential bound already occurs for a translated simplicial cone; boundedness does not account for it.

This facet bound concerns the description of a single split closure. Extended formulations may be smaller, as we show next. Proof length and closure rank measure different aspects: the number of individual cuts and the number of rounds of applying all cuts from the chosen family, respectively. Dash~\cite{Dash10} gives exponential lower bounds for split-cut proof length, and Dash and Dubey~\cite{DD24} give linear rank lower bounds for generalizations of the split closure. Their lower-bound families have growing dimension, so these results do not give a facet lower bound for split closures in fixed dimension.

\subsection{Bounds on split normals}\label{sec:normals}

To describe the whole split closure, we need a finite set of normals whose cuts imply every other cut. Following Vielma~\cite{V07}, we use addition of normals to obtain implications between their cuts. This reduces the problem to finding bounded sets of integer generators.

Let $v\in\Q^d$ and $R\in\Q^{d\times m}$, with columns $r_1,\ldots,r_m$. Consider the system $x=v+Rs$, $s\ge0$, in the variables $(x,s)$. Designate the coordinates $x_i$, $i\in I$, as integer; all other variables are continuous. The columns of $R$ need not be independent, so $s$ need not be determined by $x$. For each permitted normal $\pi\in\Z_I^d\setminus\{0\}$ with $\{\pi^Tv\}\ne0$, formula~\eqref{eq:cut} describes the disjunctive hull in these variables~\cite[Lemma~1]{ACL}.

Two normals $\pi_1,\pi_2$ are \emph{sign-compatible} if, for each $j$, the scalar products $\pi_1^Tr_j$ and $\pi_2^Tr_j$ are both nonnegative or both nonpositive. The choice may differ with $j$, and zero is allowed. We also note that increasing the coefficients of $\alpha^Ts\ge1$ makes the inequality weaker on $s\ge0$, a fact used repeatedly below.

\begin{lemma}\label{lem:addition}
For the system just defined, let $\pi_1,\pi_2\in\Z_I^d$ be sign-compatible, with $\{(\pi_1+\pi_2)^Tv\}\ne0$. On $s\ge0$, their cuts with nonzero fractional part at $v$ imply the cut with normal $\pi_1+\pi_2$. A summand with integer value at $v$ requires no cut.
\end{lemma}

\begin{proof}
Write $f_i=\{\pi_i^Tv\}$ and $f=\{(\pi_1+\pi_2)^Tv\}\ne0$. If $f_1=0$, then $f=f_2$ and sign compatibility gives $\alpha(\pi_1+\pi_2)\ge\alpha(\pi_2)$ componentwise. The case $f_2=0$ is symmetric.

Otherwise we may assume $f_1+f_2<1$. Indeed, if the sum exceeds one, replace both normals by their negatives: their cuts are unchanged, the normals remain sign-compatible, and their new fractional parts sum to $2-f_1-f_2<1$. Equality is excluded by $f\ne0$.

Now $f=f_1+f_2$, and the weights $\lambda_i=f_i/f$ sum to one. For each $j$ with $\pi_i^Tr_j\le0$ for $i=1,2$, $\sum_i\lambda_i\alpha_j(\pi_i)=\alpha_j(\pi_1+\pi_2)$. For each $j$ with both products nonnegative, $f_i\le f$ gives
\[
\frac{f_i}{f(1-f_i)}\le\frac1{1-f},
\]
so, componentwise,
\[
\alpha(\pi_1+\pi_2)\ge\lambda_1\alpha(\pi_1)+\lambda_2\alpha(\pi_2).
\]
Consequently the two cuts give $\alpha(\pi_1+\pi_2)^Ts\ge\lambda_1+\lambda_2=1$.
\end{proof}

We now allow a lineality space in addition to the translated simplicial cone.

Let
\begin{equation}\label{eq:generalcone}
C=v+R\R_+^r+L\subseteq\R^d,
\end{equation}
where $v$ and $R$ are rational, $R$ has independent columns $r_1,\ldots,r_r$, and $L$ is a rational subspace such that $\R^d=\operatorname{im}(R)\oplus L$. We assume $r\ge1$. The subspace $L$ is the \emph{lineality space}: moving any distance in either direction along a vector of $L$ stays in the cone. The coefficients $s$ in $x=v+Rs+\ell$, $\ell\in L$, are unique. Fix a set $I$ of integer-variable indices and define
\[
W=\{\pi\in L^\perp:\pi_j=0\text{ for }j\notin I\}.
\]
The integer points $W\cap\Z^d$ form a lattice, whose rank is the number of vectors in a lattice basis. A split whose normal is not orthogonal to $L$ leaves $C$ unchanged: every point of $C$ lies on a line in $C$ meeting both branches of the disjunction, and is a convex combination of points in those branches. For $\pi\in W\cap\Z^d$ with $\{\pi^Tv\}\ne0$, the split formula~\eqref{eq:cut} remains valid, with $z=R^T\pi$ and $r$ variables $s_j$; see~\cite[Lemma~1]{ACL}. For each $\sigma\in\{-1,1\}^r$, set
\begin{equation}\label{eq:signcone}
K_\sigma=\{\pi\in W:\sigma_j r_j^T\pi\ge0\ (j=1,\ldots,r)\}.
\end{equation}
Each $K_\sigma$ is a pointed polyhedral cone: it contains no line. Indeed, a vector in $K_\sigma\cap(-K_\sigma)$ is orthogonal to both $\operatorname{im}(R)$ and $L$, and is therefore zero.

To generate all integer normals from a bounded region, write a normal as a nonnegative combination of integer ray generators and remove the integer parts of its coefficients. The remainder is still an integer vector and lies in a bounded parallelepiped. Lemma~\ref{lem:addition} then derives the original normal's cut from those of the generators and the remainder.

\begin{lemma}\label{lem:normalbound}
For fixed $d$, one can compute in polynomial time integers $M_1,\ldots,M_d$ of polynomial encoding length such that $\SC_I(C)$ is obtained using only the split normals satisfying $|\pi_j|\le M_j$. More precisely, if $K_\sigma$ is a simplicial cone with integer extreme-ray generators $g_1,\ldots,g_k$, its cuts are implied by those with
\begin{equation}\label{eq:simplebox}
|\pi_j|\le\sum_{i=1}^k |(g_i)_j|\qquad(j=1,\ldots,d).
\end{equation}
\end{lemma}

\begin{proof}
Fix a nonzero cone $K_\sigma$, and list its primitive integer extreme-ray generators $G_\sigma$. In fixed dimension, these generators are computed from the constant number of inequalities and equations in~\eqref{eq:signcone}. Rational elimination and determinant bounds yield polynomial encoding length for every generator.

Every $\pi\in K_\sigma\cap\Z^d$ is a nonnegative real combination of at most $d$ linearly independent generators. Write
\[
\pi=\sum_{i=1}^k\theta_i g_i
=\sum_{i=1}^k\lfloor\theta_i\rfloor g_i+\rho,
\qquad
\rho=\sum_{i=1}^k\{\theta_i\}g_i\in K_\sigma\cap\Z^d.
\]
The generators and the remainder have their $j$th coordinates bounded in magnitude by
\[
M_{\sigma,j}=d\max_{g\in G_\sigma}|g_j|.
\]
Consequently, the integer points of this box in $K_\sigma$ generate all integer points of $K_\sigma$ by nonnegative integer sums. To see that their cuts suffice, induct on the number of summands. A partial sum whose value at $v$ is integral imposes no cut. When a new partial sum has noninteger value at $v$, at least one of its two summands does as well. Lemma~\ref{lem:addition} then derives its cut from the available cuts, including the case where the other summand is integral. Thus the bounded generators imply every cut with nonzero fractional part at $v$ in $K_\sigma$.

Take the maximum of these bounds over the $2^r$ cones $K_\sigma$, omitting cones equal to $\{0\}$. If all cones are $\{0\}$, no split changes $C$; set every $M_j=0$. The resulting bounds have polynomial encoding length and are computable in polynomial time. When $K_\sigma$ is simplicial with the stated integer generators, its unique conic representation yields the sharper bound~\eqref{eq:simplebox} directly.
\end{proof}

The box may contain exponentially many normals, so the encoding bound does not permit us to list all their cuts in polynomial time. We next use the geometry of particular families to give small extended formulations. Section~\ref{sec:separation} instead searches the box by integer quadratic optimization.

\subsection{Formulations by reflection}\label{sec:ef}

Despite their many facets, the split closures of our tetrahedra have small extended formulations.

\begin{theorem}\label{thm:ef}
For the polytopes in Theorem~\ref{thm:main},
\[
\xc(\SC(P_N))=\Theta(\log N),\qquad
\xc(\SC_{\{1,2\}}(P_N))=\Theta(\log N).
\]
For either closure, a rational extended formulation with $O(\log N)$ inequalities and auxiliary variables can be constructed in time polynomial in $\log N$. Every scalar coefficient has encoding length $O(\log N)$, and the total encoding length is polynomial in $\log N$.
\end{theorem}

The upper bound uses the parabolic-polygon construction of~\cite[Lemma~1]{BFMP}. An inequality $c_0+c_1q+c_2q^2\ge0$ for every integer $q$ in an interval asks for a linear function to be nonnegative on the convex hull of the points $(q,q^2)$. Repeatedly halving the interval yields a small formulation of that hull. We give this formulation and its LP dual explicitly: the dual represents the whole family of inequalities using continuous auxiliary variables.

\begin{lemma}\label{lem:interval}
Let $a\le b$ be integers, and let $c_0,c_1,c_2$ be rational affine functions of real variables $y$. The inequalities
\begin{equation}\label{eq:interval}
c_0(y)+c_1(y)q+c_2(y)q^2\ge0
\qquad(q\in[a,b]\cap\Z)
\end{equation}
have a linear extended formulation with $O(\log(b-a+2))$ inequalities and auxiliary variables. The formulation is constructible in polynomial time in the encoding of $a,b,c_0,c_1,c_2$, and has polynomial encoding length.
\end{lemma}

\begin{proof}
Set $n=b-a$. If $n=0$, use the single inequality at $q=a$. Suppose $n\ge1$.

In this proof, $j$ ranges over integers. Let $M_n=\conv\{(j,j^2):0\le j\le n\}$. Set $m=\lfloor n/2\rfloor$ and consider
\[
\rho_n(x,z)=(n-x,z-2nx+n^2).
\]
This affine map sends $(j,j^2)$ to $(n-j,(n-j)^2)$, so $M_n=\conv(M_m\cup\rho_n(M_m))$. Since $x\le m\le n/2$ on $M_m$, the relation
\[
(x',z')=(x+t,z+nt),\qquad (x,z)\in M_m,\qquad 0\le t\le n-2x
\]
has image exactly this convex hull. For each $(x,z)$, it traces the segment to $\rho_n(x,z)$, proving one inclusion. Its image is convex and contains both $M_m$ and $\rho_n(M_m)$, proving the other. 

To build $M_n$ from $M_0$, perform these halvings in reverse order. Let $n_1,\ldots,n_k$ be the positive terms of $n,\lfloor n/2\rfloor,\lfloor n/4\rfloor,\ldots$, listed in reverse. Thus $n_1=1$, $n_k=n$, $\lfloor n_i/2\rfloor=n_{i-1}$ for $i>1$, and $k=\lfloor\log_2 n\rfloor+1$. Define
\begin{equation}\label{eq:momentlift}
T_n=\left\{t\in\R_+^k:
t_i+2\sum_{j<i}t_j\le n_i\quad(i=1,\ldots,k)\right\}.
\end{equation}
Applying these relations successively from $M_0=\{(0,0)\}$ shows that the linear map
\begin{equation}\label{eq:momentmap}
x=\sum_i t_i,\qquad z=\sum_i n_i t_i
\end{equation}
maps $T_n$ onto $M_n$.

After substituting $q=a+j$, set
\[
\beta_0=c_0+ac_1+a^2c_2,\qquad
\beta_1=c_1+2ac_2,\qquad \beta_2=c_2.
\]
Condition~\eqref{eq:interval} is equivalent to
\[
\beta_0+\min_{t\in T_n}\sum_i(\beta_1+n_i\beta_2)t_i\ge0.
\]
Since $T_n$ is nonempty and bounded, LP duality expresses this lower bound by nonnegative multipliers $\lambda_i$ for the inequalities $t_i+2\sum_{j<i}t_j\le n_i$ in~\eqref{eq:momentlift}. The entire indexed family is therefore equivalent to
\begin{align}
\lambda_i&\ge0 &&(i=1,\ldots,k),\label{eq:dual1}\\
\beta_1+n_i\beta_2+\lambda_i+2\sum_{j>i}\lambda_j&\ge0
&&(i=1,\ldots,k),\label{eq:dual2}\\
\beta_0-\sum_i n_i\lambda_i&\ge0.\label{eq:dual3}
\end{align}
These are $2k+1$ linear inequalities in $y$ and $\lambda$. The integers $n_i$ are obtained by binary halving, and all displayed coefficients have polynomial encoding length in the given data.
\end{proof}

\begin{figure}[t]
\centering
\begin{adjustbox}{max width=\linewidth}
\begin{tikzpicture}[font=\small]
\begin{scope}[x=.72cm,y=.12cm]
\fill[lightgray] (0,0)--(1,1)--(2,4)--(3,9)--(4,16)--(5,25)--(6,36)--cycle;
\fill[cutteal!25] (0,0)--(1,1)--(2,4)--(3,9)--cycle;
\fill[warmrust!20] (3,9)--(4,16)--(5,25)--(6,36)--cycle;
\draw[boundary] (0,0)--(1,1)--(2,4)--(3,9)--(4,16)--(5,25)--(6,36)--cycle;
\draw[cut] (0,0)--(1,1)--(2,4)--(3,9)--cycle;
\draw[warmrust,line width=1pt] (3,9)--(4,16)--(5,25)--(6,36)--cycle;
\draw[gray!65,dashed,thin] (1,1)--(5,25) node[midway,below right,font=\scriptsize] {$0\le t\le4$};
\foreach \j in {0,...,6} {\pgfmathsetmacro{\jj}{\j*\j}\fill[inkblue] (\j,\jj) circle(1.5pt);}
\draw[axis] (-.25,0)--(6.65,0) node[right] {$x$};
\draw[axis] (0,-1.5)--(0,39) node[above] {$z$};
\foreach \j in {2,4,6} {\node[below,font=\scriptsize] at (\j,0) {$\j$};}
\node[text=cutteal,anchor=east] at (2,5) {$M_3$};
\node[text=warmrust,anchor=west] at (4.7,22) {$\rho_6(M_3)$};
\node[text=inkblue] at (2.1,24) {$M_6$};
\end{scope}
\begin{scope}[xshift=6.3cm]
\node[anchor=north west,text=inkblue,align=left,text width=6.3cm,font=\scriptsize] at (0,4.5) {\textbf{One new variable joins reflected points.}\\[7pt]
The parameter $t$ traverses the segment from $(x,z)$ to $\rho_n(x,z)$.\\[7pt]
The proof derives this step and then repeats it at the successive half-lengths.};
\end{scope}
\end{tikzpicture}
\end{adjustbox}
\caption{The parabolic formulation of~\cite{BFMP}, in the form used in Lemma~\ref{lem:interval}. Here $M_j=\conv\{(i,i^2):0\le i\le j\}$. For $n=6$, the affine map $\rho_6(x,z)=(6-x,z-12x+36)$ sends the left polygon to the right polygon. Adding the joining segments gives their convex hull, $M_6$. The two axes use different scales; all plotted coordinates are exact.}
\label{fig:parabola}
\end{figure}

By~\eqref{eq:corneridentity}, it suffices to formulate the split closures of the four vertex cones of $P_N$. We bound a sufficient set of normals at every vertex, then partition those normals into a constant number of intervals of an integer index. Each interval has length $O(N)$, and its cut inequalities after multiplication by $f(1-f)$ have degree at most two in the index. Lemma~\ref{lem:interval} therefore gives total formulation size $O(\log N)$.

\begin{proof}[Proof of Theorem~\ref{thm:ef}]
\emph{Bound the normals at every vertex.}
Use $v,R,\eps$ from~\eqref{eq:construction}. Up to sign, the four primitive integer facet normals of $P_N$ are
\[
a_1=(5,5N,1),\quad a_2=(-1,N,0),\quad
a_3=(0,0,1),\quad a_4=(0,5N,1).
\]
At each vertex, the normals of the three tight facet inequalities form the dual basis of its generating vectors, up to positive scaling and signs. Every cone $K_\sigma$ of normals is therefore generated by signed copies of those three normals. Lemma~\ref{lem:normalbound} bounds each coordinate by the sum of its absolute values on these generators. Summing over all four displayed vectors yields the common sufficient box
\begin{equation}\label{eq:familybox}
\pi=(p,q,r),\qquad |p|\le6,\quad |q|\le11N,\quad |r|\le3.
\end{equation}
Taking every normal in this box adds only valid split cuts. There are $91$ choices of $(p,r)$, so only the index $q$ has a range growing with $N$.

For the mixed-integer closure, we must derive all cuts using only normals with $r=0$. The generating vectors of the cones at the vertices are either $r_i$, or $-r_j$ and $r_i-r_j$ for $i\ne j$. Restricted to $r=0$, their products vanish only on the lines
\[
p=0,\qquad q=0,\qquad q=Np,\qquad q=-Np.
\]
Thus the primitive generators of the extreme rays of each cone $K_\sigma\cap\{r=0\}$ belong to
\[
\{\pm(1,0),\ \pm(0,1),\ \pm(1,N),\ \pm(1,-N)\}.
\]
Each cone has at most two such generators. Lemma~\ref{lem:normalbound} gives the sufficient bounds $|p|\le2$, $|q|\le2N$. In particular, the larger box~\eqref{eq:familybox} with $r=0$ suffices at all four vertices.

\emph{Formulate the cuts on each index interval.}
Write the vertices as $v_0=v$ and $v_j=v+4r_j$ for $j=1,2,3$. At vertex zero the coordinates and scalar products are $s_i$ and $z_i(q)$ from~\eqref{eq:products}. At vertex $j\ge1$ they are
\[
\begin{array}{c|c|c}
\text{ray}&\text{coordinate}&\text{scalar product}\\ \hline
-r_j&4-\sum_i s_i&-z_j(q)\\
r_i-r_j\ (i\ne j)&s_i&z_i(q)-z_j(q).
\end{array}
\]
The table gives three facts at every vertex: the ray coordinates are nonnegative on~\eqref{eq:box}, the scalar products are affine in $q$, and the vertex differs from $v$ by an integer vector. Thus the value of the normal at every vertex has the same fractional part $f(q)=\{p/2+\eps q\}$.

Fix $(p,r)$ in~\eqref{eq:familybox} and one vertex. On $|q|\le11N$ we have $|\eps q|\le11/200<1/2$. For odd $p$, $f(q)=1/2+\eps q$. For even $p$, it is $\eps q$ when $q>0$ and $1+\eps q$ when $q<0$; the index $q=0$ gives no cut. Each of the three scalar products is affine in $q$, so it has at most one zero unless it is identically zero. Partitioning the integer interval at these zeros and at zero, retaining a breakpoint as a singleton when necessary, gives at most nine consecutive integer intervals, from which we discard indices with $f=0$. On each remaining interval, $f$ is affine and each scalar product is nonnegative throughout or nonpositive throughout.

Let $t_i(s)$ denote the coefficients along the generating vectors at this vertex and $w_i(q)$ the corresponding products. Set
\[
U(q,s)=\sum_{w_i\ge0}t_i(s)w_i(q),\qquad
V(q,s)=\sum_{w_i<0}t_i(s)(-w_i(q)).
\]
Zero products may be assigned either sign. By~\eqref{eq:clearedsplit}, the cut is a quadratic inequality in $q$ with coefficients affine in $s$. Lemma~\ref{lem:interval} therefore represents all cuts in one interval using $O(\log N)$ inequalities and variables.

\emph{Assemble the full closure.}
Use separate auxiliary variables for each interval, each of the $91$ pairs $(p,r)$, and each of the four vertices. Add $s\ge0$, $\sum_i s_i\le4$, and $x=v+Rs$. The cuts with normals in~\eqref{eq:familybox} imply all other cuts at each vertex. Identity~\eqref{eq:corneridentity} therefore shows that the projection is exactly $\SC(P_N)$. Restricting the normal set to $r=0$ yields exactly $\SC_{\{1,2\}}(P_N)$.

There are a constant number of intervals, each of length at most $22N$ and with endpoints of magnitude $O(N)$. Every coefficient in~\eqref{eq:clearedsplit}, its translation to an interval starting at zero, and the dual system~\eqref{eq:dual1}--\eqref{eq:dual3} has $O(\log N)$ bits. This proves the claimed size, encoding, and construction bounds. The extension is bounded: the simplex bounds $s$, and then $\lambda_i\ge0$ and $\sum_i n_i\lambda_i\le\beta_0(s)$, with $n_i\ge1$, bound every auxiliary variable.

\emph{Prove the matching lower bound.}
A polyhedron described by $k$ inequalities has at most $2^k$ faces. Each facet of a projection has a distinct nonempty inverse-image face in the extension. Theorem~\ref{thm:main} therefore yields $2^k\ge2N-1$ for either closure, and hence $k\ge\log_2(2N-1)$.
\end{proof}

Our auxiliary variables are continuous, and projection gives exactly the split closure of the original polyhedron with its original designation of integer coordinates. Bodur, Dash, and G\"unl\"uk~\cite{BDG17} instead study how extensions of a relaxation strengthen cuts; Dash, G\"unl\"uk, and Hildebrand~\cite{DGH18} study binary extensions. Eisenbrand, Fiorini, Rohwedder, and Wei~\cite{EFRW26} give small polyhedral extended formulations approximating Gomory closures of packing problems, rather than exact formulations of split closures.

\subsection{Two-row systems}\label{sec:twointeger}

For rational $f\in\Q^2$ and $R\in\Q^{2\times n}$, recall the polyhedron
\[
Q(f,R)=\{(x,s)\in\R^2\times\R_+^n:x=f+Rs\}.
\]
The coordinates $x_1,x_2$ are designated integer and $s$ is continuous. We seek an extended formulation of the split closure with $O(L)$ inequalities and auxiliary variables, where $L$ is the total encoding length of $f,R$.

We cover $\R^2$ by cones of sign-compatible normals. In each cone, integer generators can be arranged in arithmetic progressions along the bounded edges of an integer hull. Lemma~\ref{lem:addition} reduces the split closure to their cuts, and a small linear system represents the cuts in each progression.

\paragraph{Lattice points on bounded edges.}

The next lemma decomposes a primitive integer vector with positive coordinates into two nonnegative integer vectors whose determinant has absolute value one. We use it later in the paper as well. This classical lattice geometry is closely related to continued fractions~\cite{SV18}.

\begin{lemma}\label{lem:primitivepair}\label{lem:bezouthalving}
Let $v=(a,b)$ have positive coprime integer coordinates. There are nonzero vectors $r,s\in\Z_+^2$ such that
\[
 v=r+s,\qquad |\det(r,s)|=1.
\]
Unless $v=(1,1)$, the vectors are comparable coordinatewise, so the smaller is at most $v/2$. If $b\ge2$, this smaller vector $(a',b')$ satisfies $a'\ge0$, $b'\ge1$, and $|ab'-a'b|=1$. 
\end{lemma}
\begin{proof}
For $b\ge2$, the extended Euclidean algorithm gives $d\in\{1,\ldots,b-1\}$ with $ad\equiv1\pmod b$. Let $c=(ad-1)/b$, and take $s=(c,d)$ and $r=v-s$. Then $0\le c<a$, both second coordinates are positive, and $\det(r,s)=1$. For $b=1$, use $(1,0)$ and $(a-1,1)$.

If two nonnegative integer vectors are incomparable, after interchanging them they have the form $(x+e,y)$ and $(x,y+f)$, with $e,f\ge1$ and $x,y\ge0$. Their absolute determinant $xf+ey+ef$ equals one only when $x=y=0$ and $e=f=1$, giving sum $(1,1)$. Otherwise the smaller vector is at most half the sum; its determinant with the sum is $\pm1$.
\end{proof}

We next work in a rank-two lattice, identified with $\Z^2$ by a lattice basis. For a pointed full-dimensional rational cone $K\subseteq\R^2$, let $p_0,p_k$ be its primitive extreme-ray generators, ordered so that $\Delta=\det(p_0,p_k)>0$. The union of the bounded edges of
\[
H=\conv\bigl((K\cap\Z^2)\setminus\{0\}\bigr)
\]
is a polygonal chain from $p_0$ to $p_k$. List the vertices of this chain as $p_0,p_1,\ldots,p_k$. We also use the lattice points in the relative interiors of its edges.

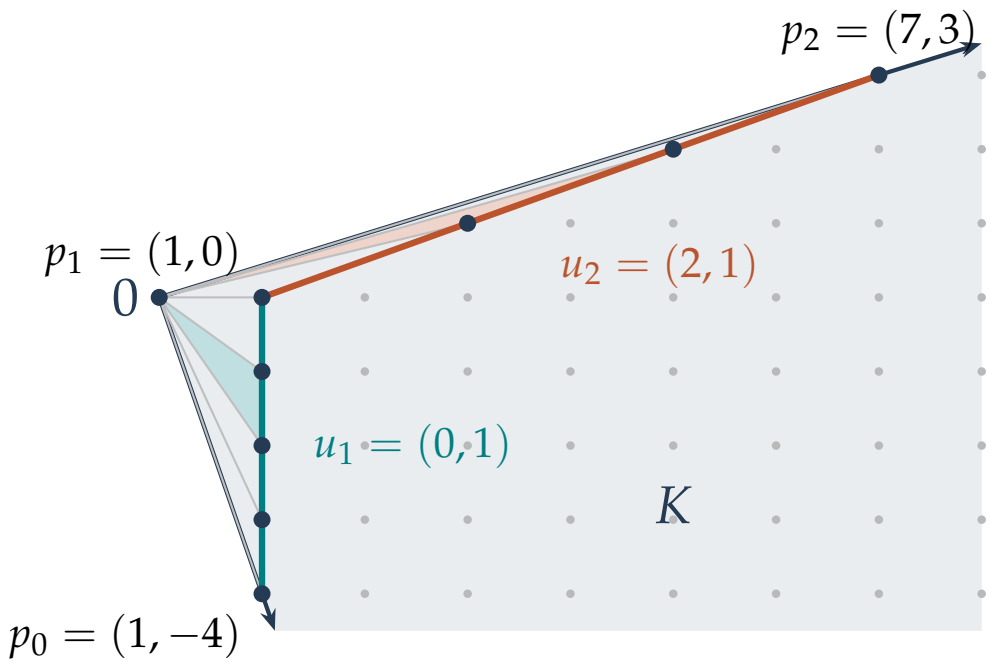
\begin{figure}[htbp]
\centering
\begin{adjustbox}{max width=\linewidth}
\begin{tikzpicture}[font=\small]
\begin{scope}[x=.68cm,y=.49cm]
\fill[lightgray] (0,0)--(1.125,-4.5)--(8,-4.5)--(8,3.428571)--cycle;
\fill[cutteal!25] (0,0)--(1,-2)--(1,-1)--cycle;
\fill[warmrust!25] (0,0)--(3,1)--(5,2)--cycle;
\fill[gray!55] (1.00000000,-4.00000000) circle(.8pt);
\fill[gray!55] (1.00000000,-3.00000000) circle(.8pt);
\fill[gray!55] (1.00000000,-2.00000000) circle(.8pt);
\fill[gray!55] (1.00000000,-1.00000000) circle(.8pt);
\fill[gray!55] (1.00000000,0.00000000) circle(.8pt);
\fill[gray!55] (2.00000000,-4.00000000) circle(.8pt);
\fill[gray!55] (2.00000000,-3.00000000) circle(.8pt);
\fill[gray!55] (2.00000000,-2.00000000) circle(.8pt);
\fill[gray!55] (2.00000000,-1.00000000) circle(.8pt);
\fill[gray!55] (2.00000000,0.00000000) circle(.8pt);
\fill[gray!55] (3.00000000,-4.00000000) circle(.8pt);
\fill[gray!55] (3.00000000,-3.00000000) circle(.8pt);
\fill[gray!55] (3.00000000,-2.00000000) circle(.8pt);
\fill[gray!55] (3.00000000,-1.00000000) circle(.8pt);
\fill[gray!55] (3.00000000,0.00000000) circle(.8pt);
\fill[gray!55] (3.00000000,1.00000000) circle(.8pt);
\fill[gray!55] (4.00000000,-4.00000000) circle(.8pt);
\fill[gray!55] (4.00000000,-3.00000000) circle(.8pt);
\fill[gray!55] (4.00000000,-2.00000000) circle(.8pt);
\fill[gray!55] (4.00000000,-1.00000000) circle(.8pt);
\fill[gray!55] (4.00000000,0.00000000) circle(.8pt);
\fill[gray!55] (4.00000000,1.00000000) circle(.8pt);
\fill[gray!55] (5.00000000,-4.00000000) circle(.8pt);
\fill[gray!55] (5.00000000,-3.00000000) circle(.8pt);
\fill[gray!55] (5.00000000,-2.00000000) circle(.8pt);
\fill[gray!55] (5.00000000,-1.00000000) circle(.8pt);
\fill[gray!55] (5.00000000,0.00000000) circle(.8pt);
\fill[gray!55] (5.00000000,1.00000000) circle(.8pt);
\fill[gray!55] (5.00000000,2.00000000) circle(.8pt);
\fill[gray!55] (6.00000000,-4.00000000) circle(.8pt);
\fill[gray!55] (6.00000000,-3.00000000) circle(.8pt);
\fill[gray!55] (6.00000000,-2.00000000) circle(.8pt);
\fill[gray!55] (6.00000000,-1.00000000) circle(.8pt);
\fill[gray!55] (6.00000000,0.00000000) circle(.8pt);
\fill[gray!55] (6.00000000,1.00000000) circle(.8pt);
\fill[gray!55] (6.00000000,2.00000000) circle(.8pt);
\fill[gray!55] (7.00000000,-4.00000000) circle(.8pt);
\fill[gray!55] (7.00000000,-3.00000000) circle(.8pt);
\fill[gray!55] (7.00000000,-2.00000000) circle(.8pt);
\fill[gray!55] (7.00000000,-1.00000000) circle(.8pt);
\fill[gray!55] (7.00000000,0.00000000) circle(.8pt);
\fill[gray!55] (7.00000000,1.00000000) circle(.8pt);
\fill[gray!55] (7.00000000,2.00000000) circle(.8pt);
\fill[gray!55] (7.00000000,3.00000000) circle(.8pt);
\fill[gray!55] (8.00000000,-4.00000000) circle(.8pt);
\fill[gray!55] (8.00000000,-3.00000000) circle(.8pt);
\fill[gray!55] (8.00000000,-2.00000000) circle(.8pt);
\fill[gray!55] (8.00000000,-1.00000000) circle(.8pt);
\fill[gray!55] (8.00000000,0.00000000) circle(.8pt);
\fill[gray!55] (8.00000000,1.00000000) circle(.8pt);
\fill[gray!55] (8.00000000,2.00000000) circle(.8pt);
\fill[gray!55] (8.00000000,3.00000000) circle(.8pt);
\draw[boundary,-{Stealth[length=4pt]}] (0,0)--(1.125,-4.5);
\draw[boundary,-{Stealth[length=4pt]}] (0,0)--(8,3.428571);
\foreach \y in {-4,-3,-2,-1,0} {\draw[gray!50,thin] (0,0)--(1,\y);}
\foreach \p in {(3,1),(5,2),(7,3)} {\draw[gray!50,thin] (0,0)--\p;}
\draw[cut] (1,-4)--(1,0);
\draw[warmrust,line width=1.2pt] (1,0)--(7,3);
\foreach \p in {(1,-4),(1,-3),(1,-2),(1,-1),(1,0),(3,1),(5,2),(7,3)} {\fill[inkblue] \p circle(1.6pt);}
\fill[inkblue] (0,0) circle(1.6pt) node[left] {$0$};
\node[below left,font=\scriptsize] at (1,-4) {$p_0=(1,-4)$};
\node[above left,font=\scriptsize] at (1,0) {$p_1=(1,0)$};
\node[above,font=\scriptsize] at (7,3) {$p_2=(7,3)$};
\node[font=\scriptsize,text=cutteal,anchor=west] at (1.3,-2) {$u_1=(0,1)$};
\node[font=\scriptsize,text=warmrust,anchor=north west] at (3.7,1) {$u_2=(2,1)$};
\node[text=inkblue] at (5,-2.8) {$K$};
\end{scope}
\begin{scope}[xshift=6.65cm]
\node[anchor=north west,text=inkblue,align=left,text width=5.8cm,font=\scriptsize] at (0,1.9) {\textbf{Growth of the determinants.}\\[6pt]
Consecutive lattice points on the chain have determinant one. At the internal vertex,
\[
u_2-u_1=(2,0)=2p_1.
\]
With $D_i=\det(p_0,p_i)$,
\[
D_1=4,\qquad D_2=31\ge(3+1)D_1.
\]
The lattice lengths satisfy\\[3pt]
$(4+1)(3+1)=20\le2\cdot31$.};
\end{scope}
\end{tikzpicture}
\end{adjustbox}
\caption{The bounded edges of $H$ in the space of split normals. The extreme rays have primitive generators $(1,-4)$ and $(7,3)$. The marked lattice points on the bounded edges generate all integer points in the cone. The lattice points on each edge form an arithmetic progression.}
\label{fig:chain}
\end{figure}

Consecutive lattice points on the chain form a lattice basis, so these points generate all integer points in $K$. For an edge from $p_{i-1}$ to $p_i$, write $p_i-p_{i-1}=\ell_i u_i$, where $u_i$ is primitive and $\ell_i$ is a positive integer. The integer $\ell_i$ is its \emph{lattice length}: the edge consists of $\ell_i$ primitive integer steps. Figure~\ref{fig:chain} gives an example. The formulation size will depend on the sum of logarithms of these lengths, which the next lemma bounds in terms of $\Delta$.

\Needspace{11\baselineskip}
\begin{lemma}\label{lem:chain}
The lattice points on this chain generate $K\cap\Z^2$ by nonnegative integer combinations. Write the lattice points on its edges as
\[
p_{i-1}+t u_i\quad(t=0,\ldots,\ell_i),\qquad i=1,\ldots,k,
\]
where $u_i$ is primitive and $\ell_i\ge1$ is the lattice length. Then
\begin{equation}\label{eq:product}
\prod_{i=1}^k(\ell_i+1)\le2\Delta,
\qquad
\sum_{i=1}^k\log_2(\ell_i+2)=O(\log(2\Delta)).
\end{equation}
The vertices $p_i$, steps $u_i$, and lengths $\ell_i$ can be computed in polynomial time from a rational description of $K$, without listing the lattice points in the relative interiors of the edges.
\end{lemma}

\begin{proof}
Write $p=p_0$ and $q=p_k$. We first obtain the boundary from a bounded polygon. Let
\[
 w(x)=\det(p,x)+\det(x,q),\qquad
 T=\{x\in K:1\le w(x)\le\Delta\}.
\]
For $x=\alpha p+\beta q\in K$, we have $w(x)=\Delta(\alpha+\beta)$. Thus $T$ is the triangle $\conv(0,p,q)$ with a neighborhood of the origin removed. Every nonzero integer point in that triangle belongs to $T$, since $w$ is positive and integer-valued there. If $x\in K$ has $w(x)\ge\Delta$, then
\[
 x=\frac{x}{\alpha+\beta}
   +\left(1-\frac1{\alpha+\beta}\right)x\in[p,q]+K.
\]
Also $H+K=H$: integer translates by $p,q$ stay in $H$, and convex interpolation gives arbitrary nonnegative real translates. Consequently
\[
 H=T_I+K,\qquad T_I=\conv(T\cap\Z^2).
\]
A planar integer-hull algorithm computes $T_I$ in polynomial time \cite{H,CHKM}; adding $K$ and retaining the bounded edges gives the chain. This also covers a segment integer hull, including $\Delta=1$.

Let $x,y$ be consecutive lattice points along one bounded edge of $H$. The triangle $\conv(0,x,y)$ contains no lattice point other than its three vertices. A supporting line of the edge separates the origin from all nonzero lattice points in $K$, and consecutiveness excludes points in the relative interior of $[x,y]$. Hence $\det(x,y)=1$. For completeness, if the determinant were larger than one, a nonzero lattice point $z=\alpha x+\beta y$ with $0\le\alpha,\beta<1$ would lie in the fundamental parallelogram. Either $z$ or $x+y-z$ would then be an additional lattice point in the triangle.

It follows that $\det(p_{i-1},u_i)=\det(p_i,u_i)=1$. The cones between consecutive lattice points on the bounded edges of $H$ cover $K$, and each such pair is a lattice basis. This proves the generation assertion.

At an internal vertex $p_i$, the identities
$\det(p_i,u_i)=\det(p_i,u_{i+1})=1$ imply
\[
 u_{i+1}-u_i=j_i p_i\qquad(j_i\in\Z).
\]
Here $p_i$ is primitive, as the determinant identity shows. All points of $H$ satisfy $\det(x,u_i)\ge1$, and $p_{i+1}$ is not on this supporting line. Hence $\det(p_{i+1},u_i)=1+\ell_{i+1}j_i>1$. Thus $j_i\ge1$.

Set $D_i=\det(p,p_i)$ and $b_i=\det(p,u_i)$. Then $b_1=1$, $D_1=\ell_1$, and inductively $b_i>0$, with
\[
 b_{i+1}=b_i+j_iD_i\ge D_i,
 \qquad
 D_{i+1}=D_i+\ell_{i+1}b_{i+1}
       \ge(\ell_{i+1}+1)D_i.
\]
Since $D_k=\Delta$, multiplication yields
\[
 \prod_{i=1}^k(\ell_i+1)
 \le\frac{\ell_1+1}{\ell_1}\Delta\le2\Delta.
\]
Each factor is at least two, and $\ell_i+2\le2(\ell_i+1)$, proving the logarithmic bound. Finally, the primitive generators of the extreme rays are computable in polynomial time from the rational description of $K$, and $T$ has a constant-size rational description of polynomial encoding. Its computed vertices, and the primitive edge directions and lattice lengths obtained by gcd calculations, have polynomial encoding. No interior edge points are enumerated.
\end{proof}

\paragraph{Formulating the cuts in a progression.}

Return to the system $x=v+Rs$, $s\ge0$, of Lemma~\ref{lem:addition}. Let $a,u\in\Z_I^d$ with $u\ne0$, and let $M\ge0$ be an integer. Consider the normals $a+tu$, $t=0,\ldots,M$, and suppose they are pairwise sign-compatible. Their values at $v$ may cross many integers, so partitioning the progression at every crossing could require too many intervals.

The cut with normal $u$ avoids this problem: when $u^Tv$ is nonintegral, it implies all progression cuts sufficiently far from the endpoints. Near either endpoint, we can represent the remaining cuts using at most two intervals on which the fractional part is affine. Equation~\eqref{eq:clearedsplit} and Lemma~\ref{lem:interval} then yield a small formulation.

The normal $u$ need not be sign-compatible with the normals in the progression. Its cut need not follow from their cuts, but it is valid for the full split closure. If $u^Tv$ is integral, no cut with normal $u$ is required, and the endpoint cuts suffice.

\begin{lemma}\label{lem:progression}
For the system $x=v+Rs$, $s\ge0$, with rational $v,R$ and set $I$ of coordinates designated integer, let $a,u\in\Z_I^d$, with $u\ne0$, and let $M\ge0$ be an integer. Suppose that, for each $j$, the function $t\mapsto(a+tu)^Tr_j$ is either nonnegative throughout $[0,M]$ or nonpositive throughout $[0,M]$. On $s\ge0$, the split cuts with normals
\[
a+tu\quad(t=0,\ldots,M),\qquad u
\]
have a rational linear extended formulation with $O(\log(M+2))$ inequalities and auxiliary variables. Its total encoding length and construction time are polynomial in the encoding of $v,R,a,u,M$.
\end{lemma}

\begin{proof}
\emph{Remove the middle of the progression.}
If $M=0$, there are at most two cuts. Suppose $M\ge1$, and define $g=\{u^Tv\}$. If $g=0$, all progression normals have the same fractional part at $v$. When this value is nonzero, their coefficient vectors are affine in $t$, because the signs of the scalar products are fixed. The endpoint cuts suffice. When the common value is zero, no cut is required.

Suppose $0<g<1$, and set
\[
\delta=\min\{g,1-g\},\qquad h=\lceil1/\delta\rceil,
\qquad \beta=\alpha(u).
\]
For each ray, $\beta_j\le |u^Tr_j|/\delta$. Also $z_j(t)=(a+tu)^Tr_j$ has one sign throughout $[0,M]$, so $|z_j(t)|$ is nonnegative and affine there. According to whether it increases or decreases, it is at least $t|u^Tr_j|$ or $(M-t)|u^Tr_j|$. The same bound is immediate when $u^Tr_j=0$. Hence
\[
|z_j(t)|\ge\min\{t,M-t\}|u^Tr_j|.
\]
For every integer $h\le t\le M-h$ with nonzero fractional part at $v$,
\[
\alpha_j(a+tu)\ge |z_j(t)|\ge |u^Tr_j|/\delta\ge\beta_j.
\]
Thus $\beta^Ts\ge1$, the cut with normal $u$, implies all cuts in this middle interval. It remains to retain the two integer intervals
\begin{equation}\label{eq:endwindows}
[0,\min\{M,h-1\}],\qquad
[\max\{0,M-h+1\},M].
\end{equation}
They may overlap, which has no effect on the intersection.

\emph{Formulate the cuts on the two end intervals.} Set $\gamma=g$ when $g\le1/2$, and $\gamma=g-1$ otherwise.
For integer $t$,
\[
\{(a+tu)^Tv\}=\{a^Tv+\gamma t\},\qquad |\gamma|=\delta.
\]
Each interval in~\eqref{eq:endwindows} has diameter at most $h-1<1/\delta$. The affine expression inside the braces therefore varies by less than one. After removing any index with fractional part zero, the interval is the union of at most two consecutive integer intervals on which the fractional part $f(t)$ is affine and in $(0,1)$. Rational floor and ceiling operations find these intervals without enumerating their integer indices.

On one interval, let
\[
U(t,s)=\sum_{z_j\ge0}z_j(t)s_j,\qquad
V(t,s)=\sum_{z_j<0}-z_j(t)s_j,
\]
using the fixed signs of the scalar products on this interval. Identity~\eqref{eq:clearedsplit} makes the cut a quadratic inequality in $t$ with coefficients affine in $s$. Lemma~\ref{lem:interval} represents all integer indices in the interval with $O(\log(M+2))$ inequalities and variables. There are at most four intervals. Add the cut with normal $u$ to these systems.

The projection of the constructed system is the intersection of the retained progression cuts and the cut with normal $u$. Every omitted progression cut follows from the cut with normal $u$, so this projection equals the intersection of all cuts in the stated family. The quantities $h$, the interval endpoints, and all coefficients have polynomial bit length. The explicit construction in Lemma~\ref{lem:interval} proves the remaining claims.
\end{proof}

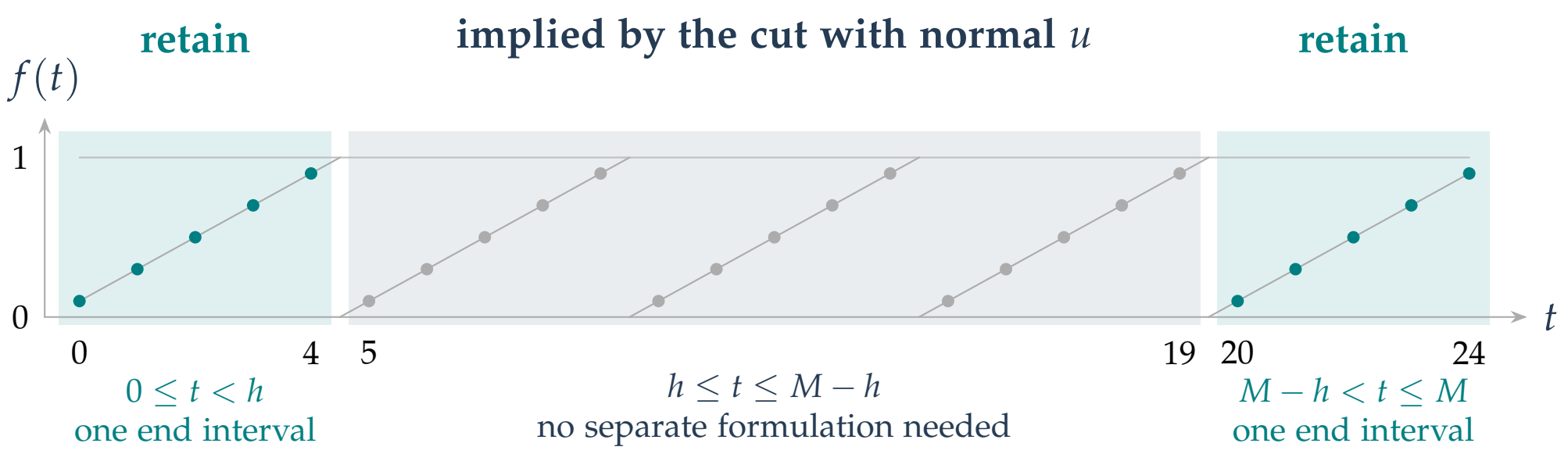
\begin{figure}[t]
\centering
\begin{adjustbox}{max width=\linewidth}
\begin{tikzpicture}[x=.49cm,y=1.35cm,font=\small]
\fill[cutteal!12] (-.35,-.05) rectangle (4.35,1.16);
\fill[lightgray] (4.65,-.05) rectangle (19.35,1.16);
\fill[cutteal!12] (19.65,-.05) rectangle (24.35,1.16);
\node[text=cutteal,font=\small\bfseries] at (2,1.75) {retain};
\node[text=inkblue,font=\small\bfseries] at (12,1.75) {implied by the cut with normal $u$};
\node[text=cutteal,font=\small\bfseries] at (22,1.75) {retain};
\draw[axis] (-.6,0)--(25,0) node[right] {$t$};
\draw[axis] (-.6,0)--(-.6,1.25) node[above] {$f(t)$};
\node[left,font=\scriptsize] at (-.6,1) {$1$};
\node[left,font=\scriptsize] at (-.6,0) {$0$};
\draw[gray!50,thin] (0,1)--(24,1);
\draw[gray!65,thin] (0.00000000,0.10000000)--(4.50000000,1.00000000);
\draw[gray!65,thin] (4.50000000,0.00000000)--(9.50000000,1.00000000);
\draw[gray!65,thin] (9.50000000,0.00000000)--(14.50000000,1.00000000);
\draw[gray!65,thin] (14.50000000,0.00000000)--(19.50000000,1.00000000);
\draw[gray!65,thin] (19.50000000,0.00000000)--(24.00000000,0.90000000);
\fill[cutteal] (0.00000000,0.10000000) circle(1.5pt);
\fill[cutteal] (1.00000000,0.30000000) circle(1.5pt);
\fill[cutteal] (2.00000000,0.50000000) circle(1.5pt);
\fill[cutteal] (3.00000000,0.70000000) circle(1.5pt);
\fill[cutteal] (4.00000000,0.90000000) circle(1.5pt);
\fill[gray!65] (5.00000000,0.10000000) circle(1.5pt);
\fill[gray!65] (6.00000000,0.30000000) circle(1.5pt);
\fill[gray!65] (7.00000000,0.50000000) circle(1.5pt);
\fill[gray!65] (8.00000000,0.70000000) circle(1.5pt);
\fill[gray!65] (9.00000000,0.90000000) circle(1.5pt);
\fill[gray!65] (10.00000000,0.10000000) circle(1.5pt);
\fill[gray!65] (11.00000000,0.30000000) circle(1.5pt);
\fill[gray!65] (12.00000000,0.50000000) circle(1.5pt);
\fill[gray!65] (13.00000000,0.70000000) circle(1.5pt);
\fill[gray!65] (14.00000000,0.90000000) circle(1.5pt);
\fill[gray!65] (15.00000000,0.10000000) circle(1.5pt);
\fill[gray!65] (16.00000000,0.30000000) circle(1.5pt);
\fill[gray!65] (17.00000000,0.50000000) circle(1.5pt);
\fill[gray!65] (18.00000000,0.70000000) circle(1.5pt);
\fill[gray!65] (19.00000000,0.90000000) circle(1.5pt);
\fill[cutteal] (20.00000000,0.10000000) circle(1.5pt);
\fill[cutteal] (21.00000000,0.30000000) circle(1.5pt);
\fill[cutteal] (22.00000000,0.50000000) circle(1.5pt);
\fill[cutteal] (23.00000000,0.70000000) circle(1.5pt);
\fill[cutteal] (24.00000000,0.90000000) circle(1.5pt);
\node[below,font=\scriptsize] at (0.00000000,-0.05000000) {$0$};
\node[below,font=\scriptsize] at (4.00000000,-0.05000000) {$4$};
\node[below,font=\scriptsize] at (5.00000000,-0.05000000) {$5$};
\node[below,font=\scriptsize] at (19.00000000,-0.05000000) {$19$};
\node[below,font=\scriptsize] at (20.00000000,-0.05000000) {$20$};
\node[below,font=\scriptsize] at (24.00000000,-0.05000000) {$24$};
\node[align=center,font=\scriptsize,text=cutteal] at (2,-.58) {$0\le t<h$\\one end interval};
\node[align=center,font=\scriptsize,text=inkblue] at (12,-.58) {$h\le t\le M-h$\\no separate formulation needed};
\node[align=center,font=\scriptsize,text=cutteal] at (22,-.58) {$M-h<t\le M$\\one end interval};
\end{tikzpicture}
\end{adjustbox}
\caption{Retaining the two end intervals of a progression of split normals. Here $M=24$, $g=1/5$, $h=\lceil1/g\rceil=5$, and $f(t)=\{1/10+t/5\}$; only the marked integer indices correspond to cuts. The cut with normal $u$ implies all cuts with $5\le t\le19$, leaving the intervals $0\le t\le4$ and $20\le t\le24$ to be represented separately.}
\label{fig:windows}
\end{figure}

For example, if $g=4/5$, use $\gamma=-1/5$: the integer-indexed fractions $\{9/10+4t/5\}$ equal $\{9/10-t/5\}$. Here $h=5$, so each end interval has diameter four and the affine function $9/10-t/5$ varies by only $4/5$. Thus the same construction applies when $\gamma<0$.

\paragraph{Total formulation size.}

Lemma~\ref{lem:chain} bounds the sum $\sum_i\log_2(\ell_i+2)$ for a cone by $O(\log(2\Delta))$. To sum over all sectors, we use the fact that each boundary ray belongs to just two of them.

\begin{lemma}\label{lem:arrangementbudget}
Let $r_1,\ldots,r_n\in\Q^2$ span $\R^2$. For each nonzero column choose a primitive integer vector $a_j$ parallel to $r_j$, and let $v_1,\ldots,v_t$ be the primitive integer generators of the distinct rays on the lines $r_j^T\pi=0$, in counterclockwise order. Put $v_{t+1}=v_1$ and $\Delta_i=\det(v_i,v_{i+1})$. The sectors $K_i=\cone(v_i,v_{i+1})$ can be computed in polynomial time. They are pointed, $t\le2n$, and
\begin{equation}\label{eq:arrangementbudget}
\sum_{i=1}^t\log_2(2\Delta_i)
=O\left(n+\sum_{j:r_j\ne0}\log_2(2+\|a_j\|_\infty)\right).
\end{equation}

\end{lemma}

Thus the sum in~\eqref{eq:arrangementbudget} is at most a constant times the total binary encoding length of the columns.

\begin{proof}
Each boundary ray belongs to exactly two neighboring sectors. Thus each ray norm appears twice when we bound the determinants by the products of their endpoint norms. Since there are at least two distinct lines, consecutive boundary rays make an angle strictly between zero and $\pi$, and $\Delta_i$ is a positive integer. The determinant estimate
\[
\Delta_i\le2\|v_i\|_\infty\|v_{i+1}\|_\infty
\]
gives
\[
\sum_i\log_2(2\Delta_i)
\le 2t+2\sum_i\log_2\|v_i\|_\infty.
\]
Each distinct input line contributes its two opposite primitive rays. Their norms equal the norm of a primitive integer vector parallel to any column defining that line. Choose one such column per line to bound the last sum by twice the sum over all nonzero columns. This proves~\eqref{eq:arrangementbudget}. Clearing the two denominators of each column separately and dividing by a greatest common divisor computes $a_j$ with encoding length linear in that column's encoding. Half-plane tests and signs of two-by-two determinants sort the boundary rays exactly in polynomial time.
\end{proof}

\begin{proof}[Proof of Theorem~\ref{thm:uniformcorner}]
Apply formula~\eqref{eq:cut} to $Q(f,R)$, as justified for dependent columns in Section~\ref{sec:normals}. It remains to choose enough normals and formulate their cuts.

We may suppose that $R$ has rank two. To ensure this, append columns $e_1,e_2$ and nonnegative variables $t_1,t_2$, and ultimately impose $t_1=t_2=0$. The coefficient of every original $s_j$ in every split cut is unchanged. Thus, by the displayed cut formula, restricting the enlarged split closure to $t=0$ recovers exactly the original split closure. This adds only a constant number of entries and variables and handles every rank-deficient case. The same device of adding rays and restricting their coordinates to zero is used for corner closures in~\cite[proof of Theorem~1.3]{BHK}.

\emph{Combine the formulations over all sectors.}
Construct the closed sectors $K_i$ of Lemma~\ref{lem:arrangementbudget}. They cover $\R^2$, and the normals in each sector are pairwise sign-compatible. For each $i$, let $H_i=\conv((K_i\cap\Z^2)\setminus\{0\})$. Lemma~\ref{lem:chain} shows that the lattice points on the bounded edges of $H_i$ generate every integer normal in $K_i$. Lemma~\ref{lem:addition} then shows that their cuts with nonzero fractional part at $f$ suffice, with integer values at $f$ handled as in Lemma~\ref{lem:normalbound}.

For each bounded edge of $H_i$ of lattice length $\ell$, use Lemma~\ref{lem:progression} to formulate its cuts, including the cut whose normal is its primitive direction $u$. Lemma~\ref{lem:chain} bounds the total size for sector $K_i$ by $O(\log(2\Delta_i))$.

Use separate auxiliary variables for each formulation, add $s\ge0$, and impose $x=f+Rs$ (and $t=0$ if rays were appended). Every point of the split closure satisfies all retained progression cuts and all cuts with normals $u$ and therefore extends to a feasible point of each formulation. Conversely, the formulations imply all cuts from lattice points on the bounded edges of the $H_i$ and then every split cut. The projection is exactly the split closure. Lemma~\ref{lem:arrangementbudget} bounds the sum of the sizes of these formulations by $O(L)$, with an absolute constant. Notice that the number of inequalities is controlled by the ray directions; the encoding of $f$ affects coefficient lengths and running time, not this sum.

\emph{Construction and encoding.}
There are $O(n)$ sectors. Each is a two-dimensional rational cone with input encoding polynomial in $L$, so Lemma~\ref{lem:chain} computes the vertices of $H_i$ and the primitive directions and lattice lengths of its bounded edges in polynomial time without enumerating lattice points in the relative interiors of the edges. These data have polynomial bit length. The progression construction uses rational floor and ceiling operations and the explicit binary-halving formulation in Lemma~\ref{lem:interval}; it therefore has polynomial total encoding and construction time. Each affine coefficient can involve all $n$ variables, which increases the encoding length polynomially but does not multiply the number of inequalities or auxiliary variables by $n$.

\emph{Sharpness.}
Take $f=(1/2,\eps)^T$ and columns $(N,1)^T$, $(-N,1)^T$, $(-N,-1)^T$, with $\eps=1/(100N^2)$ as in~\eqref{eq:construction}. The system $Q(f,R)$ for these data, with three continuous variables, is affinely isomorphic to $C_N$ with integer coordinates $I=\{1,2\}$. The isomorphism preserves those coordinates. Its split closure has at least $2N-1$ facets by the construction in Section~\ref{sec:cones}, and its input length is $\Theta(\log N)$. The face-counting argument used for Theorem~\ref{thm:ef} forces every extension to have at least $\log_2(2N-1)=\Omega(L)$ inequalities.
\end{proof}

\subsection{General polyhedra}

For the cone~\eqref{eq:generalcone}, suppose the lattice $W\cap\Z^d$ has rank at most two. Expressing its normals in a lattice basis identifies the split inequalities in the variables $s$ with those of a system having at most two integer variables. Theorem~\ref{thm:uniformcorner} then supplies a formulation; the original coordinates and lineality directions are restored by linear equations.

\begin{proposition}\label{prop:ranktwo}
Suppose the cone~\eqref{eq:generalcone} is supplied by rational data of total encoding length $L_C$, including a rational basis of its lineality space. If $\operatorname{rank}(W\cap\Z^d)\le2$, its split closure has a rational linear extended formulation with polynomial total encoding length, constructible in polynomial time, without fixing $d$. For fixed $d$, the formulation uses $O_d(L_C)$ inequalities and auxiliary variables.
\end{proposition}

\begin{proof}
Compute an integer basis matrix $B$ for $W\cap\Z^d$ by polynomial-time Hermite normal form computations~\cite{KB}. Let $k\le2$ be its number of columns. Permitted normals outside $W$ leave the cone unchanged, because they vary along its lineality space. Every remaining permitted normal is uniquely $B\mu$, with $\mu\in\Z^k$, and its value at $v$ and scalar products are
\[
\mu^T(B^Tv),\qquad \mu^T(B^TR).
\]
Thus its cuts in $s$ are exactly those of the system $x'=B^Tv+B^TRs$, $s\ge0$, where all $k$ coordinates of $x'$ are designated integer. For $k=0$ there are no cuts; for $k=1$ append an identically zero second row and a zero second entry of $B^Tv$. Theorem~\ref{thm:uniformcorner} gives the required formulation. Restore the original coordinates by $x=v+Rs+\ell$, with $\ell$ expressed in the supplied lineality basis.

All normal-form computations and transformed data have polynomial encoding length in $L_C$, uniformly in $d$. For fixed $d$, there are only a constant number of entries and every transformed entry has $O_d(L_C)$ bits. The linear bound in Theorem~\ref{thm:uniformcorner} then yields $O_d(L_C)$ inequalities and auxiliary variables, including the original cone representation.
\end{proof}

We can now apply the rank-two cone formulation to each relaxation in~\eqref{eq:corneridentity}. Counting these relaxations gives both bounds in Theorem~\ref{thm:twointeger}.

\begin{proof}[Proof of Theorem~\ref{thm:twointeger}]
Write $P=\{x:Ax\ge b\}$, where $A$ has $m$ rows, and set $r=\operatorname{rank}(A)$. Linear programming detects emptiness, in which case a contradictory inequality is an extended formulation. If $I$ is empty, use the input system itself. A nonempty rank-zero system defines $\R^d$ and needs no further cuts. Suppose therefore that $r\ge1$.

Use the relaxations $C_B$ and identity~\eqref{eq:corneridentity} from Section~\ref{sec:bounded}. Each relaxation has the form
\[
C_B=v_B+R_B\R_+^r+\ker A_B,
\qquad A_BR_B=I_r,
\qquad v_B=R_Bb_B,
\]
with polynomially encoded rational data. Its relevant normal lattice has rank at most $|I|\le2$. Apply Proposition~\ref{prop:ranktwo} to each $C_B$, using the coordinates $s=A_Bx-b_B$. Intersect the resulting systems, using separate auxiliary variables for each. The projection equals~\eqref{eq:corneridentity}. Gaussian elimination computes each right inverse $R_B$ and a basis of $\ker A_B$ with polynomial encoding and running time, uniformly in $d$. Proposition~\ref{prop:ranktwo} therefore gives an absolute polynomial bound for each system. Testing the $\binom mr$ row sets for independence and combining the systems costs $\binom mr\operatorname{poly}(L)$ in time and total encoding.

For fixed $d$, we have $r\le d$, so the number of systems is polynomial. When $m-r\le k$ for a fixed $k$, it is at most $m^k$, even when $d$ varies. A full-dimensional simplex has $m=d+1$ and $r=d$, giving $d+1$ systems. These conclusions include lineality and lower-dimensional inputs, because neither identity~\eqref{eq:corneridentity} nor the construction required full-dimensional $P$.
\end{proof}

The cones $C_N$ from~\eqref{eq:construction}, with only $x_1,x_2$ designated integer, show that the linear dependence on $L_C$ in Proposition~\ref{prop:ranktwo} is necessary already in ambient dimension three. Their input length is $\Theta(\log N)$, while the facet count and the face-counting argument used in Theorem~\ref{thm:ef} require $\Omega(\log N)$ inequalities in every extended formulation of their split closures.

For $Q(f,R)$, the \emph{triangle closure} is the intersection, over all lattice-free triangles $T\subseteq\R^2$, of the convex hulls of the points $(x,s)\in Q(f,R)$ with $x\notin\operatorname{int}T$. Basu, Hildebrand, and K\"oppe~\cite{BHK} give polynomial facet bounds for this closure and for the mixed-integer hull. Both sets are contained in the split closure, but their small descriptions do not by themselves give a formulation of that larger set. Theorem~\ref{thm:uniformcorner} represents the split closure itself with $O(L)$ inequalities and continuous auxiliary variables.

Dash, Dey, and G\"unl\"uk~\cite{DDG11} characterize the mixed-integer hull with two integer variables using a class of disjunctions more general than splits. Their characterization does not identify the split closure.

\section{Split proofs for lattice-free cuts}\label{sec:cornersimulation}

We now prove Theorem~\ref{thm:planarsimulation}. Let $P$ be a rational polytope with all variables designated integer, let $u:P\to\R^2$ be the restriction of an affine map with integer coefficients, and let $h$ be rational affine. Suppose that $h(z)>0$ implies $u(z)\in\operatorname{int}K$, where $K\subseteq\R^2$ is a supplied full-dimensional rational lattice-free polyhedron. Then $h\le0$ is valid for the integer points of $P$. Our goal is to derive this inequality by a polynomially encoded sequence of split cuts that retains every original point with $h\le0$.

For example, designate all three variables as integer in
\[
P=\{(p,q,h):0\le h\le1,\ p\ge h,\ q\ge h,\ p+q\le2-h\}.
\]
Every point with $h>0$ has $(p,q)$ in the triangle $p>0$, $q>0$, $p+q<2$. On the branch $p\le0$ we have $h=0$; on $p\ge1$, the last two inequalities give $h\le1/2$. Thus the split $p\le0$ or $p\ge1$ permits the cut $h\le1/2$. The split $h\le0$ or $h\ge1$ then permits $h\le0$. The general argument likewise decreases an upper bound until integer rounding applies. We must control the number and encoding of these cuts while retaining every point with $h\le0$.

We first prove a fixed-factor decrease for a triangle in a specified position, and then establish exact preservation for cross cuts, which we use to treat general triangles. Enlarging $K$ to a maximal lattice-free set reduces the remaining cases to a split strip or a quadrilateral. Throughout, $h$ is a rational affine function of the original variables; in general it is not a separate coordinate, as it was in the example. For an affine function $p$ with integer coefficients, a ``split on $p$'' below means $p\le0$ or $p\ge1$.

The dimension of $P$ need not be fixed. Integrality of all its original variables permits the final rounding; integrality of the two components of $u$ alone is insufficient. This distinguishes the theorem from the mixed-integer finite-rank results of Dey and Louveaux~\cite{DL11}. For the applications to the two-row system $Q(f,R)$, we also use the polynomial facet enumeration of Basu, Hildebrand, and K\"oppe~\cite[Theorems~1.4--1.5]{BHK}.

\paragraph{Keeping earlier cuts.}
Several parts of the proof construct a sequence of cuts on the original polytope, or on a polytope containing the current relaxation. These cuts can be added while retaining all earlier cuts. Indeed, suppose $Q\subseteq R$. Apply to $Q$ a split proof constructed on $R$, and let $Q_j,R_j$ be the sets obtained after imposing its first $j$ inequalities. Inductively, $Q_j\subseteq R_j$. The containment holds on both branches of the next split, so the inequality is also a split cut for $Q_j$. If the split proof has $m$ cuts, the final set on $Q$ is contained in $Q\cap R_m$ and contains every integer point of $Q$.

To control encoding length, construct each sequence on its specified starting polytope rather than recomputing it on the accumulated system. Polynomially many sequences have polynomial total encoding if their construction bounds are polynomial in the total encoding of those starting polytopes and supplied data. We use this observation throughout Sections~\ref{sec:cornersimulation}--\ref{sec:fixedproof}.

\subsection{Lattice-free triangles}

If $h>0$ implies that the image lies in the interior of a triangle, each inequality defining the triangle is strict at such points. We need quantitative versions of this fact: a positive lower bound proportional to $h$, and an affine bound valid also where $h\le0$. The next two lemmas provide them.

The first lemma can preserve a supplied polytope $G$ in addition to the integer points. It also applies with continuous variables, provided that the function defining the split has zero coefficients on those variables. The same qualification applies to the fixed-factor decrease proved below; exact termination uses integrality of all original variables.

\begin{lemma}\label{lem:supportnormalize}
Let $G\subseteq S$ be rational polytopes, let $g$ be rational affine with $g\le U$ on $S$, where $U>0$ is rational, and suppose $g\le0$ on $G$. Let $\ell$ be affine with integer coefficients and suppose
\[
 g>0\ \Longrightarrow\ \ell>0\qquad\text{on }S.
\]
At most one split cut, preserving $G$, ensures
\[
 g>0\ \Longrightarrow\ \ell\ge g/U.
\]
The cut can be found by rational LP with polynomial encoding.
\end{lemma}
\begin{proof}
The candidate cut is $F\ge0$, where $F=\ell-g/U$. We show that either it is valid on both branches $\ell\le0$, $\ell\ge1$ and on $G$, or the desired bound already holds wherever $g>0$. Suppose the validity test fails, and take a point $z^-$ in a branch or in $G$ at which $F<0$. The upper branch cannot contain it, since $\ell\ge1$ and $g\le U$ there. On the lower branch the hypothesis gives $g\le0$, and the same holds on $G$. Hence $g(z^-)\le0$ and $F(z^-)<0$.

If a point $z^+$ with $g(z^+)>0$ also had $F(z^+)<0$, the segment $[z^-,z^+]$ would have $F<0$. At its crossing of $g=0$ we would have $\ell<0$, and this would persist at nearby points with $g>0$, contradicting the hypothesis. Thus failure of the validity test means that no cut is needed.

Rational LP tests $F\ge0$ on the two branches and on $G$. If all tests pass, add it as a split cut preserving $G$; otherwise no cut is needed. If $\ell$ is a constant integer, one branch is all of $S$, so any inequality that passes these tests already holds on $S$. The tests use the supplied rational data and have polynomial encoding.
\end{proof}

The same argument strengthens an inequality valid throughout $S$. If initially $\ell\ge\beta g$ on $S$ for some $\beta>0$, the same procedure ensures $\ell\ge\beta' g$ globally for some $\beta' U\ge1$, without computing $\beta$. Indeed, if $\beta\ge1/U$, the original inequality already suffices. If $\beta<1/U$, the original bound implies $\ell\ge g/U$ wherever $g\le0$, and the lemma ensures it wherever $g>0$. In the latter case take $\beta'=1/U$.

\begin{lemma}\label{lem:slackcertificate}
Let $P$ be a rational polytope and let $h,g_1,\ldots,g_m$ be rational affine. If $h\le0$ on $P\cap\{g_1\le0,\ldots,g_m\le0\}$, one can compute nonnegative rational multipliers $\lambda_1,\ldots,\lambda_m$ of polynomial encoding such that
\[
 h\le\sum_i\lambda_i g_i\qquad\text{on }P.
\]
\end{lemma}
\begin{proof}
If $P$ is empty, take all multipliers to be zero. Otherwise, if the constrained set is nonempty, LP duality gives the inequality after dropping a nonpositive constant. If it is empty, an infeasibility certificate gives $\sum_i\mu_i g_i\ge1$ on $P$, with $\mu_i\ge0$. Multiplying by $\max\{0,\max_P h\}$ gives the required bound on $h$. Rational LP supplies all certificates and values with polynomial encoding.
\end{proof}

In particular, if $\max_P h>0$ and $h>0$ implies $\ell>0$, the lemma gives $h\le\lambda\ell$ with $\lambda>0$. Thus $\ell\ge\gamma h$ globally for some positive rational $\gamma$, including where $h<0$.

\paragraph{Rotating the sides.}

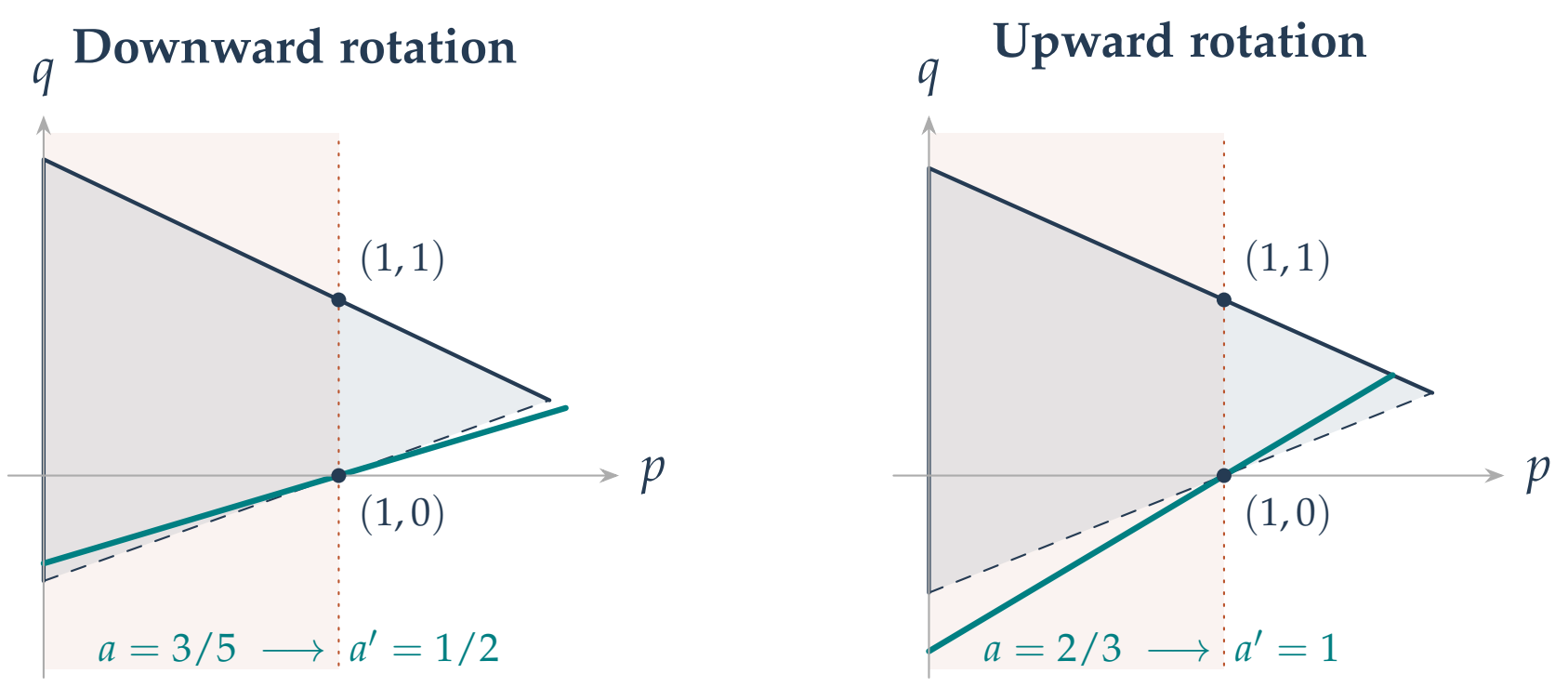
\begin{figure}[htbp]
\centering
\begin{adjustbox}{max width=\linewidth}
\begin{tikzpicture}[x=2.1cm,y=1.25cm,font=\small]
\foreach \shift/\a/\c/\newa/\oldlabel/\newlabel/\ttl in {0/.6/.8/.5/{3/5}/{1/2}/{Downward rotation},6.3/.666666667/.75/1/{2/3}/{1}/{Upward rotation}} {
\begin{scope}[xshift=\shift cm]
\pgfmathsetmacro{\tipx}{1+1/(\a+\c)}
\pgfmathsetmacro{\tipy}{\a/(\a+\c)}
\pgfmathsetmacro{\newx}{1+1/(\newa+\c)}
\pgfmathsetmacro{\newy}{\newa/(\newa+\c)}
\fill[lightgray] (0,-\a)--(0,{1+\c})--(\tipx,\tipy)--cycle;
\fill[warmrust,opacity=.07] (0,-1.1) rectangle (1,1.95);
\draw[boundary] (0,-\a)--(0,{1+\c})--(\tipx,\tipy);
\draw[inkblue,dashed] (0,-\a)--(\tipx,\tipy);
\draw[cut] (0,-\newa)--(\newx,\newy);
\draw[axis] (-.12,0)--(1.95,0) node[right] {$p$};
\draw[axis] (0,-1.15)--(0,2.05) node[above] {$q$};
\draw[warmrust,dotted] (1,-1.08)--(1,1.95);
\fill[inkblue] (1,0) circle(1.5pt) node[below right,font=\scriptsize] {$(1,0)$};
\fill[inkblue] (1,1) circle(1.5pt) node[above right,font=\scriptsize] {$(1,1)$};
\node[text=inkblue,font=\small\bfseries] at (.85,2.45) {\ttl};
\node[text=cutteal,font=\scriptsize,anchor=west] at (.12,-.98) {$a=\oldlabel\ \longrightarrow\ a'=\newlabel$};
\end{scope}}
\end{tikzpicture}
\end{adjustbox}
\caption{Change the lower slope from $a$ to $a'$, keeping $(1,0)$ fixed. The upper line remains fixed through $(1,1)$. The left picture decreases $a$; the right increases it. The shaded strip is $0<p<1$; $h$ is not a coordinate in the picture.}
\label{fig:triangle-rotation}
\end{figure}

Consider the inequalities $p>0$, $q>a(p-1)$, and $1-q>c(p-1)$ in the $(p,q)$ plane, where $a,c\ge0$. When $a=c=0$, the region lies in the split strip $0<q<1$. The lower and upper boundary lines pass through $(1,0)$ and $(1,1)$, respectively. If $a,c$ are integers, Lemma~\ref{lem:supportnormalize} applied to the three inequalities lets us decrease the bound on $h$ with a split on $p$. For rational $a,c$, we change one boundary line at a time, keeping its specified point fixed, until the slopes are integers. Lemma~\ref{lem:bezouthalving} reduces both the numerator and denominator of each changed slope by at least a factor of two, while our bounds ensure that the final split still decreases the upper bound on $h$ by a fixed fraction. Figure~\ref{fig:triangle-rotation} shows a decrease and an increase of the lower slope $a$.

During these changes, the proof maintains lower bounds on the left-hand sides obtained by clearing denominators in the three defining inequalities. Each bound is proportional to $h-t$, where $t$ increases during the construction. Increasing $t$ preserves the bounds for the unchanged sides, so only the new side can require another cut.

\begin{lemma}\label{lem:rationaltrianglecompiler}
Let $Q$ be a rational polytope, let $p,q$ be affine with integer coefficients, and let $h$ be rational affine with $h\le H$ on $Q$, for a supplied rational $H>0$. Suppose that, on $Q$, the supplied rational numbers $a,c\ge0$ satisfy
\[
 h>0\ \Longrightarrow\ p>0,\quad q>a(p-1),\quad 1-q>c(p-1).
\]
One can construct a split proof of $h\le\rho H$ with the absolute constant $\rho=511/512$. Construction time and total encoding are polynomial. If all original variables are integer, a polynomial split proof establishes $h\le0$. Every original point with $h\le0$ is preserved.
\end{lemma}

\begin{proof}
We change the rational slopes until both are integers. At that point, adding the lower bounds associated with the two sides will bound $h$ on the branch $p\ge1$, while the bound associated with the base will handle $p\le0$. We first describe the inequalities maintained during these changes.

\emph{Maintain lower bounds for the three sides.}
If $Q$ is empty or $\max_Qh\le0$, stop. Retain $G=Q\cap\{h\le0\}$ throughout. Write the current slopes as $a=A/D$, $c=C/E$ in lowest terms, with $A,C\ge0$ and $D,E\ge1$ integers. For a rational number $t$ with $0\le t<H$, let
\begin{equation}\label{eq:trianglecertificate}
 U=H-t,\qquad \lambda=(h-t)/U.
\end{equation}
At every point of the current relaxation with $h>t$, we maintain
\begin{equation}\label{eq:trianglenormalization}
 p\ge\lambda,\qquad F:=Dq-A(p-1)\ge\lambda,\qquad
 J:=E(1-q)-C(p-1)\ge\lambda.
\end{equation}
At $t=0$, apply Lemma~\ref{lem:supportnormalize} to the functions $p,F,J$, which have integer coefficients. At most three cuts establish these bounds while preserving $G$.

\emph{Rotate one side.}
Only nonintegral slopes need to change. For a nonintegral lower slope, the line $F=1$ meets $q=1$ at $p=1+(D-1)/A$. For a nonintegral upper slope, $J=1$ meets $q=0$ at $p=1+(E-1)/C$. If both slopes are nonintegral, choose the one with the smaller ratio $A/(D-1)$ or $C/(E-1)$. If only one is nonintegral, choose that one. Suppose we choose $a$. Then $D\ge2$ and
\[
 A(E-1)\le C(D-1),
\]
which also holds when $c$ is integral, since then $E=1$. Lemma~\ref{lem:bezouthalving} supplies a new slope $a'=A'/D'$ with
\[
 0\le A'\le A/2,\qquad 1\le D'\le D/2,\qquad
 |AD'-A'D|=1.
\]
Let $\varepsilon=AD'-A'D$ and $F'=D'q-A'(p-1)$. Set
\[
 \theta=\begin{cases}
 1/(m+1),\quad m=\min\{A',D'\},&A'>0,\\
 1/2,&A'=0.
 \end{cases}
\]
We show directly that $F'>0$ whenever $\lambda>\theta$.

If $\varepsilon=1$, the bounds on $p$ and $F$ give
\[
 DF'=D'F+p-1\ge(D'+1)\lambda-1>0.
\]
This also covers $A'=0$: then $AD'=1$, so $A=D'=1$ and $\theta=1/2$.

If $\varepsilon=-1$, then $A'>0$. Let $S=a+c>0$. Use the identity $F/D+J/E=1-S(p-1)$ to obtain
\[
 DSF'=(A'+cD')F+J/E-1.
\]
Here $AD'=A'D-1\ge D-1$, so the choice of side implies $CD'\ge E-1$. Consequently
\[
 A'+cD'+1/E\ge A'+1,
 \qquad DSF'\ge(A'+1)\lambda-1>0.
\]
The first calculation uses the base and lower-side inequalities; the second uses both side inequalities. To change $c$, interchange $q$ and $1-q$ and apply the same argument.

Increase $t$ to $t'=t+\theta U$. At points with $h>t'$, $\lambda'=(h-t')/(H-t')=(\lambda-\theta)/(1-\theta)\le\lambda$. The bounds for the unchanged sides therefore remain valid. Apply Lemma~\ref{lem:supportnormalize} with $g=h-t'$, upper bound $H-t'$, and function $F'$, which has integer coefficients. At most one cut establishes the bound for the new side and restores~\eqref{eq:trianglenormalization}. Since $h-t'\le0$ on $G$, this cut preserves $G$.

\emph{Bound the total increase in $t$.}
Along either side, the successive positive values of $m$ decrease by at least a factor of two, so their reciprocals sum to at most two. Each such change multiplies $U$ by $m/(m+1)$. Since $\log(1+1/m)\le1/m$, the product of these factors over both sides is at least $e^{-4}$. A final change to slope zero, if needed, contributes one further factor of $1/2$ per side. Thus the final value $U_*=H-t$ satisfies
\[
 U_*\ge\frac{H}{4e^4}>\frac{H}{256}.
\]
Both slopes are now integers, so $D=E=1$. On $p\le0$, the maintained bounds exclude $h>t$; on $p\ge1$ they give
\[
 2(h-t)/U_*\le F+J=1-(a+c)(p-1)\le1.
\]
Consequently the $p$ split permits $h\le H-U_*/2\le\rho H$, preserving $G$. If $p$ is constant, its sole branch already proves the bound and no split is needed.

\emph{Encoding and exact termination.}
Let $L$ be the input length. Denominators halve, so one fixed-factor decrease requires $O(L)$ rotations. All slope numerators and denominators decrease. The rational factors $1-\theta$ have $O(L)$ bits, giving $O(L^2)$ bits for $t,U$ and the inequalities added by Lemma~\ref{lem:supportnormalize}. Their validity tests are rational LPs on the original data and retained cuts.

Choose a positive integer $M$ clearing the denominators of $h$, including its constant term. Repeat the construction from the original slopes and $t=0$ with bounds $H_j=\rho^jH$, $j=0,1,\ldots$, retaining all cuts and the same $G$. After $O(1+\log(1+MH))=O(L)$ repetitions, $Mh<1$. Since $\rho$ is fixed, these $H_j$ have $O(L)$ bits; the entire proof has $O(L^2)$ cuts and polynomial total encoding. If all variables are integer and $h$ is nonconstant, the split $Mh\le0$ or $Mh\ge1$ gives $h\le0$. A positive constant $h$ instead forces emptiness, since $Mh$ is a positive integer; a nonpositive constant needed no cuts.
\end{proof}

\subsection{Cross cuts}

A cross cut uses two split disjunctions at once: its inequality is valid on each of their four branch intersections. We will prove it using ordinary split cuts while retaining every point that already satisfies it.

Applying Lemma~\ref{lem:slackcertificate} to the four branches gives four affine functions that are positive wherever $h>0$. We reduce their integer coefficients one pair at a time, increasing $t$ as needed to keep the new function positive wherever $h>t$. The following calculation bounds this increase and ensures that every point with $h\le0$ is retained.

\paragraph{Bounds valid on the whole relaxation.}
Let $h\le H$ on the current relaxation, where $H>0$, and set $U=H-t$, with $0\le t<H$. Suppose the current affine functions with integer coefficients satisfy
\[
 L_i\ge\alpha_i(h-t),\qquad \alpha_i>0.
\]
We can ensure $\alpha_i\ge1/U$ with at most one split cut for each function. If its current coefficient is already at least $1/U$, retain it. Otherwise the inequality $L_i\ge(h-t)/U$ is valid on both branches of the split on $L_i$: the lower branch forces $h-t\le0$, and the upper branch has $L_i\ge1$ and $h-t\le U$. At $h-t\le0$ the old bound implies the new one. Thus the cut also preserves every original point with $h\le0$. We retain a larger coefficient when one is already available, since decreasing it need not preserve the bound where $h<t$.

Suppose a new function $L$ satisfies an identity
\begin{equation}\label{eq:preservingtransfer}
 A L=\sum_i w_iL_i-B,\qquad A>0,\quad w_i\ge0,\quad B\ge0.
\end{equation}
Writing $\alpha=\sum_iw_i\alpha_i/A$ and $\beta=B/A$ gives $L\ge\alpha(h-t)-\beta$ globally. If
\[
 \theta=\frac{B}{\sum_iw_i}<1,\qquad t'=t+\theta U,
\]
then $\alpha\theta U\ge\beta$, so $L\ge\alpha(h-t')$. The unchanged functions retain their bounds relative to $h-t'$ as well. Apply the preceding split cuts as needed to ensure that each coefficient is at least $1/(H-t')$. Thus $H-t'=(1-\theta)(H-t)$; this is the quantity we must keep large enough to obtain a useful decrease in the upper bound on $h$.

Once the coefficients have been reduced, an exact split proof of $h\le t$ retains the original set $h\le0$. If $H-t\ge\delta H$, this reduces the upper bound by the factor $1-\delta$. Repeat from the original polytope with predetermined bounds $H_j=(1-\delta)^jH$, retaining all actual cuts. Each repetition constructs its cuts on the original polytope intersected with its bound $h\le H_j$, which contains the current relaxation, so the cuts remain valid on it. If $\delta^{-1}$ is polynomial in the input length, polynomially many repetitions give $Mh<1$, where $M$ clears all coefficients of $h$. The split on $Mh$ then gives $h\le0$. For positive constant $h$, the smaller upper bound already proves emptiness. These repetitions never use the accumulated cut description as the input to a new coefficient calculation.

One arithmetic case occurs in both the cross and triangle proofs, so we give its details once. All additional inequalities defining $P$ are retained in the following lemma.

\begin{lemma}\label{lem:residuepreservation}
Suppose affine functions $p,q,F,G$ with integer coefficients satisfy
\[
 F=(kb+1)p+bq\ge\gamma_F h,\qquad
 G=c(1-p)+dq\ge\gamma_G h
\]
on a rational polytope $P$ with all coordinates designated integer, where $h$ is rational affine, $k,b,c,d$ are positive integers and $\gamma_F,\gamma_G>0$ are rational. Suppose the input also supplies affine functions $L_i$ with integer coefficients and bounds $L_i\ge\gamma_i h$ on $P$, where $\gamma_i>0$. Put $U_0=kp+q$. Assume that whenever each of $U_0,G,L_i$ is bounded below throughout the current polytope by a positive multiple of $h-v$, there is a polynomially constructible split proof ending exactly at the current polytope intersected with $\{h\le v\}$. Then there is such a proof for $h\le0$. The number and total encoding of the applications of the supplied construction and the additional cuts are polynomial in the original supplied data.
\end{lemma}
\begin{proof}
We first obtain an upper bound $V_0$ on $h$, then reduce it by a factor of $4/5$ at each repetition. Set $a=kb+1$, $E=ad+bc$, and $D=1-p+q$. The identities we need are
\begin{align}
 EU_0&=(kd+c)F+G-c,\\
 ED&=(c-d)F+(a+b)G-a(c-d),\\
 (a+b)U_0&=(k+1)F+D-1.\label{eq:residueidentities}
\end{align}
The first identity gives $U_0\ge\alpha(h-V_0)$ on $P$ for some $\alpha>0$, where
\[
 V_0=\frac{c}{(kd+c)\gamma_F+\gamma_G},\qquad \gamma_F V_0<1.
\]
The bounds for $G,L_i$ remain valid with $h$ replaced by $h-V_0$. Apply the supplied construction to obtain $h\le V_0$, unless that bound already holds.

At a predetermined bound $h\le V\le V_0$, first suppose $c<d$. The split on $F$ permits
\[
 U_0\ge\left(\frac{kd+c}{EV}+\frac{\gamma_G}{E}\right)h-\frac cE.
\]
For $h\le0$, this follows from the first identity because its new coefficient is at least the old one. Positive $h$ is impossible on $F\le0$. On $F\ge1$, the identity gives $EU_0\ge kd+\gamma_G h$, which implies the proposed inequality since $h\le V$. The right-hand side vanishes at a value of $h$ below $V/2$. Thus the inequality gives $U_0\ge\alpha(h-4V/5)$ for some $\alpha>0$ throughout the relaxation. Apply the supplied exact construction at that target.

If $c\ge d$, the second identity gives $D\ge C_Dh-\beta$, where
\[
 C_D=\frac{(c-d)\gamma_F+(a+b)\gamma_G}{E},\qquad
 \beta=\frac{a(c-d)}E<k+1.
\]
Moreover $C_DV_0\le1+\beta<k+2$. Consequently the inequality
\[
 D\ge\frac{k+2}{V}h-(k+2)
\]
holds on every original point with $h\le0$. It is a split cut for the split on $p$. At $0\le h\le V$ its right-hand side is nonpositive; on $p\le0$ the bound for $F$ gives $D\ge0$, and on $p\ge1$ the bound for $G$ gives $D\ge0$ because $c\ge d$. Apply the split cut described above to ensure $F\ge A_Fh/V$ throughout the relaxation, with $A_F\ge1$. Identity~\eqref{eq:residueidentities} yields
\[
 (a+b)U_0\ge\bigl((k+1)A_F+k+2\bigr)h/V-(k+3).
\]
At $h=4V/5$ its right-hand side is at least $3(k-1)/5\ge0$. Its slope is positive, including when $k=1$, so again the supplied exact construction applies at $4V/5$.

Repeat with $V_j=(4/5)^jV_0$ and finish with the split on $Mh$. Every additional cut preserves $h\le0$, and each application of the supplied construction preserves every point satisfying its bound on $h$. Retain these additional inequalities in subsequent auxiliary polytopes. Once an exact application is complete, later steps need only the bound $h\le v$ it establishes. The number of additional inequalities and applications is polynomial, and their coefficients are explicit expressions of polynomial encoding in the original data and $V_j$. Thus the encoding bound depends polynomially on the original data, even though all cuts remain in the proof.
\end{proof}

\begin{proposition}\label{prop:crosssimulation}
Let $P\subseteq\R^n$ be a rational polytope with all variables designated integer, let $p,q$ be supplied affine functions with integer coefficients, and let $h$ be rational affine. If $h\le0$ is valid on all four branch intersections of the splits on $p$ and $q$, one can construct a split proof ending exactly at $P\cap\{h\le0\}$, in polynomial time and with polynomial total encoding in $(P,p,q,h)$. The dimension may vary.
\end{proposition}
\begin{proof}
Assume $\max_Ph>0$. Applying Lemma~\ref{lem:slackcertificate} to each of the four branch intersections gives nonnegative rational multipliers for the two branch inequalities. If both multipliers in a pair vanish, we already have $h\le0$. Otherwise clear denominators and divide each pair by its gcd. Rescaling the corresponding positive coefficient of $h$ gives inequalities $F\ge\gamma_Fh$, $G\ge\gamma_Gh$, $J\ge\gamma_Jh$, $K\ge\gamma_Kh$ on $P$, with $\gamma_F,\gamma_G,\gamma_J,\gamma_K>0$, where
\[
 F=ap+bq,\quad G=c(1-p)+dq,\quad
 J=e(1-p)+f(1-q),\quad K=gp+j(1-q).
\]
First suppose all eight coefficients are positive integers. We induct on the number $m\le4$ of pairs $(a,b),(c,d),(e,f),(g,j)$ in which both entries exceed one.

For $m=0$, first suppose at least one of the coefficient pairs $(a,c)$, $(d,f)$, $(e,g)$, $(j,b)$ equals $(1,1)$. By a symmetry of the square, we may take $a=c=1$. Then $F+G=1+(b+d)q$, and~\eqref{eq:preservingtransfer} makes $q>0$ whenever $h>t+U/2$. Increase $t$ by $U/2$. The inequalities for $J,K$ give
\[
 p>\frac jg(q-1),\qquad 1-p>\frac fe(q-1)
\]
whenever $h>t$. Lemma~\ref{lem:rationaltrianglecompiler}, with base $q$, proves $h\le t$ exactly. If none of the four pairs listed above equals $(1,1)$, square symmetry gives
\[
 F=ap+q,\quad G=1-p+dq,\quad
 J=e(1-p)+1-q,\quad K=p+j(1-q),
\]
where $a,d,e,j>1$. Choose a smallest coefficient as $a$, so $a\le d$. The identity
\[
 (ad+1)(p+q)=(d+1)F+(a-1)G-(a-1)
\]
allows us to replace $F$ by $p+q$ with $\theta\le1/2$. The coefficients of $p$ and $1-p$ in the new function and $G$ are now both one, so the preceding argument applies. A complete application reduces $H$ to at most $3H/4$. The predetermined repetitions described above prove the assertion for $m=0$.

For $m>0$, select a pair with both coefficients greater than one and use square symmetry to write it as $(a,b)$ with $a\ge b\ge2$. If $F'=rp+sq$ and $\Delta=as-br$, then
\begin{align}
 (ad+bc)F'&=(rd+sc)F+\Delta G-\Delta c &&(\Delta>0),\\
 (aj+bg)F'&=(rj+sg)F+(-\Delta)K-(-\Delta)j &&(\Delta<0).
 \label{eq:cross-transfers}
\end{align}
By~\eqref{eq:preservingtransfer}, these identities give $\theta\le\Delta/s$ and $\theta\le(-\Delta)/r$, respectively. We will choose $r/s$ so that $\theta\le1/2$ and the denominator decreases rapidly. Continued fractions provide both properties.

Consider the convergents of the regular continued fraction of $a/b$. Choose the first $r_i/s_i$ with $s_i^2\ge2b$. If it is not already $a/b$, the convergent estimate gives
\[
 |as_i-br_i|\le b/s_{i+1}\le b/s_i\le s_i/2.
\]
Since $r_i\ge s_i$, the change to this convergent has $\theta$ at most $1/2$. We then use the predecessor of this convergent. Consecutive convergents have determinant of magnitude one, and the predecessor denominator is less than $\sqrt{2b}$. If that denominator is at least two, this change also has $\theta$ at most $1/2$. 

If the predecessor denominator is one and the determinant is positive, the current pair is $(kb+1,b)$; stop before making the change and apply Lemma~\ref{lem:residuepreservation}. Its supplied exact construction uses the pair $(k,1)$ and hence the induction hypothesis for $m-1$. If the determinant is negative, the pair is $(kb-1,b)$ with $k\ge2$; $\theta$ is at most $1/k$, and the new pair has a coefficient equal to one. Apply the induction hypothesis to $h-t\le0$. 

For $2\le b\le7$, take successive predecessors $r/s$ directly. Their denominators strictly decrease and $as-br=\pm1$. When $s\ge2$, the respective values of $\theta$ are at most $1/s$ and $1/r$, both at most $1/2$. After at most five such changes, $s=1$: positive determinant gives $(kb+1,b)$ and Lemma~\ref{lem:residuepreservation}; negative determinant gives $(kb-1,b)$, $k\ge2$, and a final change with $\theta$ at most $1/k$. Thus at most six changes suffice.

While $b\ge8$, each pair of changes before the final case replaces $b$ by a number below $\sqrt{2b}$. Thus there are $O(\log\log(b+2))=O(\log L)$ changes, where $L$ is the supplied input length. If their number is $N$, the final value of $t$ satisfies $H-t\ge2^{-N}H$. An exact proof of $h\le t$ therefore decreases the bound $H$ by an inverse-polynomial fraction. The repetitions described above give $h\le0$ in polynomially many applications. The induction depth is at most four, because a recursive application is needed only when $m$ decreases.

All changed integer coefficients are bounded by the original ones. The values of $t$ and the coefficients in $L_i\ge\alpha_i(h-t)$ are obtained by rational sums and products along $O(\log L)$ changes, so their encodings are polynomial. Each outer repetition uses its original polytope and a predetermined bound. Lemma~\ref{lem:residuepreservation} retains only polynomially many explicit additional inequalities between its exact applications. At most four nested applications therefore give a polynomial bound in the original input length.

It remains to allow zero coefficients. An identically zero function among $F,G,J,K$ already proves $h\le0$. Otherwise a symmetry of the square gives $F=ap$ with $a>0$. If either $G$ or $J$ is a positive multiple of $1-p$, one split on $p$ proves the target. In the remaining case the inequalities for $G,J$ give
\[
 p>0,\qquad q>(c/d)(p-1),\qquad 1-q>(e/f)(p-1)
\]
whenever $h>0$, with $d,f>0$ and $c,e\ge0$. Lemma~\ref{lem:rationaltrianglecompiler} applies and preserves the set $h\le0$. Empty branch intersections were already covered by Lemma~\ref{lem:slackcertificate}. No step requires the linear parts of $p,q$ to be independent. If a function defining a split is an integer constant, one branch contains the whole relaxation and no cut is needed. Every construction retains the original set $h\le0$ and ultimately imposes $h\le0$, proving exact equality.
\end{proof}

\subsection{Lattice-free cuts}

Cross cuts do not cover every lattice-free triangle: a maximal triangle of type~3 need not be contained in the union of two split strips~\cite[Lemma~3.3]{DDG12}. We therefore need a separate argument for triangles. The difficulty is to change their boundary lines without losing either the points to be retained or a polynomial bound on the number of cuts. We choose coordinates at integer points on the edges, so these changes can be expressed by integer identities.

We first recall how to enlarge the supplied full-dimensional rational lattice-free polyhedron to a rational maximal lattice-free set of polynomial encoding. Write its inequalities as $a_i^Tu\le b_i$ with integral $a_i$. Delete inequalities while the system with every inequality strict still has no integer solution. On integer points, $a_i^Tu<b_i$ is equivalent to $a_i^Tu\le\lceil b_i\rceil-1$. Doignon's theorem and fixed-dimensional integer programming~\cite{Doignon73,L83} therefore leave at most four rows in polynomial time.

For each remaining row, minimize $a_i^Tu$ over the integer points satisfying every other row strictly. Delete the row if there are none; otherwise replace $b_i$ by this minimum. The objective is bounded below by the old $b_i$, since the set is lattice-free. The minimum and an attaining integer point have polynomial encoding. This point is in the relative interior of the corresponding facet. Later outward movements of other facets preserve that property. After at most four such steps every facet contains an integer point in its relative interior, so the resulting set is maximal lattice-free. It is a split strip, a triangle, or a quadrilateral~\cite{BCCZ10}. A strip requires one split cut. For a quadrilateral, the integer points in the relative interiors of its four edges form a parallelogram of area one; in unimodular coordinates its interior is contained in $\{0<p<1\}\cup\{0<q<1\}$~\cite[Section~2.2 and Lemma~3.2]{DDG12}. Proposition~\ref{prop:crosssimulation} applies.

For a triangle, choose one integer point in the relative interior of each edge. Every other point of their convex hull is interior to the lattice-free triangle. Their convex hull has area $1/2$ by Pick's theorem, so a unimodular change sends these points to $(0,0),(1,0),(0,1)$. Write the three defining inequalities as $F\ge0,G\ge0,K\ge0$, with primitive integer normal vectors, where
\begin{equation}\label{eq:contacttriangle}
 F=ap+bq,\qquad G=c(1-p)+dq,\qquad K=ep+f(1-q),
\end{equation}
with $a,b,c,f>0$, $c+d>0$, and $e+f>0$. The functions $F,G,K$ are positive whenever $h>0$. Lemma~\ref{lem:slackcertificate} gives $F\ge\gamma_Fh$, $G\ge\gamma_Gh$, and $K\ge\gamma_Kh$ on $P$, with $\gamma_F,\gamma_G,\gamma_K>0$.

Several positions reduce directly to arguments already proved. If $a=b$, primitivity gives $F=p+q$. In coordinates $r=p+q$, $s=q$, the other two inequalities are
\[
 G/(c+d)=s-A(r-1),\qquad
 K/(e+f)=1-s-C(r-1),
\]
with $A=c/(c+d)>0$ and $C=-e/(e+f)>-1$. If $e\le0$, both slopes are nonnegative and Lemma~\ref{lem:rationaltrianglecompiler} applies. This branch needs only the coefficient conditions and these three inequalities, so it also applies after changing the defining inequalities or replacing $h$ by $h-t$, with additional inequalities retained; it does not use lattice-freeness. If $e>0$, lattice-freeness at $(1,1)$ forces $d\le0$, hence $A\ge1$; replacing $s$ by $s-r+1$ gives the nonnegative slopes $A-1,C+1$.

Otherwise interchange $p,q$ and relabel the sides so that $a>b$. The integer point $(1,-1)$ forces $d\ge0$. If $d>0$, the point $(1,1)$ forces $e\le0$. When $d=0$, use $r=1-p,s=q$: the other slopes are $a/b,e/f$. They are nonnegative if $e\ge0$; if $e<0$, replacing $s$ by $s-r+1$ gives $a/b-1,1+e/f\ge0$. Again the rational-triangle lemma applies. When $e=0$, positive $h$ excludes $q\ge1$; on $q\le0$, the bounds for $F,G$ exclude both branches of the split on $p$. This is a cross cut.

It remains to treat
\begin{equation}\label{eq:signedtriangle}
 a>b>0,\qquad c,d>0,\qquad -f<e<0.
\end{equation}
Let $A=a/b$, $B=(c+d)/c$, and $C=f/(f+e)$, all greater than one. The unimodular change $p'=q$, $q'=1-p-q$ cycles these three ratios. It also shows why the same coefficient calculations can be used at any of the three chosen integer points.

Replacing $A$ by any $A'\ge1$, with the other sides fixed, leaves the triangle lattice-free. Indeed, for integer $p\ge1$, $G>0$ forces $q>0$ and $K>0$ forces $q<1$. For $p=0$ the same follows from $F>0,K>0$. For integer $p=-m<0$, the new first inequality gives $q>A'm\ge m$, while $K>0$ gives $q<1+(-e/f)m<m+1$. None of these intervals contains an integer. This verifies lattice-freeness throughout the changes below.

For $F'=rp+sq$ and $\Delta=as-br$, the identities are
\begin{align}
 (ad+bc)F'&=(rd+sc)F+\Delta G-\Delta c &&(\Delta>0),\\
 (af+be)F'&=(rf+se)F+(-\Delta)K-(-\Delta)f &&(\Delta<0).
 \label{eq:triangletransfers}
\end{align}
The corresponding values of $\theta$ in~\eqref{eq:preservingtransfer} are
\[
 \frac{\Delta c}{rd+sc+\Delta},\qquad
 \frac{(-\Delta)f}{rf+se-\Delta}.
\]
Every application below uses positive coefficients on $F,G$ or $F,K$ and satisfies $\theta<1$. If $a=kb+1$, Lemma~\ref{lem:residuepreservation} applies: its additional cuts use only $F,G$, so the negative coefficient in $K$ causes no change. For $k=1$, the replacement is $F=p+q$ with $e<0$, a case proved above without a lattice-freeness assumption. We use $k>1$ only after establishing the next lemma.

\begin{lemma}\label{lem:unittriangle}
Let $P$ be a rational polytope with all variables designated integer, let $h$ be rational affine, and let $p,q$ be affine with integer coefficients. Put
\[
 F=kp+q,\qquad G=c(1-p)+dq,\qquad K=-Ep+f(1-q),
\]
where $k\ge2$, $c,d>0$, and $0<E<f$ are integers with $\gcd(E,f)=1$. If $F,G,K$ are positive on $P\cap\{h>0\}$, there is a polynomial-time split proof ending exactly at $P\cap\{h\le0\}$, with polynomial total encoding in all supplied data.
\end{lemma}
\begin{proof}
If $P$ is empty or $\max_P h\le0$, no cuts are needed. Otherwise Lemma~\ref{lem:slackcertificate} supplies $F\ge\gamma_Fh$, $G\ge\gamma_Gh$, and $K\ge\gamma_Kh$ on $P$, with positive multipliers. Let $L$ be the total supplied input length. The coefficient reductions below end in one of three cases. We prove these first, using $U=H-t$ and the bounds $L_i\ge\alpha_i(h-t)$ introduced above.

\emph{The case $E=1$.} The identity
\[
 (fk-1)(1-q)=F+kK-1
\]
makes $1-q$ positive after increasing $t$ by $U/(k+1)$. After this increase, $h>t$ is impossible on $q\ge1$ and on either branch of the $p$ split when $q\le0$. Proposition~\ref{prop:crosssimulation} proves $h\le t$ exactly. Since $k\ge2$, the increase in $t$ is at most $U/3$. Thus $H-t$ retains at least $2/3$ of its initial value, and the fixed-factor repetitions prove the assertion.

\emph{The case $E=f-1$.} Let $\ell=\lceil\log_2(f+1)\rceil$ and $K_0=\max\{k,4\ell\}$. Replace $F$ by $K_0p+q$. If $D=K_0-k>0$, the identity in~\eqref{eq:triangletransfers} for $\Delta<0$ gives $\theta$ equal to
\[
 \frac{Df}{(K_0-1)f+1+D}.
\]
Thus $1-\theta\ge(k-1)/(K_0-1)\ge1/(4\ell)$. If $D=0$, no change is needed. After the cyclic change of coordinates above, the old $K$ has coefficients $(f,1)$; the function paired with it in the identity for $\Delta>0$ has coefficients $(1,K_0-1)$. While $f>K_0+1$, replace $f$ by $r=\lceil f/2\rceil$. Here $\theta$ equals
\[
 \frac{f-r}{r(K_0-1)+1+f-r}\le\frac1{K_0}.
\]
There are at most $\ell$ such changes, retaining a fraction at least $(1-1/K_0)^\ell\ge3/4$. Once $f\le K_0+1$, the change to $(1,1)$ has $\theta=(f-1)/(K_0+f-1)\le1/2$ and gives the case $F=p+q$ with $e<0$. In these coordinates, the third defining inequality has left-hand side $-d\widetilde p+(c+d)(1-\widetilde q)$, where $(\widetilde p,\widetilde q)=(1-p-q,p)$; hence $e=-d<0$ and $e+f=c>0$, as required. 

Combining these bounds, $H-t$ retains at least $3/(32\ell)$ of its initial value. The smaller fraction $1/(16\ell)$ therefore suffices for the predetermined repetitions. There are polynomially many repetitions, and the enlarged coefficient $K_0$ has polynomial encoding length.

\emph{The case $f=mE+1$, $m\ge2$.} Put $x=1-q$, $y=-p$, $a=mE+1$, $B_0=ac+Ed$, and $C_0=c+d$. Write
\[
 T=ax+Ey=K,\quad R=C_0-dx+cy=G,\quad U_0=mx+y,\quad D_0=x+y.
\]
Then
\begin{align*}
 B_0U_0&=(mc+d)T+R-C_0,\\
 B_0D_0&=C_0T+(a-E)R-(a-E)C_0,\\
 (a-E)U_0&=(m-1)T+D_0.
\end{align*}
Write $T\ge\gamma_Th$ and $R\ge\gamma_Rh$ on $P$, with $\gamma_T,\gamma_R>0$. The first identity gives $U_0\ge\alpha(h-V_0)$ for some $\alpha>0$, where
\[
 V_0=\frac{C_0}{(mc+d)\gamma_T+\gamma_R},\qquad \gamma_TV_0\le1.
\]
The replacement $U_0=-p+m(1-q)$ has $E=1$, so the first case gives an exact proof of this preliminary bound.

Write the second identity as $D_0\ge C_Dh-\beta$, with
\[
 \beta=\frac{(a-E)C_0}{B_0}\le m,\qquad
 C_D=\frac{C_0\gamma_T+(a-E)\gamma_R}{B_0},\qquad C_DV_0\le\beta.
\]
For completeness, the last inequality follows from $C_0\le(a-E)(mc+d)$; the first follows by expanding
\[
 mB_0-(a-E)(c+d)=c\bigl((m^2-m+1)E+m-1\bigr)+d(E-1)\ge0.
\]
At any predetermined bound $h\le V\le V_0$, the split on $x$ permits $D_0\ge mh/V-m$. For $h\le0$, this follows from $m/V\ge C_D$ and $m\ge\beta$. When $h>0$, $x\le0$ and $T>0$ give $D_0>0$, while $x\ge1$ and $R>0$ also give $D_0>0$. The proposed right-hand side is nonpositive because $h\le V$.

Apply Lemma~\ref{lem:supportnormalize} to ensure $T\ge h/V$ where $h>0$. The inequality also holds where $h\le0$: the original bound and $\gamma_TV_0\le1$ imply $T\ge\gamma_T h\ge h/V$ there. Thus the last identity gives globally
\[
 (a-E)U_0\ge(2m-1)h/V-m.
\]
At $h=4V/5$, the right-hand side is $(3m-4)/5>0$. The case $E=1$, already proved, therefore gives an exact split proof of $h\le4V/5$. Predetermined repetitions finish the proof. Retain all explicitly constructed additional inequalities; once an exact application is complete, later steps need only the bound it establishes. The exceptional value $m=1$ was covered by $E=f-1$, so no circular appeal is involved.

\emph{Reducing the general case.} For $E\ge2$, Lemma~\ref{lem:bezouthalving} gives $0\le E'\le E/2$, $1\le f'\le f/2$ and $\Delta=Ef'-E'f=\pm1$. Put $K'=-E'p+f'(1-q)$. Direct calculation gives
\begin{align*}
 (fk-E)K'&=F+(kf'-E')K-1 &&(\Delta=1),\\
 (fc+dE)K'&=G+(cf'+dE')K-(c+d) &&(\Delta=-1).
\end{align*}
The corresponding values of $\theta$ are $1/(1+kf'-E')\le1/(f'+1)$ and $(c+d)/(1+cf'+dE')\le1/E'$. A change with $\Delta=-1$ and $E'=1$ means $f=f'E+1$; stop before that change and use one of the preceding cases. A change with $\Delta=1$ and $E'=1$ has $\theta$ at most $1/2$ and reaches the first case. The possibility $E'=0$ occurs only when $E=1$, already treated. Equality $E'=f'$ can occur only at $(1,1)$ with $\Delta=-1$, where we have already stopped before the change.

All other changes have $E'\ge2$, $\theta\le1/E'\le1/2$, and preserve $0<E'<f'$. Since successive $E'$ halve, the sum of the reciprocal bounds is at most one. The product of the factors $1-\theta$ is at least $e^{-2}$; allowing the final change with $\Delta=1$ still leaves at least $e^{-2}/2>1/256$. Only after these reductions do we apply the proof of the resulting special case. The resulting upper bound is at most $(1-1/256)H$. Repeating from the original functions and finishing with the split on $Mh$ proves the assertion.

There are $O(L)$ coefficient reductions per repetition. To bound the encoding, put $H$, $U$, $U^{-1}$ and the three coefficients $\alpha_i$ over a common denominator, where $L_i\ge\alpha_i(h-t)$ and $t=H-U$. The identity~\eqref{eq:preservingtransfer} replaces a coefficient by $\sum_iw_i\alpha_i/A$. Increasing $t$ to $t+\theta U$ replaces $U$ by $(1-\theta)U$ and $U^{-1}$ by $U^{-1}/(1-\theta)$. The subsequent split cuts either retain a coefficient or replace it by the new $U^{-1}$. Thus, once these choices are made, each operation is a rational linear update of six quantities.

Suppose the initial common denominator and numerators have absolute value at most $2^{b_0}$, and each update coefficient has numerator and positive denominator bounded by $2^{b_1}$. After $T$ updates, the common denominator and numerators are bounded by
\[
 2^{\,b_0+T(37b_1+3)}.
\]
Indeed, clearing the at most $36$ denominators in one update increases the denominator bound by a factor of at most $2^{36b_1}$. Each new numerator is a sum of at most six terms, each enlarged by at most $2^{37b_1}$. Recovering $t=H-U$ requires at most one further bit. 

In our construction, $b_0$ is polynomial in $L$, $b_1=O(L)$ by the displayed identities, including the enlargement to $K_0$, and $T=O(L)$ per reduction sequence. Predetermined bounds are products of explicit rational factors, so their bit lengths also grow additively with the number of repetitions. The number of nested applications is bounded by a constant, and all additional inequalities have the explicit encoding bounds just obtained. The total running time and encoding are therefore polynomial.
\end{proof}

We can now finish the general triangle. The argument uses short continued-fraction reductions, because losing a fixed fraction at each of $O(L)$ changes would leave an exponentially small improvement.

\begin{proof}[Proof of Theorem~\ref{thm:planarsimulation}]
Let $L$ be the total input length. The enlargement and elementary cases above leave only~\eqref{eq:signedtriangle}. Choose cyclic coordinates so that $A=a/b$ is the largest of the three original ratios; in particular $C\le A$. If $b=1$, Lemma~\ref{lem:unittriangle} applies. Otherwise fix a regular continued-fraction expansion of the original $A$ and retain its initial segments, including those ending in one.

Every noninteger convergent $r/s$ of $A>1$ satisfies $r/s\ge2-1/A$. If $\lfloor A\rfloor\ge2$ this is immediate. Otherwise write $A=[1;a_1,\ldots]$. Every noninteger convergent is at least $1+1/(a_1+1)\ge2-1/A$. Since $C\le A$, it follows that $r/s\ge2-1/C$. For $\Delta<0$, set $D=-\Delta$. The value of $\theta$ from~\eqref{eq:triangletransfers} is
\[
 \frac{D}{r-s+s/C+D/f}\le\frac Ds.
\]
The identity coefficients are positive because $rf+se=f(r-s+s/C)>0$. For $\Delta>0$, we have $\theta\le\Delta/s$ since $d>0$.

At a current convergent $a/b$, use its first convergent with $s_i^2\ge2b$, followed by its predecessor in the fixed expansion, as in Proposition~\ref{prop:crosssimulation}. The first change has $\theta$ at most $1/2$ by the convergent error bound. A noninteger predecessor has determinant of magnitude one and denominator at least two, so its value of $\theta$ is also at most $1/2$. Its denominator is below $\sqrt{2b}$.

For a predecessor with denominator one and positive determinant, the current pair is $(kb+1,b)$. Stop before the change and apply Lemma~\ref{lem:residuepreservation}, using Lemma~\ref{lem:unittriangle} for the replacement $(k,1)$, or the $e<0$ case $F=p+q$ if $k=1$, since $G$ and $K$ are unchanged. A negative determinant gives $(kb-1,b)$, with $k\ge2$, and $\theta$ equal to
\[
 \frac1{k-1+1/C+1/f}.
\]
This is at most $1/2$ if $k\ge3$. If $k=2$, $2/1$ occurs as a convergent of the original expansion, so $A<3$ and $C<3$; $\theta$ is then less than $3/4$. Make this change and apply Lemma~\ref{lem:unittriangle}. 

For $2\le b\le7$, use successive predecessors in that expansion as in Proposition~\ref{prop:crosssimulation}. Noninteger predecessors remain convergents of the original expansion, so at most five changes with $\theta$ at most $1/2$ reach a denominator-one predecessor. The endpoint rules above apply, with $\theta<3/4$ for the final change when $\Delta<0$.

There are $N=O(\log L)$ changes, each retaining at least a quarter of $H-t$. Thus the final value of $t$ satisfies $H-t\ge4^{-N}H$. The exact final proof gives $h\le t$, reducing the initial upper bound by at least the fraction $4^{-N}$. This fraction depends only on the original coefficients and is inverse polynomial in $L$. Use it for predetermined repetitions on the original polytope with the corresponding upper bounds. 

It remains to bound the encoding. All changed coefficients have polynomial encoding; the only increase in the preceding lemmas was the explicitly bounded enlargement to $K_0$. The explicit encoding bound in Lemma~\ref{lem:unittriangle} applies with $T=O(\log L)$ and $b_0,b_1$ polynomial in the original input length. The number of nested applications is bounded by a constant, and no expanded system of cut inequalities is used as input to a later repetition. Thus construction time and total encoding are polynomial in all supplied data, even when the original dimension varies.

Every cut preserves the original set $P\cap\{h\le0\}$, and the final split on $Mh$ imposes $h\le0$. The final relaxation is exactly this set. Empty polytopes and nonpositive maxima need no cuts; a positive constant $h$ is handled by the contraction itself. Constant affine split functions with integer coefficients have one branch equal to the whole relaxation and contribute only implied inequalities.
\end{proof}

\subsection{Two-row systems}

We first apply the theorem to a system defined by a lattice-free triangle. Let $P$ be a rational polytope with all variables designated integer. Let $T=\conv(w_1,w_2,w_3)$ be any rational lattice-free triangle and $f\in\operatorname{int}T$ rational. Let $u$ be an affine map into $\R^2$ with integer coefficients and $s_1,s_2,s_3$ rational affine on $P$, with $s_i\ge0$ and $u=f+\sum_i(w_i-f)s_i$ on $P$. With $h=1-\sum_i s_i$,
\[
 u=hf+\sum_i s_iw_i.
\]
When $h>0$, this convex combination has positive weight on $f$, so $u\in\operatorname{int}T$. Theorem~\ref{thm:planarsimulation} therefore derives $\sum_i s_i\ge1$ from these supplied affine functions. The affine functions $s_i$ need not take integer values at integer points.

We next enforce membership in the mixed-integer hull of a supplied system $x=f+Rs$, $s\ge0$, with $x\in\Z^2$. Polynomial facet enumeration bounds the number of applications of Theorem~\ref{thm:planarsimulation}.

\begin{theorem}\label{thm:cornersimulation}
Let $P\subseteq\R^n$ be a rational polytope with all variables designated integer. Suppose the input supplies an affine map $u:\R^n\to\R^2$ with integer coefficients, a rational affine map $s:\R^n\to\R^k$, and rational data $f\in\Q^2$, $R\in\Q^{2\times k}$ such that
\[
s(z)\ge0,\qquad u(z)=f+Rs(z)\quad(z\in P).
\]
One can construct a split proof in the original $z$ variables whose final relaxation satisfies
\[
s(z)\in\conv\{s\ge0:f+Rs\in\Z^2\}.
\]
The construction time, number of cuts, and total binary encoding length are polynomial in the encoding of $P,u,s,f,R$. Neither $n$ nor $k$ is fixed.
\end{theorem}

\begin{proof}[Proof of Theorem~\ref{thm:cornersimulation}]
Append the four columns $\pm e_1,\pm e_2$ to $R$, and append four identically zero functions to $s$. The nonnegative combinations of the enlarged set of columns fill $\R^2$. The face of the enlarged mixed-integer hull where the four new coordinates vanish is exactly the original hull: nonnegative coordinates cannot cancel in a convex combination. If $f$ is integral, the hull is the nonnegative orthant. Indeed, it contains the origin and arbitrarily distant feasible points on each rational coordinate ray, and their convex combinations fill the orthant. Omit zero columns during facet enumeration and retain the nonnegativity constraints on their coordinates.

Otherwise enumerate the polynomially many facet inequalities, written as $\alpha^Ts\ge1$, using~\cite[Theorem~1.5]{BHK}. Every coefficient $\alpha_j$ is nonnegative, since the positive coordinate directions are recession directions. Write $r_j$ for the columns of the enlarged matrix $R$. The images $f+Rs$ of points with $s\ge0$ and $\alpha^Ts<1$ lie in the interior of the full-dimensional rational lattice-free polyhedron
\[
 L_\alpha=\conv\{f+r_j/\alpha_j:\alpha_j>0\}
            +\cone\{r_j:\alpha_j=0\}.
\]
We verify the interior containment and lattice-freeness using the identity
\[
 L_\alpha=\{f+Rt:t\ge0,\ \alpha^Tt\le1\}.
\]
Since the nonnegative combinations of the columns fill $\R^2$, there is a vector $\lambda>0$ with $R\lambda=0$. This places $f$ in the convex hull plus cone defining $L_\alpha$; scaling and taking convex combinations then gives the identity. The same spanning property puts $f$ in the interior. If $t\ge0$ and $\alpha^Tt<1$, add a small positive multiple of $\lambda$ to $t$. Its image is unchanged, every coordinate becomes positive, and $\alpha^Tt$ remains below one. Small perturbations of two coordinates whose columns are independent show that the image is interior. Conversely, take an interior point other than $f$, extend the segment from $f$ slightly beyond that point within $L_\alpha$, and scale the resulting representation toward $f$. This gives a representation with $\alpha^Tt<1$; $f$ itself has the representation $t=0$. An integer point in the interior would therefore contradict validity of $\alpha^Ts\ge1$.

Apply Theorem~\ref{thm:planarsimulation} with $h=1-\alpha^Ts(z)$ to each facet, always starting from the original $P$. Concatenate the proofs as described at the start of this section. Polynomial facet enumeration bounds the total supplied data, and the polynomial bounds of Theorem~\ref{thm:planarsimulation} give the asserted bounds.
\end{proof}

After $r$ cuts of a split proof, its relaxation contains $\SC^r(P)$, by induction and monotonicity. Thus polynomially many full split-closure rounds also enforce membership in the supplied mixed-integer hull; this need not be the integer hull of $P$.

\begin{corollary}\label{cor:cornerproofsimulation}
For bounded pure-integer input, an explicitly encoded proof using split cuts, inequalities valid for supplied two-row mixed-integer hulls, cross cuts, and the lattice-free cuts of Theorem~\ref{thm:planarsimulation} can be replaced by a split proof with total encoding length and construction time polynomial in the encoding of the supplied proof. Each supplied step includes its defining affine maps and two-row system, split functions, or lattice-free polyhedron, as applicable.
\end{corollary}

\begin{proof}
Construct the corresponding split proof on the polytope supplied at each step of the original proof. The final set satisfies the inequality required by that step. Induction therefore gives containment in the polytope for the next step. The starting polytopes have polynomial total encoding in the supplied proof encoding, even if each description repeats all earlier inequalities. Concatenating their split proofs preserves the required containments and gives the stated bounds.
\end{proof}

Related simulation results arise in proof complexity. Here Cutting Planes denotes the proof system based on linear combinations and division, rather than arbitrary split cuts. Beame et al.~\cite{BFIKPPR18} introduced Stabbing Planes, a branching proof system that polynomially simulates Cutting Planes. Basu et al.~\cite{BCDJComplexityII} extended this simulation to split-cut proofs. Dadush and Tiwari~\cite[Theorem~4]{DT20} convert enumerative branching refutations into Chv\'atal--Gomory proofs. Fleming et al.~\cite{FGIPRTW21} give a quasipolynomial simulation by Cutting Planes of Stabbing Planes refutations of CNF formulas with quasipolynomially bounded coefficient magnitudes. Corollary~\ref{cor:cornerproofsimulation} instead replaces the specified cut operations by split cuts, with polynomial construction time and total encoding. The construction in Theorem~\ref{thm:planarsimulation} retains every point satisfying the target inequality, including fractional points. The corollary does not assert a polynomial simulation of arbitrary Stabbing Planes proofs.

\paragraph{Imposing two integrality constraints.}

In fixed dimension the planar lattice-free set need not be supplied. The projection argument of Dash, Dey, and G\"unl\"uk~\cite[Lemma~2.1]{DDG11} constructs one from an inequality valid whenever two specified functions are integer. The preservation conclusion of Theorem~\ref{thm:planarsimulation} then gives the following stronger form of simulation.

\begin{corollary}\label{cor:twoformsimulation}
Fix $d$. Let $P\subseteq\R^d$ be a rational polytope with all coordinates designated integer, and let $u=(p,q)$ consist of two supplied affine functions with integer coefficients, allowing dependent linear parts. If a supplied rational affine inequality $h\le0$ holds at every point of $P$ with $u(z)\in\Z^2$, one can construct a split proof ending exactly at $P\cap\{h\le0\}$. One can also construct a split proof ending exactly at
\[
 U_u(P)=\conv\{z\in P:u(z)\in\Z^2\}.
\]
Both constructions take polynomial time and have polynomial total encoding in the supplied data.
\end{corollary}
\begin{proof}
Let $S=u(P\cap\{h>0\})$. If this image is empty, no cuts are needed. Otherwise fixed-dimensional elimination, with strict inequalities retained, gives a polynomial description
\[
 S=\{u:A_1u\le b_1,\ A_2u<b_2\}.
\]
Multiply each inequality by a positive common denominator to make its coefficients and right-hand side integral. Replacing the nonstrict inequalities by $A_1u<b_1+\mathbf1$ leaves exactly the same integer solutions, while retaining all strict inequalities. The resulting set $O$ contains $S$, is nonempty and open, and has no integer point. Its closure is a full-dimensional rational lattice-free polyhedron with interior $O$. This also handles a lower-dimensional image $S$. Theorem~\ref{thm:planarsimulation} gives the first assertion and its exact preservation property.

To construct $U_u(P)$, let $r$ be the rank of the linear part of $u$. For a fixed $v\in u(P)$, consider the polytope $P\cap\{z:u(z)=v\}$. Each vertex is the unique solution of $r$ independent equations from $u(z)=v$ and $d-r$ tight inequalities defining $P$, taken as equations, whose combined coefficient matrix is nonsingular. Enumerate these choices. Each gives an affine map $z=T v+t$ on the affine image of $u$, with a bounded rational domain in the two output coordinates $v$. Compute the integer hull of this domain and apply the affine map $v\mapsto Tv+t$ to its vertices. If the linear parts of the two functions are dependent, retain the equations defining the affine image of $u$ and require both coordinates of $v$ to be integer. Each nonempty polytope $P\cap\{z:u(z)=v\}$ is the convex hull of its vertices. Taking these sets over $v\in\Z^2$ therefore shows that the convex hull of all lifted points is exactly $U_u(P)$. There are polynomially many systems for fixed $d$, and fixed-dimensional integer-hull algorithms~\cite{CHKM,H} give polynomial encoding and running time.

Construct a defining inequality description of this hull, including both signs of affine-hull equations. Prove each defining inequality using the first assertion, always starting from the original $P$, and concatenate the proofs by monotonicity. Each proof retains the desired hull, and together they give exactly that hull. If the hull is empty use the constant inequality $1\le0$. If both linear parts vanish, the supplied functions are integer constants and the hull is $P$.
\end{proof}

The same argument works for one supplied function, by adjoining an integer constant. More generally, if a rational polytope $G\subseteq P$ is to be retained, every defining inequality of $\conv(G\cup U_u(P))$ meets the first assertion and preserves $G$.

\section{Split proofs in fixed dimension}\label{sec:fixedproof}

We derive each defining inequality $h\le0$ of $P_I$ by induction on dimension. Fixing an integer coordinate reduces the dimension, but convex combinations of points on other sections can restore points removed from the chosen section. The preservation property from Section~\ref{sec:cornersimulation} lets us lower a bound on the chosen section while retaining these contributions.

Unboundedness introduces one additional consideration. A recession direction along which $h$ is constant permits a projection to lower dimension. Once no such direction remains, boundedness of $h$ above implies that $P\cap\{h\ge\alpha\}$ is bounded for each $\alpha$. The argument can then proceed on this polytope, provided its cuts extend to $P$ without removing points below the chosen level.

\subsection{Preservation under split cuts}

\begin{proposition}\label{prop:prescribedpreservation}
Fix $d$. Let $P\subseteq\R^d$ be a rational polyhedron and let $h$ be an affine function with integer coefficients. Suppose $h$ is bounded above on $P$ and $h\le0$ on $P\cap\Z^d$. One can construct a split proof ending exactly at $P\cap\{h\le0\}$, with every cut retaining this set. Its running time and total encoding length are polynomial in $(P,h)$.
\end{proposition}

Here an integer-free set is one containing no integer point, and a refutation is a split proof whose final relaxation is empty. Assuming the proposition in dimension $d-1$, we first reduce an integer-free polytope in dimension $d$ to a specified face. This also gives a refutation of any integer-free polytope: reduce it to a section with one integer coordinate fixed, and refute that section in dimension $d-1$. We then prove the proposition in dimension $d$. Empty inputs and dimension zero are immediate.

We use fixed-dimensional polyhedral and integer-hull algorithms~\cite{L83,H,CHKM}. Flatness gives a primitive integer direction of width at most a fixed integer $F_d\ge2$ for every nonempty integer-free polytope in $\R^d$~\cite{BLPS99}. Such a direction has polynomial encoding and can be found in polynomial time: enumerate the vertices $v_j$ and seek a nonzero integer $a$ satisfying $a^T(v_j-v_k)\le F_d$ for all pairs. To exclude zero, consider the $2d$ alternatives $a_i\ge1$ and $a_i\le-1$. Fixed-dimensional integer feasibility applies, and dividing by the gcd makes $a$ primitive without increasing width.

Extend the notation of Corollary~\ref{cor:twoformsimulation} to a tuple $u$ of one or two affine functions with integer coefficients. For rational polytopes $G\subseteq D$, write
\begin{equation}\label{eq:retainedhull}
 U_u(D)=\conv\{x\in D:u(x)\text{ is integral}\},
 \qquad \mathcal T_u^G(D)=\conv\bigl(G\cup U_u(D)\bigr).
\end{equation}
The corollary constructs this last polytope and a split proof ending exactly there in polynomial time, retaining $G$ at every cut. It also lets us impose any supplied polytope containing $\mathcal T_u^G(D)$ and contained in $D$: prove its defining inequalities from $D$ and concatenate the proofs. For one function, adjoin an integer constant as the second function.

We will also use a simple consequence of split validity. If a cut with split function $u$ changes $D$ to $D'$ and retains $G$, then
\begin{equation}\label{eq:retainedoperation}
 \mathcal T_u^G(D)\subseteq D'\subseteq D.
\end{equation}
Indeed, every point of $D$ with integer $u$ lies in one of the split branches. As in Section~\ref{sec:cornersimulation}, we may construct cuts on an auxiliary polytope containing the current relaxation. We verify this containment at each step and construct the auxiliary polytopes from the original data and bounds of polynomial encoding.

\subsection{Integer-free polytopes}

To reduce an integer-free polytope to a face, we bound the function defining that face on parallel sections. We choose the bounds to form a concave piecewise-linear function. At each section, the segment joining the neighboring bounds lies strictly below the current bound. Its value there bounds the contribution from all other sections. This leaves room for a decrease that preserves those contributions. A complete pass through the sections decreases the upper bound by a fixed factor.

\begin{lemma}\label{lem:bottomface}
Let $C\subseteq\R^d$ be a rational polytope containing no integer point, and let $t$ be an affine function with integer coefficients. Suppose $C\subseteq\{t\ge0\}$ and $L=C\cap\{t=0\}\ne\varnothing$. Assuming Proposition~\ref{prop:prescribedpreservation} in dimension $d-1$, one can construct a polynomial-time split proof ending exactly at $L$, with polynomial total encoding in $(C,t)$ and with every cut retaining $L$.
\end{lemma}
\begin{proof}
By flatness, choose a primitive integer vector $v$ so that $s(x)=v^Tx$ has range $[l,r]$ of length at most $F_d$ on $C$. Let
\[
 \begin{gathered}
 m=\lfloor l\rfloor-1,\quad n=\lceil r\rceil+1,\quad c=(m+n)/2,\\
 a=\frac1{100(F_d+4)^2},\qquad
 \phi_j=2-a(j-c)^2\quad(m\le j\le n).
 \end{gathered}
\]
Interpolating linearly between integer arguments gives a concave function $\phi$ that is greater than one on $[m,n]$ and at most two. At every interior integer $j$,
\[
 2\phi_j-\phi_{j-1}-\phi_{j+1}=2a.
\]
All integer values of $s$ attained in $C$ lie at interior arguments, and their number depends only on $d$.

Let $H=\max_Ct$. If $H=0$, the assertion already holds. Otherwise $t\le H\phi(s)$, and we will lower this bound successively to $T_j\phi(s)$, where $T_j=H\rho^j$ and $\rho=1-a/4$. At scale $T$, apply~\eqref{eq:retainedhull} to $C\cap\{t\le T\phi(s)\}$ with $u=s$ and $G=L$. During the pass through the sections, write $\psi$ for the current concave piecewise-linear upper bound, initially $\psi=T\phi$, and put $\psi_i=\psi(i)$. At the start of each section step we use the auxiliary polytope
\begin{equation}\label{eq:profilereference}
 R=\mathcal T_s^L\bigl(C\cap\{t\le\psi(s)\}\bigr).
\end{equation}
Thus $R$ is the convex hull of $L$ and its sections with integer $s$.

Process each integer level in $[l,r]$ once. At a selected level $i$, put
\[
 E_i=\conv\bigl(L\cup\{x\in R:s(x)\in\Z,\ s(x)\ne i\}\bigr),
 \quad B_i=R\cap\{s=i\},\quad G_i=E_i\cap\{s=i\}.
\]
The set $G_i$ consists of the points that the other sections and $L$ contribute at $s=i$. Any proof on $B_i$ must retain these points if it is to extend to the convex hull with $E_i$; see Figure~\ref{fig:layer-interference}.

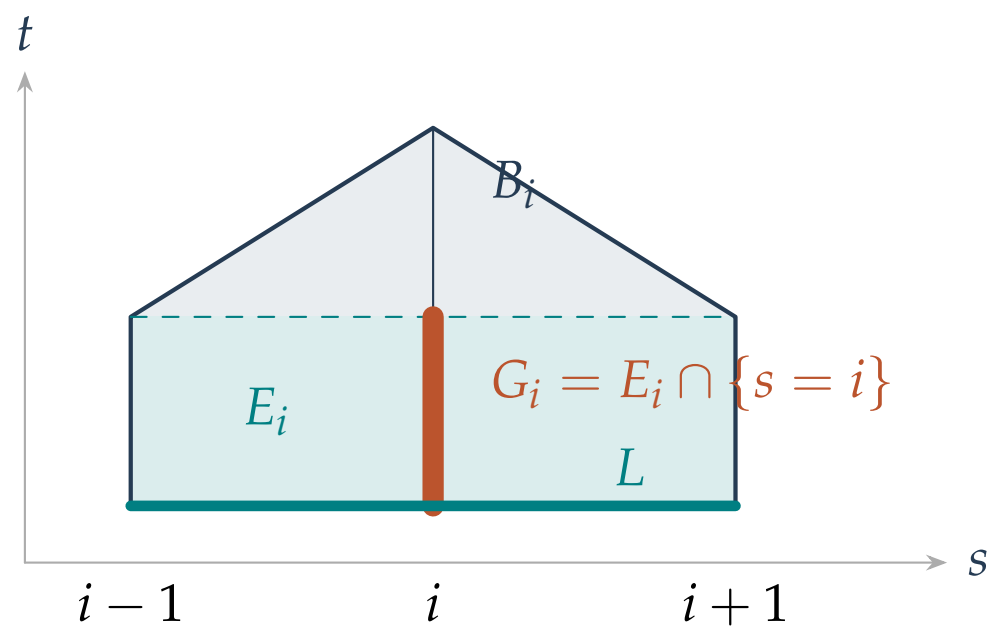
\begin{figure}[htbp]
\centering
\begin{adjustbox}{max width=\linewidth}
\begin{tikzpicture}[x=2cm,y=1.25cm,font=\small]
\fill[lightgray] (-1,0)--(-1,1)--(0,2)--(1,1)--(1,0)--cycle;
\fill[cutteal!14] (-1,0) rectangle (1,1);
\draw[boundary] (-1,0)--(-1,1)--(0,2)--(1,1)--(1,0)--cycle;
\draw[cutteal,dashed] (-1,1)--(1,1);
\draw[inkblue,thin] (0,0)--(0,2);
\draw[warmrust,line width=4pt] (0,0)--(0,1);
\draw[cutteal,line width=2pt] (-1,0)--(1,0);
\draw[axis] (-1.35,-.3)--(1.7,-.3) node[right] {$s$};
\draw[axis] (-1.35,-.3)--(-1.35,2.3) node[above] {$t$};
\node[below] at (-1,-.3) {$i-1$};
\node[below] at (0,-.3) {$i$};
\node[below] at (1,-.3) {$i+1$};
\node[text=inkblue,anchor=west] at (.12,1.7) {$B_i$};
\node[text=warmrust,anchor=west] at (.12,.65) {$G_i=E_i\cap\{s=i\}$};
\node[text=cutteal] at (-.55,.5) {$E_i$};
\node[text=cutteal,anchor=south] at (.65,0) {$L$};
\end{tikzpicture}
\end{adjustbox}
\caption{Convex combinations of points on other sections contribute $G_i$ at $s=i$. Lowering the bound on $B_i$ must retain this set. The picture suppresses the remaining coordinates; an integer value of $s$ does not make a point integral.}
\label{fig:layer-interference}
\end{figure}

Replace the graph of $\psi$ between $i-1$ and $i+1$ by the segment joining its endpoints. The resulting function is concave and still bounds $t$ on all sections contributing to $E_i$. It also bounds $t$ on $L$, since it is positive and $t=0$ there. At an unprocessed level $i$, its value is at most $T(\phi_i-a)$: the initial value is exactly this, and previous changes only lowered neighboring values. Consequently
\[
 t\le g_i:=T(\phi_i-a)\qquad\text{on }G_i.
\]
Since $s$ has a primitive coefficient vector, we may use unimodular coordinates $(s,y)$, with $y\in\R^{d-1}$. On $s=i$, the function $t$ has integer coefficients in $y$, and $B_i$ contains no integer point. Thus Proposition~\ref{prop:prescribedpreservation} in dimension $d-1$ applies to a positive integer multiple of $t-g_i$ with integer coefficients. The resulting split proof ends exactly at $B_i\cap\{t\le g_i\}$ and retains $G_i$. If $B_i$ is empty, no section proof is needed.

Retaining $G_i$ allows us to extend this proof to the whole polytope. Write $Q$ for the actual relaxation and $D$ for the current relaxation on the section. At the end of each section step we maintain
\[
 G_i\subseteq D\subseteq B_i,\qquad
 Q\subseteq R_i(D):=\conv(E_i\cup D).
\]
Initially $D=B_i$ and $R_i(D)=R$. Its section at $s=i$ is exactly $D$, since $E_i\cap\{s=i\}=G_i\subseteq D$. Every point of $R_i(D)$ with integer $s\ne i$ belongs to $E_i$, because $R_i(D)\subseteq R$. If $u$ is the affine split function used for the next cut on $D$, extend it from $y$ to $(s,y)$ independently of $s$. It has integer coefficients in the original coordinates, and
\begin{equation}\label{eq:sectionlift}
 \mathcal T_{(s,u)}^L\bigl(R_i(D)\bigr)
 \subseteq\conv\bigl(E_i\cup\mathcal T_u^{G_i}(D)\bigr).
\end{equation}
Indeed, a point with integer $s,u$ belongs to $U_u(D)$ when $s=i$ and to $E_i$ otherwise; also $L\subseteq E_i$. If the next section relaxation is $D'$, inclusion~\eqref{eq:retainedoperation} shows that the right-hand side of~\eqref{eq:sectionlift} is contained in $R_i(D')$. Corollary~\ref{cor:twoformsimulation} therefore imposes $R_i(D')$ on $R_i(D)$ by split cuts retaining $L$. The same cuts are valid on $Q$ and establish the containment for the next step.

At the end of the section proof, we can decrease $\psi_i$ by $Ta/2$, half the decrease just established. Both $D$ and $E_i$ satisfy this weaker bound. Concavity is preserved: $2\psi_i-\psi_{i-1}-\psi_{i+1}$ decreases by $Ta$ from a value at least $2Ta$, while the corresponding quantities at neighboring indices increase. All values remain positive. Reconstruct~\eqref{eq:profilereference} using the updated $\psi$. It contains the preceding $R_i(D)$, since that polytope is generated by points in $L$ and on sections with integer $s$, all satisfying the new bound. It therefore contains $Q$ as well.

After all levels have been processed, their values are at most
\[
 T(\phi_i-a/2)\le\rho T\phi_i,
\]
where we used $\phi_i\le2$. The final auxiliary polytope is generated by these sections and $L$. Concavity of $\phi$ and $t=0$ on $L$ give $t\le\rho T\phi(s)$ on their convex hull, justifying the next scale. After $O_d(1+\log(1+H))$ repetitions, $t<1$. The operation $\mathcal T_t^L$ then leaves exactly $L$, since zero is the only available integer value of $t$. This step requires integer coefficients of $t$, but does not require its coefficient vector to be primitive.

For the encoding bound, every value of $\psi$ at scale $j$ has the form
\[
 H\rho^j(\phi_i-\varepsilon_i a/2),\qquad \varepsilon_i\in\{0,1\}.
\]
If $b$ is the encoding length of $H$, these values have encoding length $O_d(b+j+1)$. The indices $m,n$ have polynomial encoding, and $n-m$ depends only on $d$. Each polytope in~\eqref{eq:profilereference} therefore uses the original $C$ and a constant number of inequalities of polynomial encoding; its sections and $E_i$ have polynomial descriptions as well. Induction bounds the entire proof on each section. To extend it, we use its supplied intermediate relaxations $D$ and the fixed $E_i$, so the coefficients of earlier lifted cuts do not enter the next construction. Polynomially many scales, with a dimension-dependent number of sections at each scale, give the asserted time and encoding bounds.
\end{proof}

\subsection{Preservation of valid inequalities}

A proof on the part above a given level can be extended while retaining every point below it, provided it retains the section at that level. This section contains the crossing point of every segment joining the two parts.

\begin{lemma}\label{lem:extendpreservation}
Fix $d$. Let $P\subseteq\R^d$ be a rational polyhedron, let $h$ be rational affine, and let $\alpha\in\Q$. Put
\[
 A=P\cap\{h\le\alpha\},\qquad C=P\cap\{h\ge\alpha\},
 \qquad F=P\cap\{h=\alpha\}.
\]
Suppose $C$ is bounded and a supplied split proof on $C$ ends exactly at $F$, retaining $F$ at every cut. One can extend it to a split proof on $P$ ending exactly at $A$, retaining every point of $A$. The construction time and total encoding are polynomial in $(P,h,\alpha)$ and the supplied proof.
\end{lemma}
\begin{proof}
If $C$ is empty, no cuts are needed. If $A$ is empty, $P=C$ and the supplied proof suffices. Otherwise let $D_0=C,\ldots,D_N=F$ be its relaxations and put $R_j=A\cup D_j$. This union is convex: a segment from $A$ to $D_j$ crosses $h=\alpha$ in $F\subseteq D_j$, so the part above that level lies in $D_j$ and the part below lies in $A$. It is also closed, as a finite union of closed sets. Thus $R_j=\conv(A\cup D_j)$ is a rational polyhedron, with $R_0=P$ and $R_N=A$.

Suppose the split producing $D_{j+1}$ has branches $H_j^-$ and $H_j^+$. For either branch,
\[
 R_j\cap H_j^\pm
 =(A\cap H_j^\pm)\cup(D_j\cap H_j^\pm)
 \subseteq A\cup D_{j+1}=R_{j+1}.
\]
Every defining inequality of $R_{j+1}$ is therefore a split cut for $R_j$, justified by this same split. Impose these inequalities one at a time, including both signs of any equations. They remain valid after earlier cuts, retain $A$, and leave exactly $R_{j+1}$.

Construct each $R_j$ from the original $A$ and the supplied $D_j$. In fixed dimension, we can represent these polyhedra by finitely many points and recession directions, and convert their closed convex hull to inequalities in polynomial time. Lineality is represented by a rational basis. The closed convex hull introduces no extra points, since $R_j$ is closed. Each description has polynomial encoding in the supplied data, so the same holds for the entire extended proof.
\end{proof}

\begin{proof}[Proof of Proposition~\ref{prop:prescribedpreservation}]
Assume the proposition in dimension $d-1$. We first obtain a refutation of any integer-free rational polytope $B\subseteq\R^d$. If $B$ is nonempty, set $k=\lceil\min_Bx_1\rceil$. The branch $x_1\le k-1$ is empty, so the split at $k-1,k$ permits $x_1\ge k$. If $B'=B\cap\{x_1\ge k\}$ is nonempty, its face $F=B'\cap\{x_1=k\}$ is nonempty. Lemma~\ref{lem:bottomface}, applied with $t=x_1-k$, reduces $B'$ exactly to $F$. On this face, fixing the integer coordinate $x_1=k$ leaves an integer-free polytope in $\R^{d-1}$. The inductive hypothesis with constant target $1\le0$ refutes it. Its cuts and splits extend directly to the original coordinates, since the entire current relaxation lies on $F$. This proves the required refutation with polynomial time and total encoding.

Now consider $(P,h)$ in the proposition. Empty $P$ and $h\le0$ on $P$ require no cuts. Otherwise write $h(x)=a^Tx+b$, and let $\operatorname{rec}(P)$ denote the recession cone of $P$. Since $h$ is bounded above,
\[
 a^Tr\le0\qquad(r\in\operatorname{rec}(P)).
\]
Suppose first that $\operatorname{rec}(P)$ contains a nonzero $r$ with $a^Tr=0$. Choose such an $r$ primitive and integral, and let $U$ be a unimodular matrix whose last column is $r$. In the coordinates $x=U(y,t)$, the function $h$ is independent of $t$, so we write it as $\bar h(y)$. Let $Q$ be the projection of $U^{-1}P$ onto $y\in\R^{d-1}$.

Every integer $y\in Q$ has an integer lift to $P$. Indeed, choose a feasible $t_0$. Since $r\in\operatorname{rec}(P)$, every $t\ge t_0$ is feasible with this $y$, so we may take $t$ integer. Unimodularity makes $U(y,t)$ integer, and hence $\bar h(y)\le0$. Also $\bar h$ has the same range on $Q$ as $h$ on $P$, so it is bounded above. Apply induction to $(Q,\bar h)$ and substitute $y$ as an integer linear function of $x$ in every cut and split. The split normals remain nonzero integer vectors. Projection commutes with imposing inequalities depending only on $y$, so the resulting proof is valid on $P$ and ends exactly at $P\cap\{h\le0\}$.

We may therefore suppose that $a^Tr<0$ for every nonzero $r\in\operatorname{rec}(P)$. If $P\cap\{h\le0\}$ is empty, $P$ contains no integer point and must be bounded: otherwise, moving along a nonzero recession direction would eventually give $h\le0$. Use the bounded refutation proved above.

It remains to treat a nonempty target. Put $H=\max_Ph>0$. If $H<1$, the split $h\le0$ or $h\ge1$ proves the assertion. Otherwise put
\[
 C=P\cap\{h\ge\tfrac12\},\qquad F=P\cap\{h=\tfrac12\}.
\]
The recession cone of $C$ is $\{0\}$, so $C$ is bounded. It contains no integer point, and $F$ is nonempty because $P$ has points with $h\le0$ and with $h\ge1$. Since $t=2h-1$ has integer coefficients, Lemma~\ref{lem:bottomface} gives a split proof on $C$ ending exactly at $F$. Extending this proof by Lemma~\ref{lem:extendpreservation} leaves exactly $P\cap\{h\le\tfrac12\}$. The final split $h\le0$ or $h\ge1$ gives exactly $P\cap\{h\le0\}$, retaining every point of this set throughout.

All data used here have polynomial encoding. To find a nonzero vector in $\operatorname{rec}(P)\cap\ker a^T$, use at most $2d$ rational LP feasibility tests, imposing $r_i\ge1$ or $r_i\le-1$. Clear denominators and divide by the gcd to obtain a primitive integer vector. Unimodular completion and projection take polynomial time in fixed dimension. Lemma~\ref{lem:bottomface} uses polynomially many lower-dimensional calls on data of polynomial encoding in the input of this induction step; the bounded refutation makes one further such call. Lemma~\ref{lem:extendpreservation} is polynomial in its supplied proof. Since recursive calls always decrease dimension, these bounds compose only a fixed number of times. This completes the induction.
\end{proof}

\subsection{The integer hull}

\begin{proof}[Proof of Theorem~\ref{thm:fixedproof}]
If $P_I$ is empty, Proposition~\ref{prop:prescribedpreservation} with $h=1$ gives a split refutation. Otherwise compute a polynomial-size inequality description of $P_I$ using fixed-dimensional integer-hull algorithms~\cite{H,CHKM}, including both signs of any affine-hull equations. These algorithms apply to unbounded rational polyhedra as well. One may also reduce the computation to a bounded polytope: in polynomial time, find a rational polytope $B\subseteq P$ such that
\[
 P=B+\operatorname{rec}(P),\qquad
 P_I=B_I+\operatorname{rec}(P);
\]
see~\cite[Lemma~2]{HOW15}. Compute $B_I$ and add the recession cone. In particular, $\operatorname{rec}(P_I)=\operatorname{rec}(P)$.

After clearing denominators, write each defining inequality of $P_I$ as $h\le0$ with $h$ affine with integer coefficients. The recession-cone identity implies that $h$ is bounded above on $P$, so Proposition~\ref{prop:prescribedpreservation} gives an exact split proof of it from the original $P$. Concatenate these proofs, keeping all earlier cuts. By monotonicity, each cut remains valid on the smaller current relaxation. Every cut retains $P_I$, and the final relaxation satisfies its entire description, so it equals $P_I$.

There are polynomially many defining inequalities of polynomial encoding. The construction time and total proof encoding are therefore polynomial in the original input length. The split-rank bound follows from the monotonicity argument in Section~\ref{sec:introproofs}.
\end{proof}

\section{Split separation and optimization}\label{sec:algorithms}\label{sec:sequential}

Finding a violated split cut does not by itself give a sequence reaching an optimal integer solution. For separation, the cut formula leads to bounded integer quadratic minimization in fixed total dimension. For an objective involving only the integer variables, we work on the current optimal face, so each cut removes an optimal LP vertex. This requires only fixed integer dimension, but the number of objective levels encountered introduces a dependence on the integer objective gap.

\subsection{Split separation}\label{sec:separation}

Section~\ref{sec:normals} bounds the coordinates of the split normals needed to define the closure. The bounds have polynomial encoding length, but the box they define may contain exponentially many integer points. We therefore need to search this box without enumerating its integer points.

By~\eqref{eq:clearedsplit}, a split inequality is violated exactly when a quadratic polynomial in the split normal and its integer right-hand side is negative. On each cone where the signs of the scalar products $\pi^Tr_j$ are fixed, bounded integer quadratic minimization therefore finds a violated cut or certifies that none exists. The quadratic problems are bounded even when the original polyhedron is unbounded.

\begin{proposition}\label{prop:reduction}
Fix $d$. Let $P=\{x\in\R^d:Ax\ge b\}$ be rational, with $m$ inequalities and $r=\operatorname{rank}(A)\ge1$, and let $I\subseteq\{1,\ldots,d\}$. Exact separation of a rational query point from $\SC_I(P)$ reduces to at most $2^r\binom mr$ bounded integer quadratic minimization problems. Each problem uses at most $|I|+1$ integer variables and has encoding length polynomial in the combined encoding of $A,b$ and the query point. A negative minimum supplies a separating split inequality of polynomial encoding length.
\end{proposition}

\begin{proof}
First test the query point $\bar x$ against the inequalities defining $P$. Any violated inequality separates $\bar x$ from $\SC_I(P)$. We may therefore assume $\bar x\in P$. A zero row with positive right-hand side certifies emptiness; all other zero rows may be discarded. If $I$ is empty, no further test is needed.

For each set $B$ of $r$ linearly independent rows, let $C_B=\{x:A_Bx\ge b_B\}$. By~\eqref{eq:corneridentity}, it suffices to test membership in $\SC_I(C_B)$ for every such $B$. We include all independent row sets, even when $A_Bx=b_B$ has no solution in $P$.

Fix $B$. Choose a rational right inverse $R$ of $A_B$ and define $v=Rb_B$, $L=\ker A_B$. Then
\[
C_B=v+R\R_+^r+L,\qquad
\bar s=A_B\bar x-b_B\ge0,
\]
and these data satisfy~\eqref{eq:generalcone}. The map $A_Bx-b_B$ gives the ray coordinates, which are continuous; integrality remains imposed on the original coordinates indexed by $I$.

Fix a cone $K_\sigma$ from~\eqref{eq:signcone}, and take the coordinate bounds $M_j$ from Lemma~\ref{lem:normalbound}. Use integer variables $(\pi,k)$, where $\pi$ is the split normal and $k$ is its right-hand side. Let $D$ be a positive common denominator of the coordinates of $v$. Impose
\begin{equation}\label{eq:iqpdomain}
\pi\in K_\sigma,\qquad
-M_j\le\pi_j\le M_j\ (j=1,\ldots,d),\qquad
1\le D(\pi^Tv-k)\le D-1.
\end{equation}
For integer $(\pi,k)$, the last two inequalities are equivalent to $0<f<1$, where $f=\pi^Tv-k$. They force $k=\lfloor\pi^Tv\rfloor$ and exclude both $\pi=0$ and integer values of $\pi^Tv$. If $D=1$, there are no feasible points. The real feasible region is bounded: the box bounds $\pi$, and the last two inequalities then bound $k$. The definition of $K_\sigma$ sets all coordinates of $\pi$ outside $I$ to zero. Eliminating them leaves at most $|I|+1$ integer variables.

Set $z=R^T\pi$ and define the linear forms
\[
u(\pi)=\sum_{\sigma_j=1}\bar s_jz_j,\qquad
w(\pi)=\sum_{\sigma_j=-1}\bar s_j(-z_j).
\]
With $U=u(\pi)$ and $V=w(\pi)$, identity~\eqref{eq:clearedsplit} says that the cut is violated exactly when
\begin{equation}\label{eq:iqpobjective}
F(\pi,k)=f u(\pi)+(1-f)w(\pi)-f(1-f)<0.
\end{equation}
This is a rational quadratic polynomial in $(\pi,k)$.
The coefficients of the domain and objective have polynomial encoding length. In particular, $D$ can be computed as a product of the denominators of the fixed number of entries of $v$.

Minimize~\eqref{eq:iqpobjective} over the integer points of~\eqref{eq:iqpdomain}. If the minimum is negative, a minimizer gives the separating inequality
\begin{equation}\label{eq:separator}
\alpha(\pi)^T(A_Bx-b_B)\ge1.
\end{equation}
It is valid for $\SC_I(C_B)$ and thus for $\SC_I(P)$, and is violated by $\bar x$. The bounds on $\pi,k$ and the rational cut formula imply polynomial encoding length for~\eqref{eq:separator}.

If every problem for this $B$ is integer-infeasible or has a nonnegative minimum, none of the bounded split normals gives an inequality violated by $\bar x$. Lemma~\ref{lem:normalbound} then implies $\bar x\in\SC_I(C_B)$. Repeat the test for every $B$ and use~\eqref{eq:corneridentity}. There are at most $2^r\binom mr$ quadratic problems. Constructing them and recovering a separating inequality require polynomially many rational arithmetic operations on polynomially encoded data.
\end{proof}

\begin{proof}[Proof of Theorem~\ref{thm:separation}]
For $r\ge1$, apply Proposition~\ref{prop:reduction} and the bounded integer quadratic minimization algorithm of Ari and Hildebrand~\cite[Theorem~1.2]{AH}. For fixed $d$, the number of subproblems is polynomial and their dimensions are bounded by $d+1$. Each can therefore be solved in polynomial time. If $r=0$, the input either contains an inconsistent zero row or defines $\R^d$, whose split closure is itself. These cases are immediate. The procedure returns a valid separating inequality or certifies membership, with exact rational comparisons throughout.
\end{proof}

Nonlinear formulations of split and mixed-integer rounding (MIR) separation are established tools~\cite{BS08,DGL}. The reduction here bounds the dimension and encoding of the integer quadratic subproblems, allowing the bounded algorithm of~\cite{AH} to apply; Lenstra's integer linear algorithm~\cite{L83} does not directly solve these possibly indefinite quadratic problems. A contemporaneous preprint of Wei~\cite{Wei26} also states a fixed-dimensional integer quadratic algorithm, including unbounded inputs.

When dimension varies, split separation is strongly NP-complete~\cite{CL}, and deciding split-closure emptiness is NP-hard even for polytopes in the unit hypercube~\cite{Lee19}. Structured applications have their own separation methods, such as the flip-inequality description for periodic timetabling~\cite{LM26}.

\subsection{Optimization with an integer objective}

Let $h$ be an objective with integer coefficients involving only the integer variables. To obtain a strong sequence, we must remove a current optimal LP vertex at each step, so we work on the current optimal face. Cuts on that face either produce a mixed-integer feasible vertex, completing optimization, or eliminate the face and allow us to lower the objective bound. Lemma~\ref{lem:quantlift} extends these cuts to the whole relaxation. We will choose the extensions so that each new integer objective level can be treated using the original polytope and one objective bound.

All earlier cuts are retained, and the algorithm may choose among tied optimal vertices. If both branches of a split are empty, the inconsistent inequality $0\le-1$ is a valid split cut.

In the sequential algorithms below, the input is a polytope given by rational inequalities. Write $(x,y)\in\R^p\times\R^q$, with $x$ designated integer and $y$ continuous. Every split normal has zero coefficients on $y$. Pure integer statements have $q=0$. All LP computations use exact rational arithmetic on the supplied inequality system; we do not need an inequality description of its projection.

Cook, Coullard, and Tur\'an~\cite[Theorem~1\textquotesingle\ and Proposition~3]{CCT87} bound the number of Chv\'atal--Gomory cuts needed to prove integer infeasibility in terms of dimension alone, and also control the proof size. Applied to the projection of a rational polytope $Q$ onto its $p$ integer coordinates, their result already gives a dimension-only bound on the length of a mixed-integer split refutation: each Chv\'atal--Gomory cut on that projection is a split cut and gives the same inequality in the original variables. This observation does not supply a polynomial-size inequality description of the projection when the continuous dimension varies. Our algorithm works directly with the supplied description of $Q$ and, when $Q$ is mixed-integer feasible, returns a feasible vertex of its final relaxation. Dadush and Tiwari~\cite[Theorem~4]{DT20} convert enumerative branching refutations into cutting-plane proofs.

We first extend cuts from a supporting face to the whole polytope. Cook, Coullard, and Tur\'an~\cite[Lemma~2]{CCT87} give a face-lifting argument for Chv\'atal--Gomory cuts; a corresponding rotation principle for split proofs appears in~\cite[Section~3.1]{BCDJComplexity}. Cornu\'ejols and Patil~\cite[Corollary~1]{CP26} prove $\SC_I(F)=\SC_I(R)\cap F$ for every face $F$ of a convex set $R$, including nonexposed faces. For a supplied sequence of cuts on an exposed polyhedral face, however, we need explicit coefficient bounds, and we also need to choose the lifting coefficients above prescribed lower bounds, while preserving the relaxation on the face after every cut.

Let $\sigma\ge0$ on the polytope, with $\sigma=0$ on the face $F$. We extend a cut $\ell(z)\le0$ on $F$ to $\ell(z)\le\lambda\sigma(z)$. On $F$, the two inequalities agree, and away from $F$, increasing $\lambda$ weakens the inequality. We first choose $\lambda$ large enough for validity on both split branches. The prescribed lower bounds let us also retain a supplied set separated from the face. For the optimization algorithm, that set will consist of the points at or below the next integer objective level.

\begin{lemma}\label{lem:quantlift}
Let $R$ be a nonempty rational polytope, write $z$ for all its coordinates, and let $\sigma$ be rational affine with $\sigma\ge0$ on $R$. Put $F=R\cap\{\sigma=0\}$. Suppose a supplied split sequence on $F$ adds the affine inequalities $\ell_i(z)\le0$, $i=1,\ldots,N$, each justified by a specified split disjunction in the coordinates of $R$, using only integer coordinates. For supplied rational lower bounds $E_i$, one can construct nonnegative integers $\lambda_i\ge E_i$ such that
\begin{equation}\label{eq:quantlift}
 \ell_i(z)\le\lambda_i\sigma(z)
\end{equation}
is a split cut for $R$ with the preceding lifted inequalities retained. After the first $i$ cuts, the relaxation obtained from $R$, intersected with $\sigma=0$, equals the relaxation obtained from $F$. The construction takes polynomial time and has polynomial encoding in the total supplied data, without fixing the total number of variables or $N$.
\end{lemma}

\begin{proof}
The coefficient $\lambda_i$ must make the $i$th lifted inequality valid after the preceding lifts. We obtain all validity certificates from the supplied, unlifted sequence. Their encoding lengths will therefore be bounded before any lifting coefficients are chosen.

Let $D_i^0,D_i^1$ be the two closed halfspaces of the split justifying $\ell_i\le0$. For each nonempty set $S=R\cap D_i^t$, $t\in\{0,1\}$, apply Lemma~\ref{lem:slackcertificate} with the inequalities $\ell_j\le0$, $j<i$, and $\sigma\le0$. Since $\sigma\ge0$ on $R$, these inequalities describe the corresponding split branch after the first $i-1$ cuts have been added to $F$. We obtain nonnegative rational coefficients of polynomial encoding such that
\begin{equation}\label{eq:liftcertificate}
 \ell_i\le\rho_i^t\sigma+\sum_{j<i}\nu_{ij}^t\ell_j\qquad\text{on }S.
\end{equation}
If $R\cap D_i^t$ is empty, validity on that branch is automatic. On each nonempty branch, substituting the preceding lifted inequalities $\ell_j\le\lambda_j\sigma$ gives
\[
 \ell_i\le\left(\rho_i^t+\sum_{j<i}\nu_{ij}^t\lambda_j\right)\sigma.
\]
It is therefore enough to choose $\lambda_i$ at least as large as the coefficient of $\sigma$ on each branch. Choose $C\ge1$ bounding all coefficients in~\eqref{eq:liftcertificate} and all supplied lower bounds $E_i$, and choose an integer $M\ge C+1$. Set
\begin{equation}\label{eq:liftrecurrence}
 \lambda_i=M^i,\qquad i=1,\ldots,N.
\end{equation}
Then
\[
 \rho_i^t+\sum_{j<i}\nu_{ij}^tM^j
 \le C\sum_{j=0}^{i-1}M^j<M^i.
\]
Since $\sigma\ge0$, the lifted inequality is valid on both branches. At $\sigma=0$ it agrees with $\ell_i\le0$, so induction gives the stated equality after every cut. Also, $\lambda_i\ge E_i$ and its bit length is $O(N\log M)$. Each separate application of the lemma starts with $i=1$; the exponents do not carry over from earlier applications.
\end{proof}

Write $F$ for the current optimal face and $\sigma$ for the slack in its supporting objective bound. Thus $\sigma\ge0$ on the relaxation and $F=\{\sigma=0\}$ within it.

Before applying Lemma~\ref{lem:quantlift}, test each supplied inequality by LP on $F$ with the preceding cuts retained, and omit it if it already holds there; this leaves the relaxation unchanged and preserves the validity of the remaining cuts. Stop when the relaxation obtained from $F$ first becomes empty. Compute the lifting certificates from this shortened sequence, so no certificate depends on an omitted inequality.

Every cut in the shortened sequence is violated by some point of the current relaxation on $F$. Maximize its violation on that face and choose a vertex attaining the positive maximum. This vertex is also a vertex of the whole relaxation and is optimal for the prescribed objective. Since $\sigma=0$ there, the lifted cut removes it. Thus every retained lifted cut is strong. A vertex returned by the algorithm on $F$ is likewise a vertex of the whole relaxation, even when the relaxation is lower-dimensional.

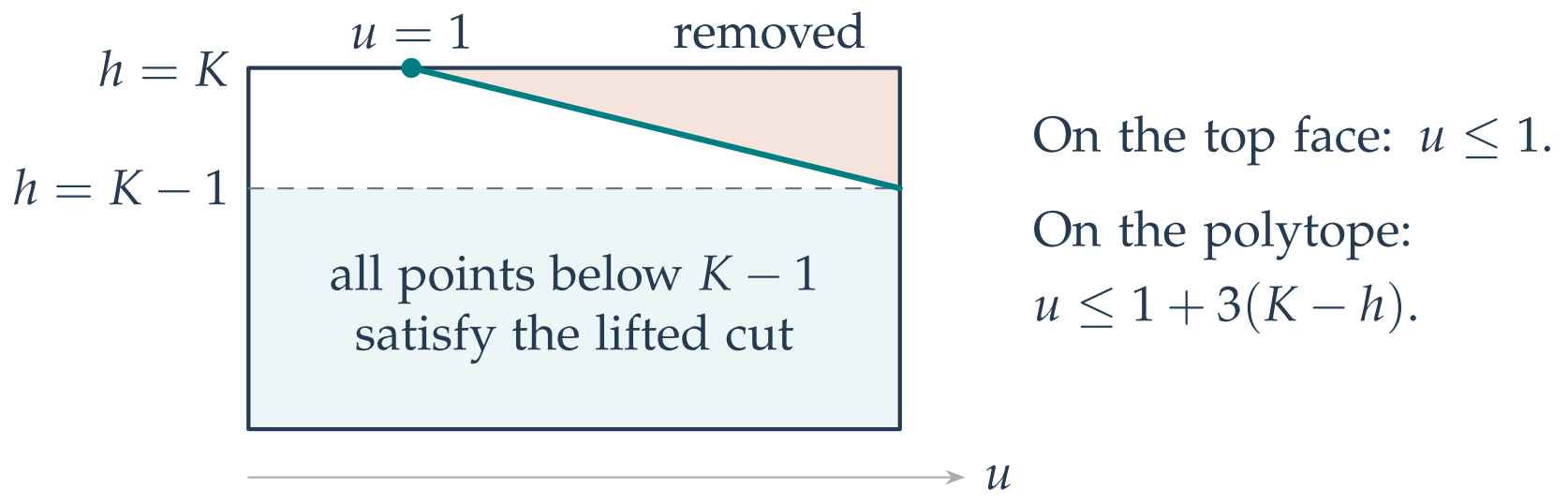
\begin{figure}[htbp]
\centering
\begin{adjustbox}{max width=\linewidth}
\begin{tikzpicture}[x=1.15cm,y=.85cm,figure label]
\fill[cutteal!7] (0,0) rectangle (4,2);
\fill[warmrust!16] (1,3)--(4,3)--(4,2)--cycle;
\draw[boundary] (0,0) rectangle (4,3);
\draw[dashed,inkblue!70] (0,2)--(4,2);
\draw[cut] (1,3)--(4,2);
\fill[cutteal] (1,3) circle (2pt);
\node[left] at (0,3) {$h=K$};
\node[left] at (0,2) {$h=K-1$};
\node[above] at (1,3) {$u=1$};
\node[above] at (3.2,3) {removed};
\node[align=center] at (2,1) {all points below $K-1$\\satisfy the lifted cut};
\draw[axis] (0,-.4)--(4.4,-.4) node[right] {$u$};
\node[align=left,anchor=west,text width=4.7cm] at (4.7,1.7) {On the top face: $u\le1$.\\[7pt]On the polytope:\\[3pt]$u\le1+3(K-h)$.};
\end{tikzpicture}
\end{adjustbox}
\caption{The cut on $h=K$ is $u\le1$, and the original polytope satisfies $0\le u\le4$. Once validity on the split branches has been established, choose the lifting coefficient to be at least three. The lifted cut still gives $u\le1$ on $h=K$, while every original point with $h\le K-1$ satisfies it.}
\label{fig:levellift}
\end{figure}

Eliminating the current optimal face only makes the objective maximum strictly smaller; it need not decrease it by a full integer. When the eliminated face has integer objective value $K$, the split $h\le K-1$ or $h\ge K$ permits $h\le K-1$. We choose every lift to retain all original points with $h\le K-1$. After the face is eliminated and this split is applied, the relaxation is exactly the original polytope intersected with the new objective bound. Figure~\ref{fig:levellift} illustrates this choice. The next lemma gives the cut count and encoding bound.

\begin{lemma}\label{lem:objectivelevel}\label{par:objectivelevel}
Let $P\subseteq\R^p\times\R^q$ be a rational polytope, with $x$ designated integer and $y$ continuous. Let $h(x,y)=c^Tx+\beta$, where $c\in\Z^p\setminus\{0\}$ and $\beta\in\Z$. Fix $K\in\Z$ and put $P_K=P\cap\{h\le K\}$ and $F_K=P\cap\{h=K\}$. Suppose a supplied weak split sequence of at most $k$ cuts on $F_K$ returns a mixed-integer feasible vertex of its final relaxation or proves the face empty. Then a strong split-cut algorithm on $P_K$ either returns that vertex as an optimal vertex for $h$, leaves exactly $P_{K-1}$, or proves emptiness. It uses at most $k$ cuts when a vertex is returned and at most $k+1$ otherwise. Its time and encoding are polynomial in $P,h,K$ and the supplied sequence on $F_K$.
\end{lemma}
\begin{proof}
If $P_K$ is empty, no cuts are needed. Otherwise omit inequalities already implied on the face and stop when its relaxation first becomes empty, as above. For each remaining inequality $\ell_i\le0$ on $F_K$, apply Lemma~\ref{lem:quantlift} with $R=P_K$ and $\sigma=K-h$, requiring
\begin{equation}\label{eq:inactivebound}
 \lambda_i\ge E_i:=\max\{0,\max_{z\in P}\ell_i(z)\}.
\end{equation}
Each lifted cut removes an optimal vertex of the current relaxation, as proved above. If the algorithm on $F_K$ returns a mixed-integer feasible vertex, the same point is a vertex of the whole relaxation and attains its maximum $h=K$.

Suppose instead that the face becomes empty. Every point of the original set $P_{K-1}$ satisfies every lifted inequality: if $h(z)\le K-1$, then
\[
 \ell_i(z)\le E_i\le\lambda_i\le\lambda_i(K-h(z)).
\]
If the whole relaxation is nonempty, compactness gives $\max h<K$. The split $h\le K-1$ or $h\ge K$ then permits the cut $h\le K-1$. Add it only if it is violated; in that case it removes every current optimal vertex. The resulting relaxation is exactly
\begin{equation}\label{eq:originalcap}
 P\cap\{h\le K-1\}.
\end{equation}
The new objective bound gives one inclusion. The other follows because every point of $P_{K-1}$ satisfies all the lifted cuts. If $F_K$ was empty initially, only the cut $h\le K-1$ may be needed.

There are at most $k$ lifted cuts and one further cut $h\le K-1$. All validity certificates and lower bounds use the supplied sequence on $F_K$ and the original $P$, so the time and encoding are polynomial in those data. All earlier cuts remain retained. Equality~\eqref{eq:originalcap} allows the next application to use the original input and the new objective level.
\end{proof}

\paragraph{Finding a feasible vertex.}

The feasibility algorithm used in Lemma~\ref{lem:objectivelevel} is free to choose its objective. A direction of approximately minimum lattice width limits the number of integer levels to be examined, and fixing a level reduces the number of integer coordinates. The recursion must return a mixed-integer feasible vertex or prove infeasibility.

For a full-dimensional compact convex set $T\subseteq\R^d$, write
\[
 w_a(T)=\max_T a^Tx-\min_T a^Tx,
 \qquad w(T)=\min_{a\in\Z^d\setminus\{0\}}w_a(T).
\]
Choose nondecreasing computable flatness bounds $f_d=O(d^{3/2})$~\cite{BLPS99} such that $w(T)>f_d$ implies $T\cap\Z^d\ne\varnothing$. We use only that $f_d$ depends on $d$. The next lemma explains why it suffices to examine a bounded number of levels near the maximum in a direction of approximately minimum width. The minimum-width version appears in Eisenbrand--Shmonin~\cite{ES08}. We use a strict inequality in the choice of $D$ so that the argument does not depend on whether equality is included in the flatness bound.

\begin{lemma}\label{lem:smallcap}
Suppose $a\in\Z^d\setminus\{0\}$ satisfies $w_a(T)\le\kappa w(T)$, and fix $D>\kappa f_d$. Let $M=\max_T a^Tx$. If $T$ contains an integer point, then it contains one with $a^Tx\ge M-D$. If $T$ contains no integer point, then $w_a(T)<D$.
\end{lemma}

\begin{proof}
Set $W=w_a(T)$. If $W\le D$, every point of $T$ satisfies $a^Tx\ge M-D$, so the first assertion is immediate. If $W>D$, choose $v\in T$ with $a^Tv=M$ and set
\[
 \alpha=D/W,\qquad T'=v+\alpha(T-v).
\]
Convexity gives $T'\subseteq T$, and $a^Tx\ge M-D$ throughout $T'$. Translation invariance and homogeneity of lattice width yield
\[
 w(T')=\alpha w(T)\ge D/\kappa>f_d.
\]
Thus $T'$ contains an integer point. Finally, if $T$ has no integer point, flatness gives $w(T)\le f_d$, hence $W\le\kappa f_d<D$.
\end{proof}

\begin{figure}[htbp]
\centering
\begin{adjustbox}{max width=\linewidth}
\begin{tikzpicture}[x=1cm,y=1cm,figure label]
\fill[inkblue!6] (0,0)--(4,0)--(2,3)--cycle;
\fill[cutteal!22] (4/3,2)--(8/3,2)--(2,3)--cycle;
\draw[boundary] (0,0)--(4,0)--(2,3)--cycle;
\draw[cut] (4/3,2)--(8/3,2)--(2,3)--cycle;
\draw[dashed,warmrust] (-.4,2)--(4.4,2);
\draw[dotted,inkblue!60] (2,3)--(4.4,3);
\fill[inkblue] (2,3) circle(2pt) node[above] {$v$};
\node at (2,.8) {$T$};
\node[text=cutteal] at (2,2.35) {$T'$};
\node[right] at (4.4,3) {$M$};
\node[right] at (4.4,2) {$M-D$};
\draw[axis] (-.6,-.1)--(-.6,3.3) node[above] {$a^Tx$};

\end{tikzpicture}
\end{adjustbox}
\caption{The set $T'$ is obtained by scaling $T$ towards a maximizer $v$ by the factor $D/W$. It lies between $a^Tx=M-D$ and $a^Tx=M$. The triangle illustrates the construction; the proof applies to every full-dimensional
compact convex set.}
\label{fig:smallcap}
\end{figure}
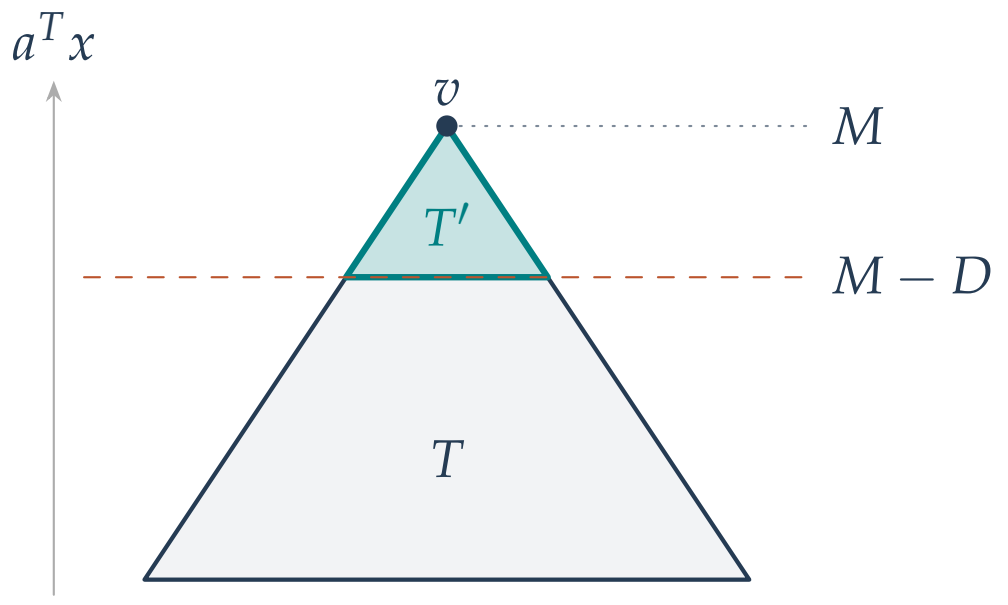

To apply the lemma with continuous variables present, we need a direction for the projection onto the integer coordinates. An inequality description of that projection may be large, but maximizing a linear function over the projection is just an LP over the original polytope. We use these LPs to construct an inscribed simplex, then apply lattice reduction to find a direction of approximately minimum width~\cite{L83,LLL82}. The lattice computations use only the integer coordinates.

The edge vectors of the simplex define a lattice in which vector length bounds directional width. A short vector in this lattice will therefore give the required integer direction.

\begin{lemma}\label{lem:widthdirection}
Fix $d\ge1$. Let $Q\subseteq\R^d\times\R^q$ be a rational polytope given by inequalities, and suppose that its projection $T$ onto the $d$ integer coordinates is full-dimensional. In time polynomial in the full input length, one can find a primitive integer vector $a$ of polynomial encoding such that
\[
 w_a(T)\le\kappa_d w(T),
\]
where $\kappa_d$ depends only on $d$. This requires only LP computations over $Q$ and lattice reduction in dimension $d$; $q$ is unrestricted.
\end{lemma}

\begin{proof}
Let $L$ be the full input length. Obtain $d+1$ affinely independent points $v_0,\ldots,v_d$ of $T$ by LP. At each stage, optimize both signs of a basis of linear functions constant on the affine span of the points already obtained. If none gives a point outside that span, it is the affine hull of $T$. Since $T$ is full-dimensional, the process produces the required points. Each optimizer can be chosen as a vertex of $Q$: restrict to its optimal face and successively optimize the $d+q$ coordinates until a singleton remains. This uses at most $d+q$ additional LPs. Every selected point is a vertex of the same original $Q$, so rational vertex bounds give polynomial encoding even when $q$ varies. Its projection need not be a vertex of $T$.

Apply the simplex construction in Lenstra~\cite[Section 2]{L83}. Let $\lambda_0,\ldots,\lambda_d$ be its barycentric affine functions, so that $x=\sum_{i=0}^d\lambda_i(x)v_i$ and $\sum_{i=0}^d\lambda_i(x)=1$. Maximize both signs of each $\lambda_i$ over $T$. If one maximum exceeds two, replace $v_i$ by a point attaining that maximum, chosen as the projection of a vertex of $Q$ as above. The simplex volume is multiplied by $|\lambda_i|>2$. Every retained point has polynomial encoding because it is the projection of a vertex of the original $Q$. Determinant bounds put each nonzero simplex volume between $2^{-\operatorname{poly}(L)}$ and $2^{\operatorname{poly}(L)}$. Since each replacement more than doubles the volume, there are only polynomially many replacements. At termination, $|\lambda_i(x)|\le2$ for every $x\in T$ and every $i$.

We now use this simplex to compare widths with lattice vector lengths. Let $B$ have columns $v_i-v_0$, $i=1,\ldots,d$, and set $u=B^Ta$ for $a\in\Z^d$. Since the simplex lies in $T$, its width is a lower bound for $w_a(T)$. For the upper bound, write each point of $T$ using its barycentric coordinates, whose absolute values are at most two. This gives
\[
 \frac{\|u\|_2}{\sqrt d}
 \le \max(0,u_1,\ldots,u_d)-\min(0,u_1,\ldots,u_d)
 \le w_a(T)\le4\|u\|_1\le4\sqrt d\,\|u\|_2.
\]
Let $\mu$ be the shortest nonzero vector length in the rational lattice $B^T\Z^d$. Lenstra--Lenstra--Lov\'asz (LLL) lattice reduction, after clearing denominators, produces $a\ne0$ with $\|B^Ta\|_2\le\rho_d\mu$, where $\rho_d=2^{(d-1)/2}$ suffices~\cite[Propositions 1.11 and 1.26]{LLL82}. The displayed lower bound, applied to a minimum-width integer direction, gives $\mu\le\sqrt d\,w(T)$, so
\[
 w_a(T)\le4d\rho_d\,w(T).
\]
Dividing $a$ by the gcd of its entries only improves the inequality. Thus one may take $\kappa_d=4d\rho_d$, or any larger dimension-dependent constant. Recover $a$ from the returned lattice vector using the rational inverse of $B^T$; both have polynomial encoding. All LPs, determinant computations, and lattice reductions use polynomially encoded data. No facets of $T$ are enumerated.
\end{proof}

\paragraph{The feasibility algorithm.}\label{subsec:vertex}

We combine the width bound with face lifting. Converting this recursion into cuts is classical for infeasibility proofs~\cite{CCT87,DT20}. For optimization, a successful recursive call must also return a vertex of the whole relaxation. Processing each level as a supporting face ensures this, by Lemma~\ref{lem:objectivelevel}.

\begin{proof}[Proof of Proposition~\ref{prop:feasiblevertex}]
We induct on the number $p$ of integer coordinates and construct the cut bound $T_p$. If the LP relaxation is empty, no cuts are needed. If $p=0$, return an LP vertex, using lexicographic LP optimization if necessary. Boundedness guarantees a vertex whenever the relaxation is nonempty. Thus $T_0=0$.

\emph{Reduction to a full-dimensional projection.} First determine the affine hull of the projection of $Q$ onto its integer coordinates using the LP construction in Lemma~\ref{lem:widthdirection}. Write this hull as $v_0+\operatorname{span}_{\R}(W)$, where $W$ is an integral $p\times r$ matrix of rank $r$, obtained by clearing denominators in independent differences of the projected points. If $r<p$, Hermite normal form~\cite{KB} yields a polynomially encoded unimodular $U$ with
\[
 UW=\begin{pmatrix}H\\0\end{pmatrix},\qquad H\in\Z^{r\times r}\text{ nonsingular}.
\]
Thus the coordinates $u=Ux$ satisfy $u_{r+1:p}=\beta:=(Uv_0)_{r+1:p}$ throughout $Q$~\cite[Section 2]{L83}. For $r=0$, take $U=I$. If some $\beta_j$ is noninteger, the split $u_j\le\lfloor\beta_j\rfloor$ or $u_j\ge\lfloor\beta_j\rfloor+1$ has both branches empty and permits $0\le-1$. Otherwise substitute
\[
 x=U^{-1}(u_1,\ldots,u_r,\beta)^T
\]
and apply induction with $r$ integer coordinates, leaving $y$ unchanged. This is an affine bijection onto the affine subspace containing $Q$, so it preserves vertices. Since $U$ and $U^{-1}$ are integral and $\beta$ is integral, it preserves mixed-integer feasibility in both directions. Extend a reduced split normal by zero entries in the fixed coordinates to obtain $\bar\pi$. Its normal in the original coordinates is $U^T\bar\pi$, which remains nonzero and integral. The right-hand sides of the split disjunctions remain integral. All transformations have polynomial encoding.

\emph{The recursive step.} We may now assume that the projection $T$ is full-dimensional in $\R^p$. Obtain $a$ from Lemma~\ref{lem:widthdirection}. Fix an integer $D_p>\kappa_p f_p$ and let $M=\max_Q a^Tx$. If $M$ is noninteger, apply the split $a^Tx\le\lfloor M\rfloor$ or $a^Tx\ge\lfloor M\rfloor+1$ to add $a^Tx\le\lfloor M\rfloor$; its upper branch is empty. Stop if the relaxation becomes empty. Otherwise begin at $K=\lfloor M\rfloor$ and consider integer levels in decreasing order.

At level $K$, the current relaxation is $Q\cap\{a^Tx\le K\}$, and if $F_K=Q\cap\{a^Tx=K\}$ is nonempty, it is the optimal face. Since $a$ is primitive, choose polynomially encoded unimodular coordinates $x=Vu$ with $a^TV=e_1^T$. On $F_K$, the coordinate $u_1=K$ is fixed, leaving at most $p-1$ integer coordinates. Apply induction to $F_K$ and express its cuts and splits in the original coordinates. Lemma~\ref{lem:objectivelevel} then either returns a mixed-integer feasible vertex of the current relaxation, proves emptiness, or leaves exactly $Q\cap\{a^Tx\le K-1\}$. It uses at most $T_{p-1}+1$ cuts and also covers an initially empty face. Thus we can proceed to the next level using the original $Q$ and the new bound, with all earlier cuts retained.

The number of levels considered depends only on $p$. If $w_a(T)\le D_p$, decreasing $K$ from $\lfloor M\rfloor$ reaches below the original minimum after at most $D_p+1$ levels. The relaxation is then empty unless a feasible vertex was returned earlier. If $w_a(T)>D_p$, Lemma~\ref{lem:smallcap} gives an integer point of $T$ with $a^Tx\ge M-D_p$. This point has a mixed-integer feasible preimage in $Q$ on one of the first $D_p+1$ levels. Every split cut retains that point, so the recursion cannot prove its level empty and must return a feasible vertex by that level. If $Q$ is mixed-integer infeasible, only the first case is possible. We may use $B_p=D_p+2$ as a common bound on the number of levels and enlarge these constants to make them nondecreasing in $p$. Then
\begin{equation}\label{eq:vertexrecurrence}
 T_0=0,\qquad T_p=1+B_p(T_{p-1}+1)
\end{equation}
is a valid coarse bound, also covering lower-dimensional reductions.

More explicitly, with $\kappa_p=4p\,2^{(p-1)/2}$, we may take
\[
D_p=\lfloor\kappa_p f_p\rfloor+1,\qquad B_p=D_p+2.
\]
These constants are nondecreasing. The recurrence gives
\[
T_p+1\le\prod_{j=1}^p(B_j+2).
\]
Since $f_j=O(j^{3/2})$, we have $\log(B_j+2)=O(j+\log(j+1))$, and hence $T_p=2^{O(p^2)}$. This bounds the number of cuts; the running-time polynomial still depends on $p$.

For fixed $p$, the number of recursive calls and cuts is bounded independently of the input length. Each call uses rational LP, unimodular coordinate changes, lattice reduction, and Lemma~\ref{lem:quantlift}. There are only a bounded number of compositions of these constructions, each with polynomial encoding, so the total encoding and running time are polynomial in the full input length. This remains true when the continuous dimension varies: the LPs use the supplied variables and inequalities, while the lattice operations have dimension at most $p$.
\end{proof}

\paragraph{The prescribed objective.}

The feasibility algorithm can choose its direction to limit the number of levels it examines, but the prescribed objective, to which we now return, may have a large range. At each integer objective level, Proposition~\ref{prop:feasiblevertex} supplies a short sequence of cuts on the face at that level. Lemma~\ref{lem:objectivelevel} extends them to strong cuts on the full relaxation. The number of levels needed is determined by the objective gap.

All cuts remain in the proof. For LP computations, however, we use the original inequalities, the current objective bound, and the cuts added at the current level. Lemma~\ref{lem:objectivelevel} shows that all cuts from earlier levels are redundant in this system. We therefore avoid repeatedly processing a growing list of inequalities.

\begin{proof}[Proof of Theorem~\ref{thm:gapalgorithm}]
Let $h=c^Tx$ and $M=\max_P h$. If $M$ is noninteger, the split $h\le\lfloor M\rfloor$ or $h\ge\lfloor M\rfloor+1$ permits the cut $h\le\lfloor M\rfloor$. Consider integer levels $K$ in decreasing order. At the start of each level, the current relaxation is exactly
\[
 P_K=P\cap\{h\le K\}.
\]
All earlier cuts are retained, but they are redundant on $P_K$. Since $c$ is primitive, choose polynomially encoded unimodular coordinates in which $h$ is one coordinate. The face $F_K=P\cap\{h=K\}$ then has at most $p-1$ free integer coordinates.

Apply Proposition~\ref{prop:feasiblevertex} to $F_K$ and express its cuts and splits in the original coordinates. It supplies at most $T_{p-1}$ weak split cuts and either returns a mixed-integer feasible vertex or proves the face empty. Apply Lemma~\ref{lem:objectivelevel} to this sequence. If a vertex is returned, it is optimal for the current relaxation. It is also optimal for the original mixed-integer problem, since every cut retains all mixed-integer feasible points. Otherwise, at most $T_{p-1}+1$ strong cuts leave exactly $P_{K-1}$, establishing the assertion for the next level. Under the feasibility assumption of the theorem, the whole relaxation cannot become empty.

The levels above $z^*$ are unsuccessful, and the level $z^*$ returns an optimal vertex. There are exactly $G$ unsuccessful levels. Including at most one initial rounding cut, the number of cuts is therefore at most
\begin{equation}\label{eq:gapcount}
 1+G(T_{p-1}+1)+T_{p-1}=(G+1)(T_{p-1}+1).
\end{equation}
Each level $K$ has polynomial encoding in $L$, because it lies in the objective range of the original rational polytope; $L$ includes the objective vector. The feasibility algorithm is always applied to $P\cap\{h=K\}$, described using the original $P$ and one equation. It uses a number of cuts depending only on $p$ and has encoding polynomial in $L$, as does any vertex it returns. The lower bounds in~\eqref{eq:inactivebound} are LP values over the original $P$. Lemma~\ref{lem:quantlift} therefore gives a $\operatorname{poly}_p(L)$ bound on every generated coefficient, independently of $G$.

At each level, solve LPs using the original inequalities, the current objective bound, and the cuts added at that level. Append every cut and its split disjunction to the retained proof. By Lemma~\ref{lem:objectivelevel}, this LP system describes the same relaxation as the full proof after every cut. Processing a level, including selecting the optimal vertices to be removed, therefore takes $\operatorname{poly}_p(L)$ time without examining cuts from earlier levels. There are $G+1$ levels, giving the stated running time.
\end{proof}

If the original LP relaxation is empty, no cuts are needed. Otherwise, if its mixed-integer set is empty, the same algorithm proves infeasibility at every integer level in its objective range. With
\[
 W=\max\{0,\lfloor\max_P h\rfloor-\lceil\min_P h\rceil+1\},
\]
the corresponding cut and time bounds are $O_p(W+1)$ and $(W+1)\operatorname{poly}_p(L)$. Test the current relaxation for emptiness before computing its objective maximum. The final objective bound may be one below the original minimum, but still has polynomial encoding.

For a nonprimitive integral objective, divide its coefficients by their gcd before measuring the gap. An objective with rational coefficients on the integer variables can first be scaled to integral coefficients and normalized in the same way, with polynomial encoding. If the objective is zero, apply Proposition~\ref{prop:feasiblevertex} directly and omit implied inequalities as above. Every remaining cut is strong because every vertex is optimal.

Consequently, for fixed integer dimension, a polynomial bound on the normalized gap $G$ in terms of the full input length gives a polynomial-time strong split-cut algorithm. It terminates at a mixed-integer feasible vertex that is optimal for both the final LP relaxation and the original problem. For infeasible inputs, the corresponding conclusion holds when the original relaxation spans only polynomially many integer objective levels.

Quantitative cutting-plane bounds are classical. The lexicographic analysis of Nourie and Venta~\cite{NV82}, as presented by Nemhauser and Wolsey~\cite[Chapter~II.4, Theorem~3.8]{NW88} and Balas, Fischetti, and Zanette~\cite[Section~4]{BFZ10}, counts integer vectors within variable and objective ranges. Neto~\cite[Theorem~3.4]{Neto12} gives $O(d^{n+1}(1+\|c\|_1))$ iterations for a bounded pure integer program in $[0,d]^n$ with integral objective vector $c$. These bounds depend on numerical variable ranges; they do not give the cut bound depending only on $p$ and $G$ above. Finite mixed-integer Gomory convergence for integer-valued objectives is also known~\cite[Section~2]{GK13}.

Owen and Mehrotra~\cite[Theorem~4]{OM01} prove finite $\varepsilon$-optimal convergence for bounded mixed-integer programs. Their algorithm generates cuts at all LP vertices satisfying a near-optimality criterion. For an integer-valued objective, choosing $\varepsilon<1$ gives exact optimality. Their theorem does not give our cut bound or require every cut to remove a current optimal LP vertex. Other finite algorithms use different cut families, including lex-inequalities~\cite{CDDR21}; see also the cut-selection survey~\cite{DM18}. Theorem~\ref{thm:gapalgorithm} gives a bound in terms of the objective gap, requires every cut to be strong, controls running time and coefficient encoding, and returns a mixed-integer feasible vertex of the final LP relaxation.

A logarithmic dependence on the objective gap is necessary in general. For integer $H\ge2$, let
\[
P_H=\conv\{(0,0,0),(2,0,0),(0,2,0),(1/2,1/2,H)\},
\]
with all three variables designated integer and objective $\max x_3$. The integer optimum is zero and the LP optimum is $H$, so $G=H$. Conforti, Del Pia, Di Summa, Faenza, and Grappe~\cite[Lemma~19]{CDDFG15} give a split-rank lower bound $\Omega(\log H)$; see also~\cite[Remark~2.10]{BCDJComplexity}. A split proof of length $k$ gives rank at most $k$, so every split proof of $x_3\le0$, including every strong sequence, has length $\Omega(\log G)$. This lower bound is logarithmic in the gap and does not rule out an algorithm polynomial in the binary input length.

\paragraph{Objectives involving continuous variables.}\label{subsec:mip}

The restriction on the objective cannot simply be dropped. A rescaled version of the classical example of Cook, Kannan, and Schrijver~\cite{CKS} maximizes a continuous variable $z$ over
\begin{equation}\label{eq:ckscontinuous}
z\ge0,\qquad z\le2x_1,\qquad z\le2x_2,\qquad x_1+x_2+z\le2,
\end{equation}
with $x_1,x_2$ integer. This tetrahedron has LP optimum $1$, attained at $(1/2,1/2,1)$, and mixed-integer optimum $0$: $z>0$ would require $x_1,x_2>0$ and $x_1+x_2<2$. Yet no finite sequence of split cuts establishes $z\le0$.

To see why sequential split cuts cannot establish $z\le0$, let $T=\conv\{(0,0),(2,0),(0,2)\}$. Every relaxation retains $T\times\{0\}$, because its vertices are mixed-integer feasible. Suppose a point with $z>0$ also remains. For any split $a^Tx\le k$ or $a^Tx\ge k+1$, with $a\in\Z^2\setminus\{0\}$, the values at the three base vertices are $0,2a_1,2a_2$. Their range has length at least two, so some base vertex lies strictly outside the closed strip $k\le a^Tx\le k+1$. Points sufficiently close to that vertex on the segment joining it to the point with $z>0$ still have $z>0$ and belong to one split branch. They survive every cut justified by this split. Induction shows that a point with $z>0$ remains after any finite sequence.

Taking the product with an integer interval $0\le r\le1$ gives the same obstruction with three integer variables and one continuous variable. On the face $r=0$, a split in $(x,r)$ either restricts to a split on $x$ or has the entire face in one branch, so the same argument applies. Thus even a gap of one does not permit an extension to arbitrary continuous objectives. The obstruction is exact termination: a small positive value of a continuous objective cannot be rounded down to zero.

There is no contradiction with Proposition~\ref{prop:feasiblevertex}. That proposition gives a short refutation when the mixed-integer feasible set is empty. Here it is nonempty, and the obstruction is proving the exact continuous objective bound $z\le0$. Refuting a relaxation with the additional constraint $z\ge\varepsilon$ for each fixed $\varepsilon>0$ does not give a finite split proof of that bound on the original relaxation.

\subsection{Strong split proofs}\label{sec:proofconversion}

The preceding algorithm constructs the cuts needed at each objective level. Suppose instead that a weak split proof of a desired bound is already supplied. On an integer objective level above that bound, the supplied proof ends with an empty relaxation, so we can apply Lemma~\ref{lem:objectivelevel} to turn its cuts into strong cuts. The same supplied proof may be used at several levels, so the resulting bound depends on the number of levels above the target.

\begin{proposition}\label{prop:gapconversion}
Let $P\subseteq\R^n$ be a nonempty rational polytope with all coordinates designated integer, let $c\in\Z^n\setminus\{0\}$, and let $\tau\in\Q$. Suppose a supplied weak split proof of length $k$ establishes $c^Tx\le\tau$, or emptiness. Set
\[
M=\max_{x\in P}c^Tx,\qquad G=\max\{0,\lfloor M\rfloor-\lfloor\tau\rfloor\}.
\]
There is a strong split-cut algorithm establishing $c^Tx\le\lfloor\tau\rfloor$, or emptiness, with at most $1+G(k+1)$ cuts. Its running time is $(G+1)\operatorname{poly}(S)$, and each generated inequality and split has polynomial encoding in $S$, the total encoding of $P,c,\tau$ and the supplied proof. The constructed sequence may also use splits $c^Tx\le K-1$ or $c^Tx\ge K$, with $K\in\Z$.
\end{proposition}

\begin{proof}
Write $h=c^Tx$. If needed, use the split $h\le\lfloor M\rfloor$ or $h\ge\lfloor M\rfloor+1$ to add $h\le\lfloor M\rfloor$; this removes every current optimal vertex. For each integer level $K>\lfloor\tau\rfloor$, apply the supplied sequence of cuts to $F_K=P\cap\{h=K\}$. Each cut remains valid on its split branches after restriction to the face. The resulting relaxation on $F_K$ is empty because it satisfies $h\le\tau<K$, or because the supplied proof already establishes emptiness. Thus the supplied sequence gives a split proof of infeasibility on $F_K$ with at most $k$ cuts.

By Lemma~\ref{lem:objectivelevel}, at most $k$ lifted strong cuts and the cut $h\le K-1$ leave exactly $P\cap\{h\le K-1\}$, or prove emptiness. An initially empty face requires no lifted cuts. After at most $G$ levels, the relaxation satisfies $h\le\lfloor\tau\rfloor$ unless it has already become empty. Including the initial rounding cut gives the bound $1+G(k+1)$.

Every level $K$ has polynomial encoding in $S$, as does the supplied proof restricted to $F_K$. At each level, Lemma~\ref{lem:objectivelevel} uses the original $P$ and this restricted proof. Its running time and the encoding of every generated cut and split are therefore polynomial in $S$, independently of the cuts from earlier levels. Processing at most $G$ levels, together with the initial LP and possible rounding cut, gives the stated running time.
\end{proof}

In Lemma~\ref{lem:quantlift}, the lifting coefficients have the form $M^i$, with $i\le k$ for each restricted proof, and therefore have $O(k\log M)$ bits. The polynomial bound on $\log M$ is obtained from the unlifted data at that level; the exponents restart at each new level.

The numerical-gap dependence remains: the supplied proof may have to be applied at every integer level above the target. Thus a short weak split proof alone does not give a strong algorithm with running time polynomial in the binary input length.

\subsection{Optimization in a split strip}\label{sec:interleaving}

The dependence on the objective gap arose from processing integer objective levels one at a time. When the relaxation lies in a closed split strip, every integer point belongs to one of its two boundary faces. We can then combine algorithms for these faces. Each lifted cut must preserve the opposite face, so that both supplied sequences remain applicable. A cut bounding the objective by the larger of the two current face maxima is valid on both split branches. After adding this cut, a face attaining that value is optimal for the whole relaxation. We can take its next strong cut without examining intervening objective levels.

\begin{proposition}\label{prop:interleave}
Let $P\subseteq\R^n$ be a rational polytope with all coordinates designated integer, let $h$ be a rational affine objective, and let $s$ be an affine function with integer coefficients satisfying $0\le s\le1$ on $P$. For $i=0,1$, suppose a supplied strong split-cut sequence optimizes $h$ over the integer points of $F_i=P\cap\{s=i\}$, or proves infeasibility. Let its length be $N_i$, with $N_i=0$ for an initially empty face, and assume all its splits extend to split disjunctions with integer coefficients in the coordinates of $P$. Then at most
\[
2(N_0+N_1)+1
\]
strong split cuts optimize $h$ over $P\cap\Z^n$, or prove infeasibility. The construction has polynomial running time and encoding in the total supplied data and retains every cut.
\end{proposition}

\begin{proof}
Let $F_i^{r_i}$ be the relaxation on face $i$ after the first $r_i$ cuts of its supplied sequence. To preserve the relaxation on the other face, each lifted cut must agree with its supplied cut on one face and hold throughout the opposite face.

If $P$ is empty, no cuts are needed. If $s$ has zero linear part, its constant value is the integer zero or one; the corresponding face is $P$, and its supplied sequence proves the assertion. Otherwise lift the sequence on $F_0$ using $\sigma=s$, and the sequence on $F_1$ using $\sigma=1-s$. For each inequality $a^Tx\le b$ on $F_0$, also require the lifting coefficient to be at least $\max\{0,\max_{x\in F_1}(a^Tx-b)\}$, omitting this requirement when $F_1$ is empty. Since $s=1$ on $F_1$, every point of the original $F_1$ then satisfies the lifted inequality. Impose the corresponding requirement for cuts on $F_1$. Lemma~\ref{lem:quantlift} gives polynomially encoded lifts from the two supplied sequences. We may interleave the sequences: cuts from the other face only restrict the split branches on which each lifted inequality must be valid.

Maintain the invariant
\[
Q\cap\{s=i\}=F_i^{r_i}\qquad(i=0,1)
\]
for the current relaxation $Q$. If both $F_0^{r_0}$ and $F_1^{r_1}$ are empty, the split $s\le0$ or $s\ge1$ proves infeasibility in one cut, unless $Q$ is already empty. Otherwise let $B$ be the larger of the objective maxima on the nonempty sets $F_i^{r_i}$. Because $0\le s\le1$ on $Q$, the two branches of this split are exactly its two boundary faces. Hence $h\le B$ is a split cut. Add it only if $\max_Q h>B$; it then removes every current optimal vertex and leaves both $F_i^{r_i}$ unchanged.

Choose a face whose current maximum is $B$. If its supplied algorithm has terminated at an integer optimum, return that point. It is optimal for $Q$ and hence for $P\cap\Z^n$. Otherwise add the next lifted cut from that face's sequence. The optimal vertex removed by the supplied cut is also an optimal vertex of $Q$, so the lifted cut is strong. It advances that face's sequence by one cut and leaves the other face unchanged. Earlier objective bounds remain valid on the smaller sets $F_i^{r_i}$, so the displayed equalities continue to hold.

There are at most $N_0+N_1$ lifted cuts from the two supplied sequences. Before each, we add at most one cut of the form $h\le B$. At the end, at most one further such cut or a cut proving infeasibility is needed. This gives $2(N_0+N_1)+1$. Each value $B$ is computed on a relaxation defined by one of the supplied sequences on $F_0$ or $F_1$, so it has polynomial encoding in the supplied data. Together with Lemma~\ref{lem:quantlift}, this gives the claimed running time and coefficient bounds. If $h$ is constant on a face, its supplied algorithm may use an auxiliary objective to select vertices; all of them remain optimal for $h$.
\end{proof}

Corollary~\ref{cor:capstrip} assumes only that the points violating $h\le\tau$ lie in a split strip. The whole polytope need not lie there. To use the preceding construction, we must show that every split branch containing a violating point is a supporting face of the original polytope. The strip condition gives precisely this property.

We use the planar algorithm~\cite[Section~2]{BCDJ} with its prescribed choice among tied optimal vertices. A rational affine objective is reduced to an integer linear objective by subtracting its constant term and multiplying by a positive common denominator; this preserves its optimal vertices. If the objective is constant on a face, any auxiliary linear objective may select vertices there, since every vertex remains optimal for the original objective.

The planar algorithm also handles segments and points. On a rational line containing no lattice point, an integer normal with noninteger constant value gives an empty split disjunction. Otherwise a unimodular coordinate change identifies the lattice points on the line with the integers in one coordinate. Rounding the endpoint that maximizes the objective either gives an integer optimal endpoint or proves emptiness. For a noninteger singleton, a split on a noninteger coordinate proves emptiness. Each cut can therefore be chosen to remove a current optimal vertex.

When proving $h\le\tau$, stop as soon as the current LP maximum is at most $\tau$, or the relaxation is empty. This cannot increase the number of cuts. On any face whose integer points satisfy the target, the planar algorithm must reach one of these stopping conditions.

\begin{proof}[Proof of Corollary~\ref{cor:capstrip}]
If $P$ is empty or already satisfies $h\le\tau$, no cuts are needed. Write $\ell_0=\pi^Tx-k$ and $\ell_1=k+1-\pi^Tx$. Let $J\subseteq\{0,1\}$ index the branches $P\cap\{\ell_i\le0\}$ containing a point $u$ with $h(u)>\tau$. For $i\in J$, condition~\eqref{eq:originalstrip} forces $\ell_i(u)=0$. Suppose some $y\in P$ had $\ell_i(y)<0$. Points on the segment from $u$ to $y$ sufficiently close to $u$ would still satisfy $h>\tau$, but would lie outside the strip, contradicting~\eqref{eq:originalstrip}. Thus $\ell_i\ge0$ throughout $P$, and the branch is the supporting face $P\cap\{\ell_i=0\}$. For $i\notin J$, the branch already satisfies $h\le\tau$.

For each $i\in J$, apply the planar algorithm~\cite{BCDJ} on $P\cap\{\ell_i=0\}$, stopping when its objective maximum is at most $\tau$ or its relaxation becomes empty. It reaches one of these outcomes because every integer point on the face satisfies $h\le\tau$. Since $\pi$ is primitive, a polynomially encoded unimodular coordinate change fixes one integer coordinate on the face. It also expresses the planar splits as split disjunctions with integer coefficients in the original three coordinates, as in the proof of Theorem~\ref{thm:gapalgorithm}. Lift the resulting cut sequences. If $J=\{0,1\}$, preserve the opposite original face in each lift.

Interleave these sequences as in Proposition~\ref{prop:interleave}. Let $B$ be the maximum of $\tau$ and the objective maxima on the nonempty current face relaxations indexed by $J$. Both branches of the original split satisfy $h\le B$, so this inequality is a split cut for the current relaxation. Add it if it is violated; it is then strong and leaves each current face relaxation unchanged. If $B>\tau$, a face attaining $B$ supplies the next strong lifted cut. If $B=\tau$, the desired bound already holds or is obtained by the cut $h\le B$. There are at most $2(N_0+N_1)+1$ cuts, where $N_i$ is the length of the sequence on branch $i$, with $N_i=0$ for $i\notin J$. The planar algorithms give polynomial bounds on these lengths and their total encoding. This proves the assertion.
\end{proof}

The strip condition is imposed on the points of the original relaxation that violate the desired bound. If it holds only after taking a split closure, the boundary hyperplanes need not support the original polytope, and the preceding lifting argument does not apply. Neither this condition nor the gap-dependent theorem gives a polynomial-time strong split-cut algorithm for every bounded three-dimensional input.

\section*{Conclusion}\phantomsection\label{sec:outlook}
\addcontentsline{toc}{section}{Conclusion}
\markboth{Conclusion}{}

We close with the questions that our results leave open. Does every rational polyhedron in $\R^3$ have a pure-integer split closure with a rational extended formulation of polynomial encoding length? Can such a formulation be constructed in polynomial time? Theorem~\ref{thm:twointeger} settles the two-integer-variable case, but efficient separation does not resolve these formulation questions: the perfect matching polytope illustrates the general distinction~\cite{R17}. The representation of pure-integer split closures through all two-dimensional lattice projections~\cite[Theorem~2.8]{BJ23} involves an infinite family and does not itself provide a polynomial formulation. Appendix~\ref{sec:beyond} gives a reduction to special bounded simplices and examines the possibilities and obstructions beyond planar normal lattices.

Theorem~\ref{thm:fixedproof} gives polynomial split proofs for rational pure-integer polyhedra in every fixed dimension. Can strong split algorithms also avoid the numerical-gap dependence in Theorem~\ref{thm:gapalgorithm}, even for bounded pure-integer inputs in dimension three? Corollary~\ref{cor:capstrip} gives such a bound when the points of the original relaxation violating the desired bound lie in a closed split strip. The question remains open for a general prescribed objective. The known $\Omega(\log G)$ lower bound is compatible with a bound polynomial in the input length, so a substantial gap remains between it and the $O_p(G+1)$ cut bound above. For arbitrary continuous objectives, the example in~\eqref{eq:ckscontinuous} obstructs even finite exact termination.

\section*{AI use disclosure}
The authors used OpenAI's ChatGPT to assist with literature searches, checking of proofs and calculations including Lean verification, and drafting parts of the exposition. The authors take full responsibility for the results, proofs, references, and final text.

\renewcommand{\bibfont}{\small}
\renewcommand{\bibsection}{\section*{References}\phantomsection\addcontentsline{toc}{section}{References}\markboth{References}{}}
\bibliographystyle{plainnat}
\bibliography{split_cuts_refs}

\clearpage
\appendix
\section{Extensions and limits of the formulation method}\label{sec:beyond}

For a pointed full-dimensional rational cone $K\subseteq\R^2$, the two-row formulation uses two properties: lattice points on the bounded edges of $\conv((K\cap\Z^2)\setminus\{0\})$ generate all integer points of $K$, and each edge gives an arithmetic progression of cuts with a small formulation. We examine how far these properties extend to higher rank. A reduction to bounded simplices is followed by a rank-three construction, an obstruction to using only edge normals, and a sufficient condition based on integral sections. These questions concern the formulation of a single split closure; the sequential proofs of Section~\ref{sec:fixedproof} do not settle them.

\subsection{Reduction to simplices}\label{sec:questions}

The tetrahedra of Section~\ref{sec:cones} have small split-closure formulations. How much harder, then, are arbitrary polyhedra? Lemma~\ref{lem:truncation} recovers the split closure of a simplicial cone from that of a bounded simplex by adding the recession directions of the cone. Extending this recovery to cones with lineality shows that, in every fixed dimension, it suffices to understand a restricted class of bounded simplices: polynomial-size extended formulations of their split closures, or polynomial-time constructions of such formulations, would give the same for every rational polyhedron.

\begin{theorem}\label{thm:simplex-reduction}
Fix $d$ and a set $I\subseteq\{1,\ldots,d\}$. The following two statements are equivalent:
\begin{enumerate}
\item Every rational polyhedron $P\subseteq\R^d$ has a split closure $\SC_I(P)$ with a rational linear extended formulation of encoding length polynomial in the input length.
\item The same assertion holds for the simplices
\[
T=a+B\{u\ge0:\textstyle\sum_i u_i\le d+1\},
\qquad a\in\Q^d,\quad B\in\Z^{d\times d}\text{ nonsingular}.
\]
\end{enumerate}
The equivalence also holds when both statements require polynomial-time construction of the formulation. In fixed dimension, the data $a,B$ and a rational inequality description of $T$ can be converted to one another in polynomial time.
\end{theorem}

\begin{proof}
First consider the conversion between representations. From $a,B$, invert $B$ to obtain the inequalities $B^{-1}(x-a)\ge0$ and $\mathbf1^TB^{-1}(x-a)\le d+1$. Conversely, enumerate the $d+1$ vertices of a simplex in the stated class. Choose one as $a$ and divide its differences to the other vertices by $d+1$ to obtain the columns of $B$. These columns are integral and independent. In fixed dimension, vertex enumeration and these rational operations have polynomial encoding length and running time.

Only the implication from (2) to (1) needs proof. Apply identity~\eqref{eq:corneridentity}. In fixed dimension it suffices to treat polynomially many cones $C=v+R\R_+^r+L$ of the form~\eqref{eq:generalcone}, with polynomially encoded rational data. Empty and rank-zero inputs are immediate. Scale each column of $R$ by a positive integer to make it integer without changing $C$. Choose linearly independent integer vectors $t_1,\ldots,t_{d-r}$ spanning $L$. Both operations have polynomial encoding cost in fixed dimension. Let
\[
B=[R\ t_1\ \cdots\ t_{d-r}],\qquad
a=v-\sum_{j=1}^{d-r}t_j,
\]
and define $T$ as in (2). We claim that, with $H=R\R_+^r+L$,
\begin{equation}\label{eq:lineality-recovery}
\SC_I(C)=\SC_I(T)+H.
\end{equation}
The case $r=d$ is~\eqref{eq:cone-recovery}. The shift from $v$ to $a$ is an integer vector, and it assigns coordinate one to each added direction $t_j$, which will make every additional cut hold at the point constructed below. The set of coordinates designated integer is unchanged.

Take $x=v+Rs+\ell\in\SC_I(C)$, set $s'_j=\min\{s_j,1\}$ for each $j$, and let $y=v+Rs'$. In the pointed cone $D=a+B\R_+^d$, the point $y$ has coordinates $(s',1,\ldots,1)$. A permitted normal orthogonal to $L$ has zero products with the added rays $t_j$. Its fractional part at $a$ equals its fractional part at $v$, because $v-a$ is integer. Thus its cut in $D$ is exactly the cut in the ray coordinates of $C$. This replacement preserves that cut by the proof of Lemma~\ref{lem:truncation}. A normal not orthogonal to $L$ has a nonzero integer product with some $t_j$. If its fractional part is nonzero, its cut coefficient on that ray is at least one, and the corresponding coordinate of $y$ is one. Thus $y$ satisfies every cut of $D$. Its coordinates sum to at most $d$, so Lemma~\ref{lem:truncation} gives $y\in\SC_I(T)$. Since $x-y\in H$, one inclusion in~\eqref{eq:lineality-recovery} follows. Conversely, every point of $T$ is $v+Rs+\sum_j(\eta_j-1)t_j$ for some $s,\eta\ge0$, so $T\subseteq C$. By monotonicity, $\SC_I(T)\subseteq\SC_I(C)$. The only splits that can change $C$ have normals orthogonal to $L$; their cut coefficients in $s$ are nonnegative. Adding $H=R\R_+^r+L$ therefore preserves all these cuts. This proves the identity.

A formulation of $\SC_I(T)$ now yields one of $\SC_I(C)$ by adding $r$ nonnegative ray variables and $d-r$ free lineality variables. All coefficients have polynomial encoding length. Intersect these formulations using separate auxiliary variables for the relaxations in~\eqref{eq:corneridentity}. The number and total encoding of these systems are polynomial, proving both the existence and the constructive assertions. No change of the original integrality lattice is used.
\end{proof}

In dimension three, the special simplices are tetrahedra whose vertex differences belong to $4\Z^3$, although their vertices themselves may be fractional. Thus even this restricted tetrahedral case would resolve the general three-variable formulation question.

\subsection{A rank-three family}\label{sec:threefamily}

The reduction leaves the problem of representing the split cuts of a simplicial cone. Use the cone $C$, the space $W$, and the ray vectors $r_j$ from~\eqref{eq:generalcone}. Consider a progression of sign-compatible normals whose step has scalar products with the rays all nonnegative or all nonpositive. Imposing the cut of that step makes the quadratic expression in~\eqref{eq:clearedsplit} concave on each interval where the fractional part is affine. The two endpoint cuts then imply all intermediate cuts. This observation strengthens Lemma~\ref{lem:progression} in any rank and yields the rank-three construction below.

\begin{lemma}\label{lem:onesign}
Let $a,u\in W\cap\Z^d$, and let $M\ge0$ be an integer. Suppose $a+tu$, $0\le t\le M$, lies in one cone $K_\sigma$ from~\eqref{eq:signcone}, and all products $u^Tr_j$ are nonnegative or all are nonpositive. On $s\ge0$, the cuts with normals $a$, $a+Mu$, and $u$ imply every cut with normal $a+tu$, $t=0,\ldots,M$.
\end{lemma}

\begin{proof}
The cases $M=0$ and $u=0$ are immediate. Reverse the progression if necessary, using $\alpha(-u)=\alpha(u)$, so that $u^Tr_j=-d_j$ with $d_j\ge0$. Set $g=\{u^Tv\}$ and $D(s)=\sum_jd_js_j$. If $g=0$, the fractional parts are constant and the endpoint argument in Lemma~\ref{lem:progression} applies.

Suppose $0<g<1$. The cut with normal $u$ is $D(s)\ge g$. Fix $s\ge0$ satisfying this cut and the nontrivial endpoint cuts. For real $t\in[0,M]$, define $f(t)=\{a^Tv+gt\}$. This agrees with $\{(a+tu)^Tv\}$ at integer $t$; its values at nonintegers are only an auxiliary interpolation. Use the fixed ray signs to define $U(t,s),V(t,s)$ as in
\eqref{eq:clearedsplit}. On an interval between consecutive
integer crossings of $a^Tv+gt$, the left side of that inequality has quadratic coefficient
\[
g\bigl(g-D(s)\bigr)\le0.
\]
It is therefore concave in $t$ on that interval.

On an interval between two successive crossings, the expression tends to $V\ge0$ at the left endpoint, since $f$ tends to zero. At its right endpoint the limit is $U\ge0$, since $f$ tends to one. For the truncated first or last interval, the endpoint cut gives the same nonnegative boundary value when the fraction is nonzero. If it is zero, the one-sided limit is again $U$ or $V$. Concavity gives nonnegativity throughout every interval. At every integer index with $0<f(t)<1$, divide by $f(t)(1-f(t))$ to obtain its split cut. Indices with $f(t)=0$ require no cut.
\end{proof}

The family below is governed by the congruence lattice $w_1+w_2+w_3\equiv0\pmod N$. In its nonnegative orthant, the points of coordinate sum $N$ generate all others. In the other orthants, we first subtract integer vectors with one positive and one negative coordinate. The remaining vector has coordinates of one sign. Applying Lemma~\ref{lem:onesign} twice will then show that the cuts from a triangular set of normals follow from those on one edge and four additional normals.

For $N\ge1$ and $v\in\Q^3$, define
\begin{equation}\label{eq:rankthreefamily}
\widehat R_N=\begin{pmatrix}N&0&0\\-1&1&0\\-1&0&-1\end{pmatrix},
\qquad \widehat C_N(v)=v+\widehat R_N\R_+^3.
\end{equation}
All three coordinates are integer variables. The determinant is $-N$, so the relevant normal lattice has rank three. Write $\rho(a)=\lceil a\rceil-a$, which is zero at integers.

The next proposition gives both a logarithmic-size extended formulation and a polynomial-size description in the original space. The latter requires a separate argument, and its size depends on the full input, including $v$.

\begin{proposition}\label{prop:threefamily}
In the coordinates $x=v+\widehat R_Ns$, the full split closure of $\widehat C_N(v)$ is given by $s\ge0$ and
\begin{align}
(N-t)s_1+ts_2&\ge\rho(v_1+tv_2)
&& (t=0,\ldots,N),\label{eq:rankthreebottom}\\
s_1+s_3&\ge\{v_3\},\label{eq:rankthreeone}\\
s_2+s_3&\ge\rho(v_2-v_3),\label{eq:rankthreetwo}\\
Ns_3&\ge\{v_1+Nv_3\},\label{eq:rankthreethree}
\end{align}
together with
\begin{equation}\label{eq:rankthreestep}
\frac{s_1}{f_2}+\frac{s_2}{1-f_2}\ge1
\quad\text{if } f_2=\{v_2\}>0.
\end{equation}
It admits a rational linear extended formulation with $O(\log(N+2))$ inequalities and auxiliary variables, of total encoding length and construction time polynomial in the encoding of $N,v$. It also has a polynomial-encoding inequality description in the original space, computable in polynomial time.
\end{proposition}

\begin{proof}
For $\pi=(p,q,r)$, the scalar products with the rays are
\[
\widehat R_N^T\pi=(Np-q-r,q,-r).
\]
The displayed inequalities are precisely the cuts of the normals $(1,t,0)$, $e_3$, $e_2-e_3$, $(1,0,N)$, and $e_2$, with inequalities already implied by $s\ge0$ when the fractional part is zero. We prove that these suffice for every integer normal.

Introduce $w=(Np-q-r,q,r)$. This is a bijection from $\Z^3$ onto
\[
\Lambda_N=\{w\in\Z^3:w_1+w_2+w_3\equiv0\pmod N\},
\qquad
\pi=\bigl((w_1+w_2+w_3)/N,w_2,w_3\bigr).
\]
The inverse image of each closed orthant in the $w$ coordinates is a cone of the form~\eqref{eq:signcone}, because the scalar products with the rays are $(w_1,w_2,-w_3)$. In the nonnegative orthant, every nonzero point has coordinate sum $kN$ for an integer $k\ge1$. Distributing its coordinate units into $k$ groups of $N$ proves generation by the points with sum $N$. Their normals are
\[
\mathcal H_N=\{(1,q,r):q,r\in\Z_+,\ q+r\le N\}.
\]
The nonpositive orthant gives the same cuts by negation.

In an orthant with two nonnegative coordinates $i,j$ and one nonpositive coordinate $k$, cancel opposite units using $e_i-e_k$ and $e_j-e_k$. Each cancellation vector belongs to $\Lambda_N$ and stays in the same orthant. Once cancellation stops, the remaining vector has only one sign and coordinate sum a multiple of $N$. Its generators are $te_i+(N-t)e_j$ for $0\le t\le N$, or $-Ne_k$. Orthants with two negative coordinates follow by negation. Under the inverse map, the generators whose coordinates are all nonnegative or all nonpositive lie on the boundary of $\mathcal H_N$, up to sign. The cancellation vectors give, up to sign, just $e_2,e_3,e_2-e_3$. Lemma~\ref{lem:addition} therefore shows that the cuts of $\mathcal H_N$ and these three normals imply every split cut.

For fixed $q$, the progression $(1,q,0)+re_3$, $0\le r\le N-q$, lies in $K_{(+,+,-)}$. Its step has products $(-1,0,-1)$. By Lemma~\ref{lem:onesign}, the cut of $e_3$ reduces this progression to its endpoints. One set of endpoints is $(1,q,0)$, retained in~\eqref{eq:rankthreebottom}; the other is $(1,q,N-q)$, $0\le q\le N$. The latter progression has step $e_2-e_3$ with products $(0,1,1)$. A second application of the lemma reduces it to $(1,0,N)$ and $(1,N,0)$ in the presence of
\eqref{eq:rankthreetwo}. The remaining cancellation normal $e_2$
gives~\eqref{eq:rankthreestep}. This proves the full description.

For the extended formulation, apply Lemma~\ref{lem:progression} to $(1,t,0)$, $0\le t\le N$, whose cut with normal $e_2$ is already
\eqref{eq:rankthreestep}. Add the three other displayed cuts,
nonnegativity, and $x=v+\widehat R_Ns$. This proves the logarithmic size and polynomial encoding and construction bounds, with no restriction on the denominator of $v$.

For a direct inequality description, choose a common positive denominator $D$ of $v_1,v_2$ and form the bounded rational polygon
\[
Q=\{(t,h)\in\R^2:0\le t\le N,\quad
0\le D(h-v_1-tv_2)\le D-1\}.
\]
Its integer points are exactly $(t,\lceil v_1+tv_2\rceil)$ for integer $0\le t\le N$. For such a point, \eqref{eq:rankthreebottom} is
\[
Ns_1+v_1+t(s_2-s_1+v_2)-h\ge0.
\]
The expression is affine in $(t,h)$, so it is enough to retain the vertices of $\conv(Q\cap\Z^2)$. The fixed-dimensional integer-hull results in~\cite{CHKM,H} compute this polynomially sized list in polynomial time. Add \eqref{eq:rankthreeone}--\eqref{eq:rankthreestep} and substitute $s=\widehat R_N^{-1}(x-v)$. The resulting system has polynomial total encoding length in $N,v$.
\end{proof}

\subsection{Interior normals}\label{sec:obstruction}

For the preceding family, the cuts indexed by a triangular set of generators follow from those on one edge and a few additional cuts. For a three-dimensional cone $K_\sigma$ from~\eqref{eq:signcone}, write $K=K_\sigma$. Do the lattice points on the edges of $\conv((K\cap\Z^3)\setminus\{0\})$ suffice, together with normals whose scalar products with the rays are all nonnegative or all nonpositive? The next example shows that they do not. A lattice point in the relative interior of a bounded face gives a cut that is not implied even by all the other split cuts, except the identical cut from its negative. This rules out that choice of normals, but gives no lower bound on extension complexity.

\begin{proposition}\label{prop:interior-normal}
There are a translated simplicial cone $C\subseteq\R^3$ and an integer normal $\pi_*$ such that the intersection cut~\eqref{eq:cut} of $\pi_*$ is not implied by the intersection cuts of any normals other than $\pi_*$ and $-\pi_*$. The normal $\pi_*$ lies in the relative interior of a two-dimensional bounded face of $\conv((K\cap\Z^3)\setminus\{0\})$ where $K=\{\pi:(\pi^Tr_j)(\pi_*^Tr_j)\ge0\text{ for all }j\}$. All three scalar products with the rays are nonzero and have mixed signs.
\end{proposition}

\begin{proof}
Take
\[
 R=\begin{pmatrix}-1&2&1\\-1&-1&-2\\-1&-2&1\end{pmatrix},
 \qquad v=(1/4,1/2,1/2),\qquad C=v+R\R_+^3.
\]
For $\pi=(x,y,z)\in\Z^3$, its scalar products with the rays and its fractional part are
\[
 t_1=-x-y-z,\quad t_2=2x-y-2z,\quad t_3=x-2y+z,
 \qquad f=\{(x+2y+2z)/4\}.
\]
The normal $\pi_*=(0,0,1)$ has products $(-1,-2,1)$ and fraction $1/2$, giving
\begin{equation}\label{eq:interior-cut}
 2s_1+4s_2+2s_3\ge1.
\end{equation}

For a normal in $K$, $w=(-t_1,-t_2,t_3)\ge0$ and $\sum w_i=4z$. Since $\det R=12\ne0$, a nonzero normal gives $w\ne0$, hence $z>0$. For a lattice point, this implies $z\ge1$. On $z=1$, the inequalities are $y\ge-1-x$, $y\ge2x-2$, and $y\le(x+1)/2$. They force $-1\le x\le1$, and their lattice points are exactly
\[
 A=(-1,0,1),\ B=(0,-1,1),\ \pi_*,\ C_0=(1,0,1),\ E=(1,1,1).
\]
Every nonzero integer point of $K$ has $z\ge1$, so $z=1$ is a supporting plane of $\conv((K\cap\Z^3)\setminus\{0\})$. Its section is exactly the convex hull of the five listed points: a convex combination of points with $z\ge1$ can have $z=1$ only if every point used has $z=1$. The points $A,B,C_0,E$ are the vertices of a bounded two-dimensional polygon. The identity $\pi_*=(2A+B+C_0+E)/5$ expresses it as a convex combination with positive weight on each vertex, placing it in the relative interior.

Now take $s_*=(11,3,5)/48$. The left-hand side of the target inequality equals $11/12<1$ at this point. To prove that the target cut is not implied by the others, it remains to show that $s_*$ satisfies every other cut written with right-hand side one. We do this by bounding its value for an arbitrary integer normal. Set $W=11|t_1|+3|t_2|+5|t_3|$. The identities
\[
 12x=-5t_1+3t_2+t_3,\qquad x+2y+2z=-2t_1-x
\]
give $W\ge12|x|$. For even $x$, a nontrivial cut has $f=1/2$ and value $W/24$ at $s_*$. Thus the cut is satisfied whenever $|x|\ge2$. For odd $x$, $f\in\{1/4,3/4\}$ and the value is at least $W/36$, so the cut is satisfied whenever $|x|\ge3$. Negation leaves a cut unchanged; only $x=0,1$ remain.

If $x=0$, let $k=t_1$, $l=t_3$. Then $3t_2=5k-l$, $l\equiv2k\pmod3$, and a nontrivial cut requires odd $k$. By negation take $k>0$. For $k\ge3$, $W\ge33>24$. For $k=1$,
\[
 W=11+|5-l|+5|l|,\qquad l\equiv2\pmod3.
\]
As a function of $l$, the sum of the last two terms decreases up to $l=0$ and increases thereafter. The neighboring permitted values $l=-1,2$ give $W=22,24$, respectively. The choice $k=1,l=-1$ gives exactly $-\pi_*$; every other permitted value gives $W\ge24$.

If $x=1$, again put $k=t_1$, $l=t_3$, so $3t_2=12+5k-l$ and $l\equiv2k\pmod3$. The fractional part depends on the parity of $k$. Define
\[
 \Phi_0(t)=\max(3t,-t),\qquad \Phi_1(t)=\max(t,-3t),
\]
and fix $\epsilon\in\{0,1\}$ with $k\equiv\epsilon\pmod2$. Write $\Phi=\Phi_\epsilon$ for this parity. The left-hand side of the cut inequality at $s_*$ is $N/36$, where subadditivity gives
\[
 N=11\Phi(k)+\Phi(12+5k-l)+5\Phi(l)
 \ge11\Phi(k)+\Phi(12+5k)+4\Phi(l).
\]
For $\epsilon\in\{0,1\}$, define $F_\epsilon(k)=11\Phi_\epsilon(k)+\Phi_\epsilon(12+5k)$. Each function is convex on the real line. We evaluate $F_0$ at even integers and $F_1$ at odd integers; three values of each suffice:
\[
\begin{array}{c|ccc}
\text{even indices}&F_0(-4)=52&F_0(-2)=28&F_0(0)=36\\
\text{odd indices}&F_1(-1)=40&F_1(1)=28&F_1(3)=60.
\end{array}
\]
Along either arithmetic progression, convexity makes successive differences nondecreasing. The middle value is smaller than both neighbors, so every value outside these triples is at least the nearer endpoint value. Thus $F_\epsilon(k)\ge36$ whenever $k\equiv\epsilon\pmod2$, except possibly at $k=-2,1$.

Only $k=-2,1$ remain. For $k=1$, the congruence gives $\Phi_1(l)\ge2$, hence $N\ge11+17+8=36$. For $k=-2$, it gives either $l\le-1$, when $N=28-8l\ge36$, or $l\ge2$, when $N\ge22+30>36$. Thus every cut other than those of $\pm\pi_*$ is satisfied at $s_*$. The products $\pi_*^Tr_j$ are all nonzero and include both signs, so $\pi_*$ also lies in the interior of $K$, completing the proof.
\end{proof}

\FloatBarrier
\subsection{Integral sections}\label{sec:fibers}

The preceding example shows that edge normals alone need not suffice. We now give a condition under which the inequalities indexed by the integer points of a rational polytope have a compact extended formulation. For $r\ge1$, consider a rational polytope
\[
T=\{(t,u)\in\R\times\R^r:Bu\le b-at\}.
\]
For each real $t$, let $T_t=\{u:Bu\le b-at\}$. Thus $T_t$ is the projection onto the $u$-coordinates of the intersection of $T$ with the hyperplane fixing its first coordinate at $t$.

Suppose that $T_t$ is integral for every integer $t$ for which it is nonempty. An inequality that is affine in $u$ holds at all integer points of $T_t$ if and only if it holds at all vertices of $T_t$. We use this observation for inequalities that are quadratic in $t$ and affine in $u$. Solving a fixed set of $r$ linearly independent equations from $Bu=b-at$ expresses $u$ as an affine function of $t$. Substitution then gives a quadratic expression in $t$, to which Lemma~\ref{lem:interval} applies. The following proposition makes this argument precise and gives a polynomial-time test of the integrality condition.

\begin{proposition}\label{prop:integralfibers}
Fix $r\ge1$, and let $T$ and $T_t$ be as above, with $T$ a rational polytope. Suppose that every nonempty $T_t$ with $t\in\Z$ is integral, that is, the convex hull of its integer points. Let
\begin{equation}\label{eq:fiberpoly}
\Phi_x(t,u)=A(x)t^2+B_0(x)t+C(x)+\sum_{i=1}^r\bigl(E_i(x)t+G_i(x)\bigr)u_i,
\end{equation}
where all displayed coefficients are rational affine functions of real variables $x$. Then
\[
\{x:\Phi_x(t,u)\ge0\text{ for every }(t,u)\in T\cap\Z^{r+1}\}
\]
has a rational linear extended formulation of polynomial encoding length, constructible in polynomial time. The integrality hypothesis is also testable in polynomial time.
\end{proposition}

\begin{proof}
If $T$ is empty there are no indexed inequalities. Otherwise boundedness implies $\operatorname{rank}(B)=r$. Enumerate its nonsingular $r$-row submatrices $B_J$. Their number is polynomial because $r$ is fixed. Solving the corresponding equations gives a candidate vertex of $T_t$ of the form
\[
u_J(t)=B_J^{-1}b_J-B_J^{-1}a_Jt=p_J+q_Jt.
\]
The real values of $t$ for which this point satisfies every inequality in $Bu\le b-at$ form a bounded rational interval $I_J$, possibly empty or a singleton. Discard empty $I_J$. For each remaining interval, put $l_J=\lceil\inf I_J\rceil$ and $h_J=\lfloor\sup I_J\rfloor$, and discard it if $l_J>h_J$. All these data have polynomial encoding. Feasibility of the candidate is an intersection of linear inequalities in the single variable $t$, so $I_J$ is indeed an interval. If a nonempty $I_J$ were unbounded, the points $(t,p_J+q_Jt)$ would contradict boundedness of $T$.

Every vertex of every nonempty $T_t$ appears in this enumeration, even when $T_t$ is lower-dimensional. Conversely, each feasible candidate is a vertex, since it satisfies $r$ linearly independent constraints at equality. 

To test the hypothesis, use the same intervals. If $l_J=h_J$, check that $p_J+q_Jl_J$ is integral. If $l_J<h_J$, check that both $p_J+q_Jl_J$ and $q_J$ are integral. These conditions are necessary, by subtracting the values at two consecutive integer indices, and sufficient by the affine formula. They hold exactly when every $T_t$ with $t\in\Z$ is integral. No enumeration of the integer indices is needed.

For integer $t$, the hypothesis makes every vertex of $T_t$ integral. Since~\eqref{eq:fiberpoly} is affine in $u$ for fixed $t$, it is nonnegative at every integer point of $T_t$ if and only if it is nonnegative at every vertex of $T_t$. The required inequalities are therefore equivalent to
\[
\Phi_x(t,p_J+q_Jt)\ge0\qquad(t=l_J,\ldots,h_J)
\]
for every retained $J$. Each expression is quadratic in $t$, with coefficients affine in $x$. Lemma~\ref{lem:interval} gives a rational extension with $O(\log(h_J-l_J+2))$ inequalities and variables. Concatenating these formulations proves both exactness and the polynomial construction bound.
\end{proof}

For example, the index set $0\le t\le N$, $u\ge0$, $\sum_i u_i\le t$ gives $T_t=\conv\{0,te_1,\ldots,te_r\}$ for each integer $t$ with $0\le t\le N$. Even with arbitrary mixed terms $tu_i$, Proposition~\ref{prop:integralfibers} reduces its inequality family to $r+1$ univariate quadratic families, giving $O_r(\log(N+2))$ inequalities and auxiliary variables.

For split cuts, write $f=\pi^Tv-k$ and use the notation of~\eqref{eq:iqpobjective}. Since $u(\pi)-w(\pi)=\pi^T(\bar x-v)$, that objective factors as
\[
F(\pi,k)=f(\pi^T\bar x-k-1)+w(\pi).
\]
The factor $\pi^Tv-k$ suggests a choice of the first index. Clear its denominators and divide its integer coefficient vector by its gcd. A polynomially encoded unimodular change makes the resulting primitive linear form the first index $t$~\cite{KB}. Thus $f=\eta t$ for a nonzero rational constant $\eta$. Both $\pi^T\bar x-k-1$ and $w(\pi)$ are affine in the transformed indices, with coefficients affine in the query point $\bar x$. Their substitution gives precisely~\eqref{eq:fiberpoly}, with $x$ there denoting $\bar x$ here. The proposition therefore gives a compact extended formulation for this cut family, provided that, for the transformed polytope $T$, every nonempty $T_t$ with $t\in\Z$ is integral. This condition must be checked using the test above; it need not hold for a general three-dimensional domain.

\end{document}